\documentclass[11pt,x11names]{article}
\usepackage{amsthm}

\PassOptionsToPackage{square,sort,comma,numbers}{natbib}

\usepackage[utf8]{inputenc} 
\usepackage[T1]{fontenc}    
\usepackage{hyperref}       
\usepackage{url}            
\usepackage{booktabs}       
\usepackage{amsfonts}       
\usepackage{nicefrac}       
\usepackage[T1]{fontenc}
\usepackage{lmodern}  
\usepackage{microtype}
\usepackage{framed}

\usepackage{import}
\usepackage{latexsym,amssymb,amsmath,amsthm,amsfonts}
\usepackage[dvipsnames]{xcolor}
\usepackage{graphicx,shadow}
\usepackage{float}
\usepackage[square,sort,comma,numbers]{natbib}
\usepackage[a4paper,
bindingoffset=0.2in,%
left=2.5cm,right=2.5cm,top=2.5cm,bottom=2.5cm]{geometry}
\usepackage{authblk}
\usepackage{mathrsfs}
\usepackage{caption,subcaption}
\usepackage{soul}
\usepackage{colortbl} 
\usepackage{multirow} 
\usepackage{epstopdf}
\usepackage[boxed,ruled,linesnumbered]{algorithm2e}
\usepackage{enumitem}
\usepackage[noabbrev,capitalise,nameinlink]{cleveref}
\usepackage{glossaries}

\hypersetup{colorlinks={true},linkcolor={MidnightBlue},citecolor=ForestGreen}
\usepackage{comment}
\usepackage{ulem}
\usepackage{svg}
\usepackage{booktabs}

\def\red{\color{red}}

\newcommand{\nh}[1]{}
\newcommand{\dg}[1]{}

\numberwithin{equation}{section}

\def\argmin{\displaystyle\mathop {\mbox{argmin}}}
\newcommand{\pos}[1]{\mathrm{pos}\left(#1\right)}

\def\argmax{\displaystyle\mathop {\mbox{argmax}}}

\newtheorem{assumption}{Assumption}

\newtheorem{theorem}{Theorem}[section]
\newtheorem{lem}{Lemma}[section]
\newtheorem{corollary}{Corollary}[section]
\newtheorem{proposition}{Proposition}[section]

\theoremstyle{definition}
\newtheorem{example}{Example}[section]
\newtheorem{defin}{Definition}[section]
\newtheorem{remark}{Remark}[section]

\newcommand{\aff}[1]{\displaystyle\mathop{ \rm{aff}}\left(#1\right)}  

\newcommand{\ri}[1]{\displaystyle\mathop{ \rm{ri}\left(#1\right)}}  
\newcommand{\rb}[1]{\displaystyle\mathop{ \rm{rb}\left(#1\right)}}  
\newcommand{\unif}[1]{\mathop{Unif \left(#1\right)}}

\newcommand{\cl}[1]{\displaystyle\mathop{ \rm{cl}\left(#1\right)}}  

\def\conv{\displaystyle\mathop {\mbox{\rm conv}}}  

\def\ext{\displaystyle\mathop {\mbox{\rm ext}}}    

\def\dist{\displaystyle\mathop {\mbox{\rm dist}}}  

\def\ker{\displaystyle\mathop {\mbox{\rm Ker}}}    

\def\sup{\displaystyle\mathop {\mbox{\rm sup}}}

\def\inf{\displaystyle\mathop {\mbox{\rm inf}}}

\def\argmin{\displaystyle\mathop {\mbox{\rm argmin}}}

\def\argmax{\displaystyle\mathop {\mbox{\rm argmax}}}

\def\real{\mathbb R}

\newcommand{\N}{{\mathcal{N}}}

\newcommand{\sign}{\mbox{sign}}

\newcommand{\inner}[2]{{\langle#1,#2\rangle}}
\newcommand{\norm}[1]{\left\|#1\right\|}

\DeclareMathOperator{\diam}{diam}

\newcommand{\interior}[1]{\textbf{int}(#1)}
\newcommand{\bdr}[1]{\textbf{bd}(#1)}
\newcommand{\tanCone}[1]{\mathcal{T}_{\mathcal{C}}\left(#1\right)}

\newcommand{\finitePow}[1]{\mathcal{P}_{fin}\left(#1\right)}

\usepackage{bm}

\newcommand{\Comment}[1]{}

\title{Randomized and Oracle-Based Feasible Directions Methods for Nonsmooth Nonconvex Constrained Optimization }

\author{%
    Dan Greenstein\textsuperscript{1} \quad Nadav Hallak\textsuperscript{1}\thanks{Corresponding Author} \\
    \textsuperscript{1} Technion Israel Institute of Technology, Faculty of Data and Decision Sciences, \\
    {\tt\small  sdngreen@campus.technion.ac.il, ndvhllk@technion.ac.il} \\
}

\begin{document}
\date{ }
\maketitle
\begin{abstract}

Minimizing nonsmooth nonconvex objectives over general convex sets using no derivatives of the objective, together with structural access to the constraint set, poses a fundamental geometric difficulty: at boundary points, the feasible direction cone may fail to be closed, preventing the construction of finite sets of search directions that capture all feasible descent directions.
This paper develops a feasible-directions framework for problems structured as the difference between a continuously differentiable function and a convex function over a closed convex set, and accesses the objective without knowledge of its decomposition, through exact line searches along feasible segments.
To address the geometric obstruction, we introduce $\eta$-feasible directions, which remain uniformly separated from the boundary of the feasible direction cone, thereby enabling finite-direction search even in the presence of nonlinear constraints.
We present two algorithmic realizations of this framework. The first assumes access to an $\eta$-direction oracle and extends feasible-directions methods from linear to general convex constraints. The second replaces the oracle by a noise-augmented randomization scheme that constructs spanning directions from the affine-hull inner normal and almost surely avoids nonsmooth boundary points when initialized outside that set.
Under curvature bounds and Lipschitz continuity of the smooth component's gradient, we jointly control the spanning coefficients and feasible step lengths of the oracle construction. This yields finite-time approximate $d$-stationarity along every unit feasible direction admitting a prescribed positive step length.
For both methods, we also establish a rate in the proximal directional minimization (PDM) measure, assuming uniformly bounded direction-set sizes.
Randomized exploration ensures that, almost surely, every accumulation point is $d$-stationary.

\end{abstract}

\section{Introduction}
\subsection{Problem Formulation}

This paper studies the nonsmooth nonconvex constrained problem 
\begin{equation}\label{eq:feasible_directions_beck_hallak}\tag{P}
    \min\limits_{x \in \mathcal{C}} h(x) := f(x) - g(x),
\end{equation}
where $f: \real^n \rightarrow \real$ is continuously differentiable, $g: \real^n \rightarrow \real$ is convex and possibly non-differentiable,  the feasible set $\mathcal{C} \subseteq \real^n$ is closed and convex, and $h$ is lower bounded over $\mathcal{C}$. 
We further assume that $|\mathcal{C}|>1$ and $n\geq 2$ to exclude  trivial or pathological cases. These conditions imply that the dimension of $\mathcal{C}$ is at least one, i.e., $\dim\mathcal{C}\ge 1$.

Problem \eqref{eq:feasible_directions_beck_hallak} captures a broad class of nonlinear constrained models, including the prototypical problems below.

\begin{example}[Difference-of-convex programs]
    When $f$ is convex, problem \eqref{eq:feasible_directions_beck_hallak} reduces to a difference-of-convex (DC) program. 
    DC formulations arise in a wide range of applications; see, e.g., \cite{Ahn2017} studying DC programs machine learning tasks.
\end{example}

\begin{example}[Linear constraints]
    When the feasible set $\mathcal{C}$ is a polytope, problem \eqref{eq:feasible_directions_beck_hallak} reduces to the model studied in \cite{feasible_directions_2020,BECK2022517} -- see application examples therein.
\end{example}

\begin{example}[Concave regularization]
    Problem \eqref{eq:feasible_directions_beck_hallak} includes sparsity-promoting regularizers $\lambda\sum_i\varphi(|x_i|)$ with $\varphi$ concave and non-decreasing on $[0,\infty)$ and $\varphi'(0^+)<\infty$, such as the log-sum penalty \citep{candes2008enhancing} and SCAD \citep{fan2001variable}.
    Adding variables $y$ with the linear constraints $y_i\ge\pm x_i$ and replacing $\varphi(|x_i|)$ by $\varphi(y_i)$ gives an equivalent instance of \eqref{eq:feasible_directions_beck_hallak}; the condition $\varphi'(0^+)<\infty$ lets $-\varphi$ extend to a finite convex function on $\real$.
\end{example}

\Comment{\begin{example}[Conic and Semidefinite Programming]
    Problem~\eqref{eq:feasible_directions_beck_hallak} also captures conic optimization models with closed cones, including semidefinite programming (SDP). Conic programming refers to a class of convex problems in which the feasible region is defined by a convex cone~\cite{ben2001lectures,boyd2004convex}, and includes SDPs as a prominent subclass, where the decision variable is a symmetric matrix constrained to be positive semidefinite~\cite{vandenberghe1996semidefinite}. Such problems arise in a wide range of applications, including control theory, combinatorial optimization, and quantum information~\cite{goemans1995improved,parrilo2003semidefinite,boyd1994linear,skrzypczyk2023semidefinite}.
\end{example}
\begin{example}[Conic and semidefinite constraints]
Problem~\eqref{eq:feasible_directions_beck_hallak} also accommodates optimization
over closed convex cones, including second-order and semidefinite cones.
In such instances, the cone constraint is absorbed into the feasible set $\mathcal{C}$,
and our framework applies whenever feasible directions in $\mathcal{C}$ can be generated
efficiently. This includes several structured conic and SDP models for which
the tangent and normal cones admit tractable representations.
\end{example}
}

In nonsmooth nonconvex optimization, a critical point may still admit
a feasible direction of first-order descent. Directional stationarity
provides a stronger target: it requires a nonnegative directional
derivative along every feasible direction and thus rules out all
feasible first-order improvements. Our goal is to obtain guarantees
of this kind while maintaining feasibility, without requiring
objective derivatives or access to the decomposition into $f$ and $g$.
We study methods that access the objective through exact regularized
line searches and exploit the geometry of the constraint set.

The central question is whether finitely many feasible line searches
can provide directional optimality guarantees beyond the directions
actually searched. For the objective structure considered here,
non-descent along the search directions extends to their positive
span. The geometry of this span therefore determines the scope of
the resulting optimality guarantees.

When $\mathcal{C}$ is polyhedral, the feasible-direction cone admits
a finite positive spanning set. This property underlies the
finite-time guarantees in the proximal directional minimization
(PDM) measure established in \cite{BECK2022517}, without requiring
objective derivatives. For general closed convex sets, however,
the feasible-direction cone may be nonclosed at boundary points
and therefore admit no finite positive spanning set. This obstruction
arises even for basic constraint sets such as Euclidean balls and
ellipsoids: finitely many feasible search directions cannot cover
all feasible directions through their positive span.

To retain finite-direction guarantees under general convex constraints,
we introduce $\eta$-feasible directions, which remain angularly
separated from the boundary of the feasible-direction cone.
This restriction allows finitely many feasible search directions
to cover the $\eta$-feasible directions through their positive span.
It provides the geometric basis for finite-time PDM guarantees
with full directional coverage in the relative interior and
a controlled restriction at nonlinear boundary points.

We implement this idea through two feasible-directions methods that differ in how they construct the search directions: one uses a constraint-set oracle, while the other uses a randomized construction based on the local boundary geometry. At each iteration, the method performs exact regularized line searches along the resulting directions and selects the candidate with the greatest objective decrease. The searches are restricted to feasible segments, so all iterates remain in $\mathcal{C}$. Neither objective derivatives nor access to the decomposition into $f$ and $g$ is required (\cref{rem:line-search-access}).

Directional coverage alone does not turn small objective decrease
into a quantitative bound on directional derivatives: the spanning
coefficients and feasible step lengths must also be controlled.
Under curvature bounds and Lipschitz continuity of $\nabla f$,
our oracle construction provides this control, yielding finite-time
approximate $d$-stationarity along every unit feasible direction
admitting a prescribed positive step length
(\cref{thm:finite-time-rate,cor:eps-complexity}).

To obtain full $d$-stationarity despite the restricted directional coverage at boundary points, we supplement these searches with the randomized exploration step of \citet{greenstein2026achieving}. 
At each iterate $x$, we sample a direction uniformly from the unit
sphere in $\aff{\mathcal{C}-x}$ and search along its feasible segment. 
This exploration ensures that, almost surely, every accumulation point of either method is $d$-stationary (\cref{thm:main}).

\paragraph{Outline.}
Following the preliminaries in \Cref{sec:preliminaries}, the main development is concentrated in Sections~\ref{sec:eta_feasible}--\ref{sec:analysis}: Section~\ref{sec:eta_feasible} introduces $\eta$-feasible directions; Sections~\ref{sec:directions_oracle_and_alg} and~\ref{sec:noised_alg} present the oracle-based and noised methods, respectively; and \Cref{sec:analysis} establishes PDM rates, almost-sure $d$-stationarity of accumulation points, and finite-time guarantees under curvature assumptions. \Cref{sec:noised_alg_properties} supplies the supporting spanning and boundary-avoidance results. The extensive appendices collect auxiliary proofs, detailed convergence arguments, and measure-theoretic foundations needed for completeness.

\subsection{Contribution}

\paragraph{Geometric framework.}
In \Cref{sec:eta_feasible} we develop the notion of $\eta$-feasible directions ($\eta$-FDs) and establish its properties. At smooth boundary points of nonlinear constraint sets, the $\eta$-FD cone is contained in the positive span of finitely many feasible directions, yielding a PDM guarantee over that cone.

\paragraph{Feasible-directions framework.}
We propose two methods utilizing $\eta$-FDs:
(i) an oracle-based method that extends the linearly constrained feasible-directions method of \citet{BECK2022517} to general convex constraint sets (\Cref{sec:directions_oracle_and_alg}), and
(ii) a randomized method that generates feasible spanning directions from the inward unit normal in the affine hull and, for $\alpha>0$, avoids nonsmooth boundary points almost surely when initialized outside that set (\Cref{sec:noised_alg}).
With randomized exploration adapted from \citep{greenstein2026achieving}, almost surely every accumulation point of either method is $d$-stationary (\cref{thm:main}).

\paragraph{Approximate $d$-stationarity under curvature.}
Under the stated curvature bounds on a $C^2$ constraint boundary, we use a basis adapted to the boundary geometry and inward tilting to jointly control the spanning coefficients and feasible step lengths. With Lipschitz continuity of $\nabla f$, this yields a finite-time guarantee of approximate $d$-stationarity along every unit feasible direction admitting a prescribed positive step length (\cref{thm:finite-time-rate,cor:eps-complexity}).

\paragraph{Finite-time PDM guarantees.}
For both methods, we establish an $O(1/K)$ rate in the proximal directional minimization (PDM) measure of \cite{BECK2022517}, assuming uniformly bounded direction-set sizes (\cref{thm:convergence_rate}). For the oracle method, vanishing PDM certifies $d$-stationarity at relative-interior points and non-descent along all $\eta$-feasible directions at boundary points. For sufficiently small noise, the noised method retains the interior guarantee and covers a slightly smaller cone at smooth boundary points.

\paragraph{Geometric analysis.}
We also establish results on the local geometry of feasible directions for closed convex sets, including characterizations of $\eta$-FDs at relative-interior and smooth boundary points (\Cref{lem:eta_FD_properties}). Geometry for directional spanning and trajectory behavior within general convex sets is developed in \Cref{sec:noised_alg_properties}.

\subsection{Background}
\label{sec:literature}

We review related approaches to constrained optimization, emphasizing
their information requirements, treatment of feasibility, and
stationarity guarantees.



\paragraph{Methods for linear constraints.}
The closest works to ours are \cite{feasible_directions_2020,BECK2022517}, which study the same objective structure under linear constraints. 
In that setting, the feasible direction cone is closed, explicitly computable, and finitely generated. 
These works establish convergence to stationary points and introduce the PDM optimality measure. 
The same spanning construction is classical in bound- and linearly constrained pattern search and generating-set search, which poll a finite set of generators of the tangent cone \citep{Lewis1999,Lewis2000,kolda2003optimization,Kolda2006}. 

Over a general closed convex set, the feasible direction cone need not have a finite set of generators, and sequences of feasible directions may converge to directions that are no longer feasible; \cref{alg:Oracle_GRFD} instead uses a finite set of feasible directions whose positive span contains the $\eta$-feasible directions. 
\citet{gratton2019direct} analyze direct search with sufficient decrease for smooth bound- and linearly constrained problems; we compare their complexity bound with ours after \cref{cor:eps-complexity}.

\paragraph{Difference-of-convex programming.}
When $f$ is convex, problem \eqref{eq:feasible_directions_beck_hallak} becomes a difference-of-convex (DC) program; see, e.g., \cite{tao1998dc}. 
Most DC algorithms guarantee convergence only to DC-critical points; see, e.g., \citep{abbaszadehpeivasti2024rate,phan2024difference}.
The optimality hierarchy for nonsmooth DC programs was developed in \citep{PRA17}, where directional stationarity was shown to be a stronger condition.
Under additional structure, directional stationarity can be attained by some DC methods; see, e.g., \citep{PRA17}.

Our setting differs in several aspects: we allow $f$ to be nonconvex, impose general convex constraints, and access $h$ through exact line searches, with no derivatives and without requiring access to the decomposition into $f$ and $g$. 
While DC methods benefit from first-order information, our focus is on procedures that use no derivatives and maintain feasibility.
Projection-free Frank--Wolfe methods likewise produce feasible iterates over general convex sets \citep{Frank1956,Jaggi2013}, but they require first-order information and a linear minimization oracle over $\mathcal{C}$.



\paragraph{Derivative-free constrained optimization.}
Conceptually, our approach can be affiliated with the Derivative-free methods (DFO) branch -- see \cite{larson2019derivative,conn2009introduction} for general surveys.  
DFOs for constrained optimization include direct-search approaches such as MADS \cite{audet2006mesh}, which achieve Clarke-type stationarity under general constraints. 
These methods typically rely on asymptotically dense sampling of directions. 
For nonlinear inequality constraints, \citet{Lucidi2002} keep all iterates feasible by sampling the objective along search arcs obtained by projecting rays onto the feasible set, with directions that approximately generate the cone of first-order feasible variations. 
Our methods use finite-direction constructions based on $\eta$-feasible directions for a finite-time PDM certificate; almost-sure subsequence $d$-stationarity is obtained from a randomized exploration component developed by the authors in \citep{greenstein2026achieving}.

\paragraph{Penalty and barrier methods.}
Many derivative-free approaches handle constraints via penalty formulations, allowing infeasible iterates \citep{liuzzi2010sequential,Fasano2014Linesearch}; the interior-point method of \citep{brilli2025interior} keeps iterates strictly feasible with respect to nonrelaxable inequality constraints through a log-barrier, while penalizing equality constraints.
In contrast, our framework maintains feasibility with respect to the whole set $\mathcal{C}$ at every iteration, including iterates on its boundary, leading to a different algorithmic design and analysis.

\subsection{Notation}
We use the common notation of convex analysis and continuous optimization literature.
Depending on the context, we use the index set $[q]:= \{1,2,\ldots,q\}$ for a natural number $q$, and $\Delta_{q}:= \{ x\in\real^q : \sum x_i = 1, x\geq 0\}$ is the $q$-dimensional simplex.
For a nonempty set $\mathcal{X}$, the set of all finite subsets of $\mathcal{X}$ is denoted by $\finitePow{\mathcal{X}} = \left\{C \subseteq \mathcal{X} \; | \; \left|C \right| < \infty\right\}.$
For a convex body $\mathcal{C}$ with a $C^2$ boundary and principal curvatures in $[\kappa_{\min},\kappa_{\max}]$ we write $r_{\mathrm{in}}:=1/\kappa_{\max}$ and $r_{\mathrm{out}}:=1/\kappa_{\min}$ for the inner and outer rolling-ball radii, and $\rho(x):=\mathrm{dist}(x,\rb{\mathcal{C}})$ for the depth of $x\in\mathcal{C}$.
For $x\in\mathcal{C}$ write $\Pi(x):=\arg\min_{y\in\rb{\mathcal{C}}}\norm{y-x}$.
Note that $\rb{\mathcal{C}}$ is not convex, so $\Pi(x)$ need not be a singleton; at the center of a ball it is the entire sphere.

\section{Preliminaries}\label{sec:preliminaries}
In this section we recall classical elements to be used throughout the paper, starting with the notions of positive span and cone of feasible directions. 
\begin{defin}[{\cite{feasible_directions_2020}}, Definition 2.1]\label{def:postive_span}    The positive span of a finite set of vectors $S = \{v_1, v_2, \ldots, v_k\} \subseteq \real^n$, where $k\in\mathbb{N}$, is denoted by $\pos{S}$ and is defined as the convex cone 
    \begin{equation*}
       \pos{S} := \left\{\sum\limits_{i=1}^k \lambda_i v_i \; | \; \lambda_i \ge 0 \;\; \forall i \in \{1,\ldots,k\}\right\}.
    \end{equation*}

\end{defin}

\begin{defin}[positive spanning set {\cite[Definition 2.2]{feasible_directions_2020}}]\label{def:postive_spanning_set} A finite set $S \subseteq \real^n$ is a \textbf{positive spanning set (PSS)} of a convex cone $C \subseteq \real^n$ if $\pos{S} = C$. In this case, $S$ is
said to positively span $C$.

\end{defin}

\begin{defin}[cone of feasible directions]\label{def:feasible_directions}
    Let $\mathcal{C} \subseteq \real^n$ be a convex set and let $\bar{x} \in \mathcal{C}$. 
    The \textbf{cone of feasible directions} is given by
    \begin{equation*}
        D_{\mathcal{C}}(\bar{x}) := \left\{w \; | \; \exists \lambda > 0, \;\; \bar{x} + \lambda w \in \mathcal{C}\right\}.
    \end{equation*}
\end{defin}
The next lemma records the elementary fact that feasible directions form a convex cone. A proof is included in \cref{sec:appendix_proofs_of_preliminaries} for completeness.
\begin{lem}\label{lem:feasible_directions_are_convex_cones}
    Let $\mathcal{C} \subseteq \real^n$ be a convex set, and let $\bar{x} \in \mathcal{C}$. Then $D_{\mathcal{C}}\left(\bar{x}\right)$ is a convex cone.
\end{lem}

\begin{defin}[directional ($d$-)stationarity]\label{def:d-stationary}
A point $\bar{x} \in \mathcal{C}$ is \textbf{$d$-stationary} for \eqref{eq:feasible_directions_beck_hallak} if
\begin{equation*}
    h'(\bar{x};\, d) \ge 0 \qquad \forall\, d \in D_{\mathcal{C}}(\bar{x}).
\end{equation*}
\end{defin}
Since $f$ is continuously differentiable and $g$ is convex, the map $d \mapsto h'(\bar{x};\, d)$ is continuous, so the same condition is equivalent to $h'(\bar{x};\, d) \ge 0$ for every $d \in \cl{D_{\mathcal{C}}(\bar{x})}$.

We will use the Euclidean distance from a point to a nonempty set, the induced distance between two nonempty sets, and the Hausdorff distance between compact sets.
\begin{defin}[Hausdorff distance]\label{def:hausdorff_distance}
    The Hausdorff distance between two nonempty compact sets $A, B \subseteq \real^n$ is
    \begin{equation*}
        d_H(A, B) := \max\left\{
            \sup_{a \in A} \inf_{b \in B} \norm{a - b},\ \
            \sup_{b \in B} \inf_{a \in A} \norm{a - b}
        \right\}.
    \end{equation*}
\end{defin}
\begin{defin}[Euclidean distance]\label{def:euclidean_distance}
    For $x \in \real^n$ and a nonempty set $A \subseteq \real^n$,
    \begin{equation*}
        \dist(x, A) := \inf_{a \in A} \norm{x - a}.
    \end{equation*}
    For nonempty sets $A, B \subseteq \real^n$,
    \begin{equation*}
        \dist(A, B) := \inf_{a \in A} \dist(a, B).
    \end{equation*}
\end{defin}

Our approach relies on foundational notions 
which we introduce and utilize in this section 
-- proofs not given in the text are in \cref{sec:appendix_proofs_of_preliminaries}.

\begin{defin}[affine hull]
    Let $S \subseteq \real^n$. Then the affine hull of $S$ is given by
    \begin{equation*}
        \aff{S} := \left\{ \sum\limits_{i=1}^k \alpha_i x_i\; | \; k > 0, \;\; x_1,\ldots,x_k \in S, \;\; \alpha_1,\ldots,\alpha_k \in \real, \;\; \sum\limits_{i=1}^k \alpha_i = 1 \right\}.
    \end{equation*}
\end{defin}
The relative interior of a set is defined via the affine hull. 
\begin{defin}[relative interior]
    Let $S \subseteq \real^n$. Then the relative interior of $S$ is given by
    \begin{equation*}
        \ri{S} := \left\{x \in S \; | \; \exists \epsilon > 0, \;\; \mathcal{B}\left(x, \epsilon\right) \cap \aff{S} \subseteq S\right\}.
    \end{equation*}
\end{defin}
We next define the operation of shifting a set by a reference point. 
\begin{defin}
    Let $S \subseteq \real^n$ be a nonempty set, and let $\bar{x} \in S$. Then the set $S - \bar{x}$ is given by
    \begin{equation*}
        S - \bar{x} = \left\{y \in \real^n\; | \; y + \bar{x} \in S\right\}.
    \end{equation*}
\end{defin}
\begin{remark}
    Note that for every nonempty set $S \subseteq \real^n$ and every $\bar{x} \in S$, we have $0 \in S - \bar{x}$.
\end{remark}
The next lemma provides key foundational properties we will use, some are derived from \cite[Proposition 2.40]{RW98}.
\begin{lem}[foundational properties]\label{lem:foundational properties}
    Let $S \subseteq \real^n$ be a nonempty set.
    Then $\aff{S} = \aff{\cl{S}}$, and for $\bar{x}, \bar{y} \in S$ it holds that $\aff{S - \bar{x}} = \aff{S - \bar{y}}.$
    If $S$ is in addition convex, then $\ri{S} = \ri{\cl{S}}$, $\ri{S} \ne \emptyset$, and the line segment principle extends to the relative interior, that is, for $x \in \ri{S}$, $y \in \cl{S}$, and any $\lambda \in (0, 1]$, it holds that $\lambda x + \left(1 - \lambda\right) y \in \ri{S}. $
\end{lem}

For a linear subspace $L\subseteq\real^n$, we write $\mathbb{S}_L := \{\omega \in L : \norm{\omega}=1\}$ for the unit sphere of $L$. When $L = \aff{\mathcal{C}-\bar{x}}$ for some $\bar{x}\in\mathcal{C}$, this notation is independent of the reference point by the preceding lemma.

We continue with the definition of the relative boundary.

\begin{defin}[Relative Boundary]
    Let $\mathcal{C} \subseteq \real^n$. The relative boundary of $\mathcal{C}$, $\rb {\mathcal{C}}$, is given by $\rb{\mathcal{C}} = \cl{\mathcal{C}} \setminus \ri{\mathcal{C}}. $
\end{defin}

\begin{defin}[Smoothness at Boundary Point]\label{def:1}
    Let $\mathcal{C} \subseteq \real^n$ be a closed convex set. Then $\mathcal{C}$ is
    \textit{smooth} at a boundary point $\bar{x} \in \rb{\mathcal{C}}$ if
    \begin{equation*}
        \cl{D_\mathcal{C}(\bar{x})}
        = \left\{d \in \real^n \; | \; \inner{d}{u} \ge 0\right\}
          \cap \aff{\mathcal{C} - \bar{x}}
    \end{equation*}
    for some non-zero vector $u \in \aff{\mathcal{C} - \bar{x}}$.
\end{defin}

\begin{lem}\label{lem:feasible_directions_contained_in_affine_hull}
    Let $\mathcal{C} \subseteq \real^n$ be a closed convex set, and let $\bar{x} \in \mathcal{C}$. Then,
    \begin{equation*}
        D_\mathcal{C}(\bar{x}) \subseteq \cl{D_\mathcal{C}(\bar{x})} \subseteq \aff{\mathcal{C} - \bar{x}}.
    \end{equation*}
\end{lem}
\begin{proof}
    The inclusion $D_\mathcal{C}(\bar{x}) \subseteq \cl{D_\mathcal{C}(\bar{x})}$ is immediate.
    Let $d \in D_\mathcal{C}(\bar{x})$. By definition, there exists $\epsilon_d > 0$ such that $\bar{x} + \epsilon_d d =: x_d \in \mathcal{C}$, hence $d = \frac{x_d - \bar{x}}{\epsilon_d}$.
    Since $0 \in \mathcal{C}-\bar{x}$, the set $\aff{\mathcal{C}-\bar{x}}$ is a linear subspace, and therefore closed under scalar multiplication; as $x_d - \bar{x} \in \mathcal{C}-\bar{x} \subseteq \aff{\mathcal{C}-\bar{x}}$, we conclude $d \in \aff{\mathcal{C}-\bar{x}}$.
    Hence $D_\mathcal{C}(\bar{x}) \subseteq \aff{\mathcal{C}-\bar{x}}$.

    Now, let $\bar{d} \in \cl{D_\mathcal{C}(\bar{x})}$. Then there exists a sequence $\left\{d_i\right\}_{i \ge 0} \subseteq D_\mathcal{C}(\bar{x})$, such that $d_i \xrightarrow{i\to\infty} \bar{d}$. Since $D_\mathcal{C}(\bar{x}) \subseteq \aff{\mathcal{C}-\bar{x}}$, it follows that $\left\{d_i\right\}_{i \ge 0} \subseteq \aff{\mathcal{C}-\bar{x}}$. Since $\mathcal{C}$ is a finite dimensional set, its affine hull is closed, and hence $\bar{d} \in \aff{\mathcal{C}-\bar{x}}$. We conclude that $\cl{D_\mathcal{C}(\bar{x})} \subseteq \aff{\mathcal{C}-\bar{x}}$.
\end{proof}

\begin{lem}\label{lem:feasible_directions_at_relative_interior}
    Let $\mathcal{C} \subseteq \real^n$ be a closed convex set, and let
    $\bar{x} \in \ri{\mathcal{C}}$. Then $D_\mathcal{C}(\bar{x}) = \aff{\mathcal{C} - \bar{x}}$.
\end{lem}
\begin{proof}
    The inclusion $D_\mathcal{C}(\bar{x}) \subseteq \aff{\mathcal{C}-\bar{x}}$ is
    \cref{lem:feasible_directions_contained_in_affine_hull}.
    For the reverse inclusion, let $\tilde{d} \in \aff{\mathcal{C}-\bar{x}}$,
    $\tilde{d} \ne 0$.
    Since $\bar{x} \in \ri{\mathcal{C}}$, there exists $\epsilon > 0$ with
    $\mathcal{B}(\bar{x}, \epsilon) \cap \aff{\mathcal{C}} \subseteq \mathcal{C}$.
    Setting $t = \epsilon / (2\norm{\tilde{d}})$ gives $\bar{x} + t\tilde{d} \in
    \mathcal{B}(\bar{x},\epsilon) \cap \aff{\mathcal{C}} \subseteq \mathcal{C}$,
    so $\tilde{d} \in D_\mathcal{C}(\bar{x})$.
\end{proof}

The affine hull, relative interior, and relative boundary, together with the smoothness notion of \Cref{def:1} and the basic feasible-direction properties of Lemmas \ref{lem:feasible_directions_contained_in_affine_hull}, \ref{lem:feasible_directions_at_relative_interior}, play a central role in defining $\eta$-feasible directions, which are discussed in the following section.

\section{{$\eta$}-Feasible Directions}\label{sec:eta_feasible}
To address the geometric difficulties induced by general convex constraints, specifically, the structure of the feasible-direction cone near the boundary, we introduce the concept of \textit{$\eta$-feasible directions} ($\eta$-FD).
These feasible directions are angularly separated from the relative boundary of the feasible-direction cone; a geometric interpretation is given in \cref{rem:ang eta}.
\begin{defin}[{$\eta$}-Feasible Directions ($\eta$-FD)]\label{def:eta_feasible_directions}
    Let $\eta \in \left[ 0,1\right]$ and $x \in \mathcal{C}$. 
    A direction $0 \ne d \in D_{\mathcal{C}}(x)$ is called an $\eta$-FD with respect to $\mathcal{C}$ at $x$, if 
    it holds that
    \begin{equation*}
        \dfrac{\left| \langle d, \bar{d}\rangle \right|}{\left\Vert d \right\Vert \cdot \left\Vert \bar{d} \right\Vert} \le \sqrt{1 - \eta^2}, \qquad \forall\bar{d} \in \rb{D_\mathcal{C}(x)}, \bar{d} \ne 0.
    \end{equation*}
    The set of all $\eta$-FDs at $x \in \mathcal{C}$ is denoted by $D^{\eta}_\mathcal{C}(x)$.
\end{defin}
\begin{remark}[angular interpretation]\label{rem:ang eta} 
Let $\theta=\arcsin(\eta)\in[0,\pi/2]$. Then the defining inequality is equivalent to \[ 
\arccos\left( \frac{|\langle d,\bar d\rangle|} {\|d\|\,\|\bar d\|} \right) \geq \theta \qquad \forall\bar{d} \in \rb{D_\mathcal{C}(x)}, \bar{d} \ne 0. 
\] 
Thus, an $\eta$-feasible direction is a feasible direction that remains uniformly separated, by an angle of at least $\theta$, from the lines generated by the relative-boundary directions of $D_{\mathcal C}(x)$. Equivalently, larger values of $\eta$ correspond to directions that lie deeper inside the feasible-direction cone. 
\end{remark}

\begin{example}[$\eta$-FDs in $\ell_2$-ball] \label{ex:eta FD in l2 ball}
    We illustrate the structure of $\eta$-FDs when the feasible set is the $\ell_2$-norm ball.
    If $\mathcal{C} = \mathcal{B}_2[0,1] = \left\{x \in \real^n \; | \; \left\Vert x \right\Vert_2 \le 1\right\}$.
    Then $\rb{\mathcal{C}} = \left\{x \in \real^n \; | \; \left\Vert x \right\Vert_2 = 1\right\}$, and 
    \begin{equation*}
        D^{\eta}_\mathcal{C}(x) = \begin{cases}
            \left\{0 \ne d \; | \; \left\langle x, \dfrac{d}{\left\Vert d \right\Vert} \right\rangle \le -\eta\right\}, & \text{if $\left\Vert x \right\Vert_2 = 1$},\\
            \real^n - \left\{0\right\}, & \text{otherwise.}
        \end{cases}
    \end{equation*}
    

\end{example}

The next lemma gives two useful characterizations of $\eta$-feasible directions. 
Although \Cref{def:eta_feasible_directions} is stated in terms of the relative boundary of $D_{\mathcal C}(\bar x)$, the condition simplifies in the two cases needed later.

First, if $D_{\mathcal C}(\bar x)$ coincides with the affine hull of $\mathcal C-\bar x$, then every nonzero feasible direction is an $\eta$-FD. 
Second, if $\mathrm{cl}(D_{\mathcal C}(\bar x))$ is a half-space in this affine hull, then $\eta$-feasibility is equivalent to a single normalized inner-product inequality.

These two cases are instantiated in \Cref{lem:eta_FD_properties} below, which identifies relative-interior points with part~(a) and smooth boundary points with part~(b).





\begin{lem}[on the characterizations of $\eta$-FD]\label{lem:eta_FD_properties_without_smoothness}
    Let $\bar{x} \in \mathcal{C}$, and $\eta \in [0,1]$.
    \begin{enumerate}
        \item \textbf{Entire affine hull:} If $D_\mathcal{C}(\bar{x}) = \aff{\mathcal{C}-\bar{x}}$, then every nonzero $ d \in D_\mathcal{C}(\bar{x})$ is an $\eta$-FD.
        \item \textbf{Half-space within the affine hull:} If $\cl{D_\mathcal{C}(\bar{x})} = \left\{d \in \real^n \; | \; \langle d, u \rangle \ge 0\right\} \cap \aff{\mathcal{C}-\bar{x}}$ where $u \in \aff{\mathcal{C}-\bar{x}}$ and $\left\Vert u \right\Vert = 1$, then $d \in D_\mathcal{C}(\bar{x})$ such that $\left\Vert d \right\Vert = 1$ is an $\eta$-FD if and only if $\langle d, u \rangle \ge \eta.$
    \end{enumerate}
\end{lem}

\begin{proof}
    Suppose that $D_\mathcal{C}(\bar{x}) = \aff{\mathcal{C}-\bar{x}}$. 
    Since $D_\mathcal{C}(\bar{x}) = \aff{\mathcal{C} - \bar{x}}$ is an affine space, it holds that $\cl{\aff{\mathcal{C} - \bar{x}}} = \ri{\aff{\mathcal{C} - \bar{x}}} = \aff{\mathcal{C}-\bar{x}}.$
    Hence, the first part follows trivially from 
    \begin{equation*}
        \rb{D_\mathcal{C}(\bar{x})} = \rb{\aff{\mathcal{C}-\bar{x}}} = \cl{\aff{\mathcal{C}-\bar{x}}} - \ri{\aff{\mathcal{C}-\bar{x}}} = \emptyset.
    \end{equation*} 

    We proceed to the second part. 
    Let $W = \left\{u, v_1, \ldots, v_{\dim\mathcal{C}-1}\right\}$ be some orthogonal spanning set of $\aff{\mathcal{C}-\bar{x}}$ that includes $u$, such that $\left\Vert v_1 \right\Vert = \ldots = \left\Vert v_{\dim\mathcal{C}-1}\right\Vert = 1$. 
    Note that
    \begin{equation*}
        \rb{D_\mathcal{C}(\bar{x})} = \rb{\left\{d \in \real^n \; | \; \langle u, d \rangle \ge 0\right\} \cap \aff{\mathcal{C} - \bar{x}}} = \left\{d \in \real^n \; | \; \langle u, d \rangle = 0\right\} \cap \aff{\mathcal{C} - \bar{x}},
    \end{equation*}
     assume that $\langle d, u \rangle \ge \eta$, and let $\tilde{d} \in \rb{D_\mathcal{C}(\bar{x})}$. Then there exist coefficients $\tilde{a}_1,\ldots,\tilde{a}_{\dim\mathcal{C}-1} \in \real$, such that $\tilde{d} = 0 \cdot u + \sum\limits_{i=1}^{\dim\mathcal{C}-1}\tilde{a}_i v_i.$
    Meanwhile, since $d \in D_\mathcal{C}(\bar{x}) \subseteq \cl{D_\mathcal{C}(\bar{x})}$, there exist coefficients $c, a_1, \ldots, a_{\dim\mathcal{C}-1} \in \real$, such that 
    \begin{equation*}
        d = c \cdot u + \sum\limits_{i=1}^{\dim\mathcal{C}-1} a_i v_i.
    \end{equation*}
    Since $\eta \le \langle d, u \rangle = c$, and $1 = \left\Vert d \right\Vert^2 = c^2 + \sum\limits_{i=1}^{\dim\mathcal{C}-1} a_i^2$, it follows that $\sum\limits_{i=1}^{\dim\mathcal{C}-1} a_i^2 \leq  1 - \eta^2$.
    Note that $\left| \langle d, \tilde{d} \rangle\right| = \max \left\{ \langle d, \tilde{d} \rangle, \langle d, -\tilde{d} \rangle\right\}$, and that both
    \begin{equation*}
        \langle d, \tilde{d} \rangle \le \sup\limits_{q \in \mathrm{rb}{D_\mathcal{C}(\bar{x})}} \left\{ \langle d, q \rangle \; | \; \left\Vert q \right\Vert^2 \le \left\Vert \tilde{d} \right\Vert^2\right\} \; \text{ and } \; \langle d, -\tilde{d} \rangle \le \sup\limits_{q \in \mathrm{rb}{D_\mathcal{C}(\bar{x})}} \left\{ \langle d, q \rangle \; | \; \left\Vert q \right\Vert^2 \le \left\Vert \tilde{d} \right\Vert^2\right\}
    \end{equation*}
    hold true. 
    Therefore, 
    \begin{equation*}
        \left| \langle d, \tilde{d} \rangle\right| \le \sup\limits_{q \in \mathrm{rb}{D_\mathcal{C}(\bar{x})}} \left\{ \langle d, q \rangle \; | \; \left\Vert q \right\Vert^2 \le \left\Vert \tilde{d} \right\Vert^2\right\}.
    \end{equation*}
    Rewriting the supremum in term of the orthogonal decomposition, we derive
    \begin{equation*}
        \sup\limits_{b_1,\ldots,b_{\dim\mathcal{C}-1}} \left\{ \sum\limits_{i=1}^{\dim\mathcal{C}-1} a_i b_i \; | \; \sum\limits_{i=1}^{\dim\mathcal{C}-1} b_i^2 \le \left\Vert \tilde{d} \right\Vert^2\right\}.
    \end{equation*}
    Using Cauchy-Schwarz,  $ \sum\limits_{i=1}^{\dim\mathcal{C}-1} a_i b_i \le \left\Vert a \right\Vert \cdot \left\Vert b \right\Vert \le \left\Vert a \right\Vert \cdot \left\Vert \tilde{d} \right\Vert, $
    and the upper bound is attained for $b = \left\Vert \tilde{d} \right\Vert a / \left\Vert a \right\Vert$. 
    Therefore,
    \begin{equation*}
        \left| \langle d, \tilde{d} \rangle \right| \le \left(\sqrt{\sum\limits_{i=1}^{\dim\mathcal{C}-1} a_i^2}\right) \cdot \left\Vert \tilde{d} \right\Vert \le \sqrt{1 - \eta^2} \cdot \left\Vert \tilde{d} \right\Vert.
    \end{equation*}
    Combining with the assumption that $\left\Vert d \right\Vert = 1$, we have that
    \begin{equation}\label{eq:sm145}
        \dfrac{\left| \langle d, \tilde{d} \rangle\right|}{\left\Vert d \right\Vert \cdot \left\Vert \tilde{d} \right\Vert} \le \sqrt{1 - \eta^2}.
    \end{equation}
    Thus,  $\langle d, u \rangle \ge \eta$ implies that $d$ is $\eta$-FD.
\medskip

    In the other direction, assume that \eqref{eq:sm145}
    holds true for every $\tilde{d} \in \rb{D_\mathcal{C}(\bar{x})}$. 
    As before, let $d = c \cdot u + \sum\limits_{i=1}^{\dim\mathcal{C}-1} a_i \cdot v_i$
    be the decomposition of $d$ according to the normalized orthogonal basis $W$. 
    Note that if $c < 0$, then $\langle d, u \rangle < 0$, and hence $d \notin \cl{D_\mathcal{C}(\bar{x})}$, in contradiction to the assumption that $d \in D_\mathcal{C}(\bar{x})$, and if $c > 1$, then $\left\Vert d \right\Vert > 1$, in contradiction to the assumption that $\left\Vert d \right\Vert = 1$.
    
    Consequently,  $c \in [0,1]$.
    If $c = 1$, then $\langle d, u \rangle = 1 \ge \eta$, and the claim follows. 
    Otherwise, assume $c < 1$, and hence $\sum\limits_{i=1}^{\dim\mathcal{C}-1} a_i^2 > 0$. 
    Consider the vector $\bar{d} = \sum\limits_{i=1}^{\dim\mathcal{C}-1} a_i \cdot v_i.$
    Since $\langle \bar{d}, u \rangle = 0$, and $\bar{d} \in \cl{D_\mathcal{C}(\bar{x})}$, it follows that $\bar{d} \in \rb{D_\mathcal{C}(\bar{x})}$, and therefore, $\tilde{d}=\bar{d}$ satisfies \eqref{eq:sm145}.
    Since $1 = \left\Vert d \right\Vert^2 = c^2 + \sum\limits_{i=1}^{\dim\mathcal{C}-1} a_i^2$, it follows that
    \begin{equation*}
        \left\Vert\bar{d}\right\Vert = \sqrt{\sum\limits_{i=1}^{\dim\mathcal{C}-1} a_i^2} = \sqrt{1 - c^2}.
    \end{equation*}
    Additionally, $\left| \langle d, \bar{d} \rangle \right| = \sum\limits_{i=1}^{\dim\mathcal{C}-1} a_i^2 = 1 - c^2$. 
    Combining the above then yields that 
    \begin{equation*}
        \dfrac{\left| \langle d, \tilde{d} \rangle\right|}{\left\Vert d \right\Vert \cdot \left\Vert \tilde{d} \right\Vert} = \dfrac{1 - c^2}{1 \cdot \sqrt{1 - c^2}} \le \sqrt{1 - \eta^2} \Rightarrow c \ge \eta,
    \end{equation*}
    where the last transition follows from the fact that $c, \eta \in \left[0,1\right]$. 
    Thus, $\langle d, u \rangle \ge \eta$,
    which concludes the proof for $x \in \rb{\mathcal{C}}$.
    
\end{proof}

The following corollary specializes \cref{lem:eta_FD_properties_without_smoothness} to interior points and smooth boundary points; it follows directly from \cref{def:1},
\cref{lem:feasible_directions_at_relative_interior}, and
\cref{lem:eta_FD_properties_without_smoothness}.

\begin{corollary}\label{lem:eta_FD_properties}
    Let $\bar{x} \in \mathcal{C}$, and $\eta \in [0,1]$.
    \begin{enumerate}[label=(\alph*)]
        \item If $\bar{x} \in \ri{\mathcal{C}}$, then $D_\mathcal{C}(\bar{x}) = \aff{\mathcal{C}-\bar{x}}$, and every $0 \ne d \in D_\mathcal{C}(\bar{x})$ is an $\eta$-FD.
        \item If $\bar{x} \in \rb{\mathcal{C}}$ and $\mathcal{C}$ is smooth at $\bar{x}$, then $d \in D_\mathcal{C}(\bar{x})$ with $\|d\| = 1$ is an $\eta$-FD if and only if
        \[
            \langle d, u \rangle \ge \eta,
        \]
        where $u \in \aff{\mathcal{C} - \bar{x}}$ is a unit vector satisfying
        \[
            \cl{D_\mathcal{C}(\bar{x})} = \left\{ \tilde{d} \in \real^n \; | \; \langle u, \tilde{d} \rangle \ge 0 \right\} \cap \aff{\mathcal{C} - \bar{x}}.
        \]
    \end{enumerate}
\end{corollary}
\begin{proof}
    \textit{Part (a).} By \cref{lem:feasible_directions_at_relative_interior},
    $D_\mathcal{C}(\bar{x}) = \aff{\mathcal{C}-\bar{x}}$. Since $\aff{\mathcal{C}-\bar{x}}$
    is a linear subspace, its relative boundary is empty, so the condition in
    \cref{def:eta_feasible_directions} is vacuously satisfied for every nonzero $d$.
    Part~1 of \cref{lem:eta_FD_properties_without_smoothness} then gives that every
    $0 \ne d \in D_\mathcal{C}(\bar{x})$ is an $\eta$-FD.

    \textit{Part (b).} By \cref{def:1}, smoothness at $\bar{x}$ means that
    \begin{equation*}
        \cl{D_\mathcal{C}(\bar{x})}
        = \left\{d \in \real^n \; | \; \inner{d}{u} \ge 0\right\}
          \cap \aff{\mathcal{C}-\bar{x}}
    \end{equation*}
    for some non-zero $u \in \aff{\mathcal{C}-\bar{x}}$. Rescaling $u \leftarrow u/\norm{u}$
    yields a unit vector satisfying the half-space assumption of
    Part~2 of \cref{lem:eta_FD_properties_without_smoothness}, from which the
    $\eta$-FD characterization $\inner{d}{u} \ge \eta$ follows directly.
\end{proof}

\section{The Oracle-Based Regularized Feasible Directions Method}\label{sec:directions_oracle_and_alg}
This section introduces the Oracle-Based Regularized Feasible Directions Method (\cref{alg:Oracle_GRFD}).
At each iteration, an oracle produces a positively spanning set of $\eta$-feasible directions, and exact line searches are applied to the associated univariate regularized subproblems.

The methods access $h$ only through the exact regularized line searches in line~7 of \cref{alg:Oracle_GRFD} and line~5 of \cref{alg:GRFD} (\cref{rem:line-search-access}).
\Cref{alg:Oracle_GRFD} accesses $\mathcal{C}$ through an $\eta$-direction oracle, defined next, and through the step-size map $\tau$.

\begin{defin}[$\eta$-Direction Oracle]\label{def:eta_directions_oracle}
    Let $\mathcal{C} \subseteq \real^n$ be a nonempty closed convex set.
    A mapping $\mathcal{A} : \mathcal{C} \times (0,1] \rightarrow \finitePow{\real^n}$ is called an $\eta$-Direction Oracle ($\eta$-DO) if for every $\bar{x} \in \mathcal{C}$ and $\eta \in (0,1]$, the output $X = \mathcal{A}(\bar{x},\eta)$ satisfies:
    \begin{enumerate}
        \item Each $d \in X$ lies in the feasible direction cone: $d \in D_\mathcal{C}(\bar{x})$.
        \item The positive span of $X$ contains all $\eta$-feasible directions at $\bar{x}$: $  D^\eta_\mathcal{C}(\bar{x}) \subseteq \pos{X}$.
    \end{enumerate}
\end{defin}

\begin{remark}[On the case $\eta=0$]
\label{rem:eta_zero_oracle}
The oracle is defined only for $\eta>0$. The reason is that, for general closed convex sets, the full feasible-direction cone $D_{\mathcal C}(\bar x)$ need not be closed and therefore need not admit a finite positive spanning set of feasible directions. For example, if $\mathcal C=\mathcal B_2[0,1]$ and $\|x\|_2=1$, then
\begin{equation*}
    D_{\mathcal{B}[0,1]}(x) = \left\{d \in \real^n \; | \; \langle x, d \rangle < 0\right\},
\end{equation*}
which is not closed and cannot be represented as the positive span of finitely many feasible directions. When $\mathcal C$ is polyhedral, $D_{\mathcal C}(\bar x)$ is closed and finitely generated; see, e.g., \cite{BECK2022517}. In that case, the oracle may be extended to $\eta=0$ by returning a positive spanning set of $D_{\mathcal C}(\bar x)$, and we use this extension whenever applicable.
\end{remark}

\begin{assumption}[Constraint-set access]\label{assump:constraint-access}
Let $\mathcal{C} \subseteq \real^n$ be a nonempty closed convex set.
\Cref{alg:Oracle_GRFD} is given an $\eta$-direction oracle (\cref{def:eta_directions_oracle}) and, for every $x \in \mathcal{C}$ and $v \in \real^n \setminus \{0\}$, the value
$\tau(x;v) := \sup\{t \ge 0: x + tv \in \mathcal{C}\}$.
\end{assumption}

\begin{remark}[Line-search access]\label{rem:line-search-access}
Both methods access $h$ only through exact solutions of the regularized univariate problems
\begin{equation*}
\min_{0\le t\le \min\{\tau(x;v),\,r\}} \varphi_v(t), \qquad \varphi_v(t) := h(x+tv) + \frac{\gamma}{2}t^2,
\end{equation*}
as in~\cite{BECK2022517}. For nonconvex $h$ this amounts to global minimization of a nonconvex univariate function; we adopt it as a modeling assumption and do not pursue inexact line searches here. The oracle uses neither derivatives of $h$ nor its decomposition into $f$ and $g$.
\end{remark}

\normalem
\noindent	\begin{algorithm}[h]
\caption{Oracle-Based Regularized Feasible Directions Method}\label{alg:Oracle_GRFD}
	\KwIn{$x_0 \in \mathcal{C}$, $\eta \in (0, 1]$, $\gamma>0$, $r>0$, \text{ $\eta$-DO } $\mathcal{A}: \mathcal{C} \times (0,1] \rightarrow \mathcal{P}_{fin}\left(\real^n\right)$} 
	
        \For{$k=0\ldots,K$}{

                $\tilde{v}_{\mathrm{ex}}^k \sim \unif{\mathbb{S}_L}$, \quad $L := \aff{\mathcal{C}-x_k}$

                $V^k = \mathcal{A}(x_k, \eta)$ 

                $\tilde{V}^k := \left\{\tilde{v}_1^k, \ldots, \tilde{v}_{s_k}^k\right\} = V^k \cup \left\{\tilde{v}_{\mathrm{ex}}^k\right\}$

                \For{$i \in \{1,\ldots,s_k\}$}{

                    $\tau(x_k;\tilde{v}_i^k) := \sup\{t \ge 0: x_k + t \tilde{v}_i^k \in \mathcal{C}\}$

                    $t_i^k \in \argmin\limits_{0 \le t \le \min \left\{\tau(x_k;\tilde{v}_i^k), r\right\}} h(x_k + t \tilde{v}_i^k) + \dfrac{\gamma}{2}t^2$
                
                }

                $i^* \in \argmin\limits_{i \in \{1,\ldots,s_k\}} h\left(x_k + t_i^k \tilde{v}_i^k \right)$                

                $x_{k+1} = x_k + t_{i^*}^k \tilde{v}_{i^*}^k$

        }
    
\end{algorithm}
\ULforem
\medskip

Below we provide classes of sets for which an \(\eta\)-DO can be implemented straightforwardly.

\paragraph{1. Polyhedral Constraints.}  
Suppose that $ \mathcal{C} = \left\{ x \in \real^n \; | \; Ax \le b \right\}, $ where \( A \in \real^{m \times n} \) and \( b \in \real^m \). Let \( A_i \) denote the \(i\)th row of \( A \), and define the active set at \( \bar{x} \in \mathcal{C} \) by $ \mathcal{I}_\mathcal{C}(\bar{x}) := \left\{ i \; | \; [A\bar{x}]_i = b_i \right\}. $ 
Then, the feasible directions at \( \bar{x} \) are given by
\begin{equation*}
    D_\mathcal{C}(\bar{x}) = \left\{ d \in \real^n \; | \; A_i d \le 0 \quad \forall i \in \mathcal{I}_\mathcal{C}(\bar{x}) \right\}.
\end{equation*}
Computing a finite set of vectors that positively span this cone corresponds to converting its \(\mathcal{H}\)-representation to a \(\mathcal{V}\)-representation -- a well-studied task in polyhedral theory. Efficient algorithms for this conversion exist depending on the feasible set; see \cite[Section 9]{fukuda2016lecture}. 

\paragraph{2. Smooth Convex Sets.}  
When the boundary of \( \mathcal{C} \) is smooth at every point, such as when $\mathcal{C}$ is the \( \ell_2 \)-ball (also see \Cref{ex:eta FD in l2 ball}) or an ellipsoid,  the non-exploration output of \cref{alg:construct_directions} with zero noise is an $\eta$-DO (\cref{lem:tilde_v_spans_eta_feasible}).
This holds for every $\rho_0\ge 0$. Under the curvature assumptions of \cref{lem:rolling_balls} and with $\rho_0=r_{\mathrm{in}}/2$, this $\eta$-DO also has the spanning-coefficient bounds of \cref{def:quantitative-eta-oracle} (\cref{prop:quantitative-spanning}).

\paragraph{3. Orthogonally-Structured Constraints.}
Suppose that 
\begin{equation}\label{eq:orthogonally_structured}
    \mathcal{C} = \left\{ x \in \real^n \; | \; h_i(x) \le 0, \;  i = 1, \dots, k \right\}, 
\end{equation}
where each \( h_i \) is convex and differentiable, and \( \nabla h_i(x) \ne 0 \) whenever \( h_i(x) = 0 \). 
A sufficient condition for the latter is Slater’s condition -- that is, the existence of some \( \hat{x} \) for which \( h_i(\hat{x}) < 0 \) for all \( i \). Additionally, assume that the gradients of active constraints are orthogonal.
This structure includes various settings with mutually exclusive or disjoint constraints. 
For $\eta\in (0,1]$, under these assumptions, the set of \(\eta\)-feasible directions at \( \bar{x} \) is given by
\begin{equation*}
    D^\eta_\mathcal{C}(\bar{x}) = \left\{ d \in \real^n \; | \; d \ne 0, \quad \left\langle \frac{\nabla h_i(\bar{x})}{\| \nabla h_i(\bar{x}) \|}, \frac{d}{\| d \|} \right\rangle \le -\eta \quad \forall i \in \mathcal{I}(\bar{x}) \right\},
\end{equation*}
where \( \mathcal{I}(\bar{x}) := \{ i \; | \; h_i(\bar{x}) = 0 \} \) is the set of active constraints at \( \bar{x} \).

Define \( u_i = -\nabla h_i(\bar{x}) / \| \nabla h_i(\bar{x}) \| \) and \( v = \sum_{i \in \mathcal{I}(\bar{x})} u_i \). Let \( \bar{\eta} = \eta / (2\sqrt{n}) \), and let $\{w_j\}$ 
be an orthonormal basis for the orthogonal complement of \( \{u_i\}_{i \in \mathcal{I}(\bar{x})} \). 
Then the set
\begin{align}\label{eq:spanning_directions_for_orthogonal_constraint_gradients}
    V := &\left\{ v \right\} 
    \cup \bigcup_{i \in \mathcal{I}(\bar{x})} \left\{ (1 - \bar{\eta}) u_i + \bar{\eta} v \right\} 
    \cup \bigcup_{j=1}^{n - |\mathcal{I}(\bar{x})|} \left\{ (1 - \bar{\eta}) w_j + \bar{\eta} v,\, -(1 - \bar{\eta}) w_j + \bar{\eta} v \right\}
\end{align}
contains only feasible directions, and its positive span includes \( D^\eta_\mathcal{C}(\bar{x}) \); The proof is nearly identical to that of \cref{lem:tilde_v_spans_eta_feasible}.
Thus, the mapping \( \mathcal{A}(\bar{x}, \eta) = V \) defines a valid \(\eta\)-DO.

\paragraph{Access for the constructions above.}
\Cref{alg:Oracle_GRFD} uses only the $\eta$-DO and the map $\tau$ of \cref{assump:constraint-access}.
Polyhedra and orthogonally structured constraints supply an $\eta$-DO by the constructions above, and do not use \cref{alg:construct_directions}.
For smooth convex sets, the $\eta$-DO is the non-exploration output of \cref{alg:construct_directions} with zero noise, which at $\rho_0=0$ requires, at each query point $x \in \mathcal{C}$, a relative-interior membership test, an orthonormal basis of $\aff{\mathcal{C}-x}$, and, if $x \in \rb{\mathcal{C}}$, the unique unit vector $u$ of \cref{def:1}.
We assume that these three items are available from the description of $\mathcal{C}$, as they are for the Euclidean ball and for ellipsoids; smoothness of the boundary does not, by itself, yield them for an arbitrary smooth convex set.
A membership oracle can be used to estimate them, but does not compute them exactly.
Implementation of the same construction at $\rho_0>0$ is discussed in \cref{rem:rho0-default}.
The map $\tau$ is obtained by one-dimensional root finding along the ray. For polyhedra, one may instead compute the distance to each linear constraint.


\section{The Noised Regularized Feasible Directions Method}\label{sec:noised_alg}
The \textit{Noised Regularized Feasible Directions Method}
(\cref{alg:GRFD}) replaces the $\eta$-direction oracle with an
explicit construction based on local boundary geometry
(\cref{alg:construct_directions}). It combines basis directions,
inward tilting, and random perturbations. For $\alpha>0$, the
iterates almost surely avoid nonsmooth boundary points when
initialized outside that set
(\cref{lem:always_smooth_almost_surely}).

The threshold $\rho_0\ge0$ determines where boundary geometry is
used: $\rho_0=0$ restricts its use to boundary points, whereas
$\rho_0>0$ also includes nearby interior points. The required
constraint queries are described at the end of
\cref{sec:directions_oracle_and_alg} and in
\cref{rem:rho0-default}.

\normalem
\begin{algorithm}[htp]
\caption{\textsc{ConstructSearchDirections}}\label{alg:construct_directions}
\KwIn{$x \in \mathcal{C}$, $\eta \in (0,1]$, $\alpha \in [0,1)$, $\rho_0 \ge 0$}
\KwOut{Search directions $\tilde{V}$ (including exploration direction)}
\tcp{queries $\rho(x)$; if $\rho\le\rho_0$, also $\Pi(x)$}
$\tilde{V} \leftarrow \emptyset$\;
$\bar{\eta} \leftarrow \dfrac{\eta}{2\sqrt{\dim\mathcal{C}}}$\;

\tcp{\textbf{Step 1: Exploration direction}}
$\tilde{v}_{\mathrm{ex}} \sim \unif{\mathbb{S}_L}$, \quad $L := \aff{\mathcal{C}-x}$\;

\tcp{\textbf{Step 2: Base directions}}
$\rho \leftarrow \rho(x)$\;
\uIf{$\rho > \rho_0$}{
    $W \leftarrow \{w_1, \ldots, w_{\dim \mathcal{C}}\}$, an orthonormal basis of $\aff{\mathcal{C} - x}$\;
    $U \leftarrow \{w_1, -w_1, \ldots, w_{\dim \mathcal{C}}, -w_{\dim \mathcal{C}}\}$\;
    $T \leftarrow \emptyset$\;
}
\uElseIf{$\rho = 0$ and $\mathcal{C}$ is not smooth at $x$}{
    $\tilde{V} \leftarrow \left\{\tilde{v}_{\mathrm{ex}}\right\}$\;
    \Return{$\tilde{V}$} 
}
\Else(\tcp*[f]{smooth boundary or shallow interior}){
    $\hat{x} \in \Pi(x)$ \tcp*{$\Pi(x)=\{x\}$ when $\rho=0$}
    $u \leftarrow$ inward unit normal to $\mathcal{C}$ at $\hat{x}$, $u\in\aff{\mathcal{C}-x}$\;
    Extend to $W = \{u, w_1, \ldots, w_{\dim \mathcal{C} - 1}\}$, an orthonormal basis of $\aff{\mathcal{C} - x}$\;
    $U \leftarrow \{u, -u, w_1, -w_1, \ldots, w_{\dim \mathcal{C} - 1}, -w_{\dim \mathcal{C} - 1}\} \cap \ri{D_{\mathcal{C}}(x)}$\;
    $T \leftarrow \{w_1, -w_1, \ldots, w_{\dim \mathcal{C} - 1}, -w_{\dim \mathcal{C} - 1}\}$\;
}

\tcp{\textbf{Step 3: Tilt and noise}}
Enumerate the set $U \cup T = \{v_1, \ldots, v_m\}$\;
Sample $\delta \sim \unif{\Delta_m}$ and set $\zeta \leftarrow \sum_{i=1}^{m} \delta_i v_i$\;
\For{$v \in T$}{
    $\tilde{V} \leftarrow \tilde{V}\cup \bigl\{\mathrm{normalize}\bigl((1-\alpha)\bar{\eta}\,u + \alpha\bar{\eta}\,\zeta + (1-\bar{\eta})\,v\bigr)\bigr\}$\;
}
\For{$v \in U$}{
    $\tilde{V} \leftarrow \tilde{V}\cup \bigl\{\mathrm{normalize}\bigl(\alpha\bar{\eta}\,\zeta + (1-\alpha\bar{\eta})\,v\bigr)\bigr\}$\;
}
$\tilde{V} \leftarrow \tilde{V}\cup \{\tilde{v}_{\mathrm{ex}}\}$\;
\Return{$\tilde{V}$}
\end{algorithm}

\begin{algorithm}[htp]
\caption{The Noised Regularized Feasible Directions Method}\label{alg:GRFD}
\KwIn{$x_0 \in \mathcal{C}$, $\gamma \in (0,\infty)$, $r > 0$, $\eta \in (0, 1]$, $\alpha \in [0,1)$, $\rho_0 \ge 0$}

\For{$k = 0, \ldots, K$}{

    $\tilde{V}^k \leftarrow \textsc{ConstructSearchDirections}(x_k,\, \eta,\, \alpha,\, \rho_0)$ \quad (\cref{alg:construct_directions})\;

    \For{$\tilde{v} \in \tilde{V}^k$}{
        $\tau(x_k;\tilde{v}) \leftarrow \sup\{t \ge 0: x_k + t\tilde{v} \in \mathcal{C}\}$\;
        $t_{\tilde{v}} \in \arg\min_{t \in [0,\, \min\{\tau(x_k;\tilde{v}), r\}]} h(x_k + t\tilde{v}) + \frac{\gamma}{2} t^2$\;
    }
    $\tilde{v}^* \in \arg\min_{\tilde{v} \in \tilde{V}^k} h(x_k + t_{\tilde{v}}\, \tilde{v})$\;
    $x_{k+1} \leftarrow x_k + t_{\tilde{v}^*}\, \tilde{v}^*$\;
}
\end{algorithm}
\ULforem

\begin{remark}[Access at $\rho_0>0$]\label{rem:rho0-default}
Choosing $\rho_0>0$ requires access to $\rho(x)$ and, when
$0<\rho(x)\le\rho_0$, a nearest boundary point
$\hat{x}\in\Pi(x)$ and its inward unit normal.
These quantities can be computed explicitly for Euclidean balls
and full-dimensional polyhedra. For ellipsoids, the nearest-point
problem reduces to quadratic minimization on a sphere
\citep{martinez1994local}.

For general convex sets, this geometric access is an additional
assumption; membership queries alone do not provide it exactly.
When it is unavailable, one may choose $\rho_0=0$, using only
the queries listed at the end of
\cref{sec:directions_oracle_and_alg}.
If $\rb{\mathcal{C}}=\emptyset$, we set $\rho(x):=+\infty$,
so the first branch of \cref{alg:construct_directions} applies
and no boundary query is needed.
\end{remark}

\section{Analysis}\label{sec:analysis}
\label{sec:convergence_results}

We establish PDM convergence rates and almost-sure $d$-stationarity
of accumulation points for \cref{alg:Oracle_GRFD,alg:GRFD}.
We also derive finite-time bounds for approximate $d$-stationarity
and specialize them under curvature assumptions.
The analysis relies on the feasibility, spanning, and
boundary-avoidance results proved in \cref{sec:noised_alg_properties}.

%

\subsection{The PDM measure and a finite-direction rate}\label{sec:convergence_rate_result}
The analysis of \cref{alg:Oracle_GRFD} and \cref{alg:GRFD} is based on a feasible-directions optimality measure that extends the one introduced by \cite{BECK2022517}. 
We refer to this quantity as the \textit{proximal directional minimization (PDM) measure}.
The measure quantifies the penalized decrease available along a finite set of feasible directions, and is therefore suited to methods that explore finitely many directions per iteration.

\medskip

For any feasible $\bar{x}\in \mathcal{C}$, direction set $V=\{v_i\}_{i=1}^s$, penalty parameter $\gamma\ge 0$, and interval radius $r\in(0,\infty]$, the PDM is defined as
\begin{equation}\label{eq:optimality_measure}
    \mathcal{O}^\gamma(\bar{x};V,r)
    :=\sum_{i=1}^s \bigl(\mu^\gamma(\bar{x};v_i,r)\bigr)^2,
\end{equation}
where the \textit{$v$-proximal directional minimizer} scalars are defined by
\begin{equation}\label{eq:PDM_minimizer}
    \mu^\gamma(\bar{x}; v, r) \in \argmin_{0 \le t \le \min\{\tau(\bar{x};v), r\}} \left\{h (\bar{x} + t v) + \dfrac{\gamma}{2}t^2\right\},
\end{equation}
and $\tau(\bar{x};v):=\sup\{t\ge 0:\bar{x}+tv\in\mathcal{C}\}$ is the maximal feasible step length at $\bar{x}$ in direction $v$. 
When multiple minimizers exist, we select the smallest one.

The quantity $\mu^\gamma(\bar{x};v,r)$ represents the step size that achieves the best penalized decrease along direction $v$, while $\mathcal{O}^\gamma(\bar{x};V,r)$ aggregates the descent attainable over the direction set $V$. 
In particular, $\mathcal{O}^\gamma(\bar{x};V,r)=0$ if and only if no penalized decrease is available along any direction in $V$; \cref{lem:optimality_measure_zero_implies_no_descent_in_span} below relates this to first-order conditions.


%

The following lemma is adapted from \citep[Lemma 4.1]{BECK2022517}.
\begin{lem}\label{lem:optimality_measure_zero_implies_no_descent_in_span}
    Let $\bar{x} \in \mathcal{C}$, $r\in (0,\infty]$,  $V\subseteq D_{\mathcal{C}}(\bar{x})$ be finite, and $\mathcal{K}$ be a cone such that the positive span of $V$ is a superset of $\mathcal{K}$. 
    If $\mathcal{O}^\gamma(\bar{x};V,r) = 0$, then
    \begin{equation*}
        \langle \nabla f(\bar{x}),d \rangle - g'(\bar{x};d) \ge 0, \qquad \forall d \in \mathcal{K}.
    \end{equation*}
\end{lem}
The proof is in \cref{sec:appendix_proofs_analysis}.

\begin{lem}[properties of $h'(\bar{x};\cdot)$]\label{lem:h-prime-lipschitz-span}
Let $\mathcal{C}$ be closed and convex, let $f$ be continuously differentiable, let $g$ be convex, and let $\bar{x}\in\mathcal{C}$. Then:
\begin{enumerate}
    \item If $\mathcal{C}$ is compact and $g$ is finite near $\mathcal{C}$, then
    \begin{equation*}
        \Lambda\;:=\;\sup_{x\in\mathcal{C}}\Bigl(\norm{\nabla f(x)}+\max_{s\in\partial g(x)}\norm{s}\Bigr)
    \end{equation*}
    is finite, $d\mapsto h'(\bar{x};d)$ is $\Lambda$-Lipschitz, and $|h'(\bar{x};d)|\le\Lambda\norm{d}$.
    \item If $V\subseteq D_{\mathcal{C}}(\bar{x})$ is finite and
    $d=\sum_{v\in V}\lambda_v v$ with $\lambda_v\ge 0$, then
    $-h'(\bar{x};d)\le\sum_{v\in V}\lambda_v[-h'(\bar{x};v)]_+$.
\end{enumerate}
\end{lem}

The proof is in \cref{sec:appendix_proofs_analysis}.
Part~1 is the Lipschitz refinement of the continuity recorded in \cref{lem:h_prime_properties} and supplies the constant $\Lambda$ used below.
Part~2 is the quantitative form of \cref{lem:optimality_measure_zero_implies_no_descent_in_span}: taking all $[-h'(\bar{x}; v)]_+=0$ recovers it.

The PDM measure is nonlocal: it requires the absence of penalized descent anywhere along
the feasible segment $[0,\min\{\tau(\bar{x};v),r\}]$ for every $v\in V$, and is therefore
strictly stronger than the first-order condition $h'(\bar{x};v)\ge 0$.
Even when the directional derivative is nonnegative, a step $t>0$ satisfying $h(\bar{x}+tv)+\frac{\gamma}{2}t^2<h(\bar{x})$ may still exist, and the line search then returns a positive step.

Under additional smoothness assumptions, the following proposition provides a quantitative lower bound on the PDM measure in terms of directional derivatives and feasible step lengths; see \cite[Theorem~4.1]{BECK2022517}.


\begin{proposition}[Lower bound on the PDM measure; adapted from
            {\cite[Theorem~4.1]{BECK2022517}}]
\label{prop:PDM_lower_bound}
    Suppose $g\equiv 0$ and $f\in C^{1,1}_L$ for some $L>0$.  Then for
    any $\bar{x}\in\mathcal{C}$ and finite direction set $V$ with $\norm{v}=1$ for all $v\in V$,
    \begin{equation}\label{eq:PDM_lower_bound}
        \mathcal{O}^\gamma(\bar{x};V,r)
        \;\ge\;
        \frac{1}{(L+\gamma)^2}
        \sum_{v\in V}
        \Bigl(\min\bigl\{[-\langle\nabla f(\bar{x}),v\rangle]_+,\;
                         (L+\gamma)\,\min\{\tau(\bar{x};v), r\}\bigr\}\Bigr)^2.
    \end{equation}
\end{proposition}
The proof is in \cref{sec:appendix_proofs_analysis}.
The bound shows that the PDM measure controls the magnitude of descent available along directions in $V$, as quantified by the projected gradient $[-\inner{\nabla f(\bar{x})}{v}]_+$ truncated by the available step $\min\{\tau(\bar{x};v),r\}$.

We illustrate the significance of the PDM measure and
\cref{prop:PDM_lower_bound} through three examples.


\begin{example}[Unconstrained smooth minimization]
\label{ex:PDM_unconstrained}
    Let $\mathcal{C}=\mathbb{R}^n$, $g\equiv 0$, $f\in C^{1,1}_L$, $r=\infty$, and
    take $V=\{\pm e_1,\ldots,\pm e_n\}$.  Since $\tau(\bar{x};v)=\infty$ for
    every $v$, the minimum in \eqref{eq:PDM_lower_bound} reduces to
    $[-\langle\nabla f(\bar{x}),v\rangle]_+$ for each direction.  Summing over
    $V$ then yields
    \[
        \mathcal{O}^\gamma(\bar{x};V,r)
        \;\ge\;
        \frac{1}{(L+\gamma)^2}\,\|\nabla f(\bar{x})\|^2,
    \]
    so the PDM measure provides a lower bound proportional to the standard
    gradient-norm optimality measure, and $\mathcal{O}^\gamma(\bar{x};V,r)=0$
    implies $\nabla f(\bar{x})=0$.
\end{example}

\begin{example}[Interior point]
\label{ex:PDM_interior}
    Suppose $\bar{x}\in\ri{\mathcal{C}}$, let $\{w_1,\ldots,w_{\dim\mathcal{C}}\}$
    be an orthonormal basis of $\aff{\mathcal{C}-\bar{x}}$, and set
    $V = \{w_1,-w_1,\ldots,w_{\dim\mathcal{C}},-w_{\dim\mathcal{C}}\}$.
    By \cref{lem:feasible_directions_at_relative_interior},
    $D_\mathcal{C}(\bar{x})=\aff{\mathcal{C}-\bar{x}}$, and by construction
    $\pos{V}=\aff{\mathcal{C}-\bar{x}}$ as well.
    Therefore \cref{lem:optimality_measure_zero_implies_no_descent_in_span}
    applies to all feasible directions: $\mathcal{O}^\gamma(\bar{x};V,r)=0$
    implies $\bar{x}$ is $d$-stationary, without any smoothness assumption on
    $f$ or $g$.

    Under the additional assumptions $g\equiv 0$ and $f\in C^{1,1}_L$,
    \cref{prop:PDM_lower_bound} yields a complementary interpretation.
    Since $\bar{x}\in\ri{\mathcal{C}}$, every $v\in V$ satisfies $\min\{\tau(\bar{x};v),r\}>0$,
    so each term in the lower bound \eqref{eq:PDM_lower_bound} vanishes
    only if $\langle\nabla f(\bar{x}),v\rangle\ge 0$. Since $V$ contains both
    $v$ and $-v$ for every basis vector, this forces
    $\langle\nabla f(\bar{x}),v\rangle=0$ for all $v\in V$, and hence the
    projection of $\nabla f(\bar{x})$ onto $\aff{\mathcal{C}-\bar{x}}$ is zero -- a
    necessary optimality condition for smooth minimization over
    $\mathcal{C}$ at a relative interior point.
\end{example}

\begin{example}[Boundary point with curved nonlinear constraints]
\label{ex:PDM_boundary}
Suppose $\mathcal{C}$ is compact, fix $\eta\in(0,\tfrac12]$, and let $\hat{x}\in\rb{\mathcal{C}}$ be a point at which $\mathcal{C}$ is smooth, with inward unit normal $u$.
Let $V$ positively span $D^\eta_\mathcal{C}(\hat{x})$ but not necessarily the full cone $D_\mathcal{C}(\hat{x})$; \cref{lem:tilde_v_spans_eta_feasible} shows that the non-exploration output of \cref{alg:construct_directions} with $\alpha=0$ has this property.
Let $\Lambda$ be the constant of \cref{lem:h-prime-lipschitz-span}(1), which is finite since $\mathcal{C}$ is compact and $g$ is finite near $\mathcal{C}$; then $d\mapsto h'(x;\,d)$ is $\Lambda$-Lipschitz for every $x\in\mathcal{C}$.
Let $d\in D_\mathcal{C}(\hat{x})$ be a unit feasible direction that is not an $\eta$-FD.
Then $0\le\inner{d}{u}<\eta$ by \cref{lem:eta_FD_properties}(b). Set $s:=\eta/(1-\eta)$. Then $\norm{d+su}^2=1+2s\inner{d}{u}+s^2\ge 1+s^2>0$, and the tilted direction $w:=(d+su)/\norm{d+su}$ satisfies $\inner{w}{u}\ge\eta$, so $w\in D^\eta_\mathcal{C}(\hat{x})$, while $\norm{w-d}\le 2s\le 4\eta$.
Hence
\begin{equation*}
    -h'(\hat{x};d)\;\le\;-h'(\hat{x};w)+4\Lambda\eta,
\end{equation*}
so the descent rate available outside the $\eta$-FD cone exceeds that available inside it by at most $4\Lambda\eta$.
Consequently the PDM measure, constructed using $\eta$-FDs, captures all descent directions up to an $O(\eta)$ error in the directional derivative.
\end{example}

\begin{remark}[Scope of the PDM certificate]\label{rem:pdm-vs-dstat}
By \cref{lem:optimality_measure_zero_implies_no_descent_in_span},
$\mathcal{O}^\gamma(\bar{x};V,r)=0$ certifies $h'(\bar{x};\,d)\ge 0$ only on $\pos{V}$.
\Cref{tab:pdm-scope} lists, for each method and type of point, a cone contained in $\pos{V}$, where $V$ is the set of non-exploration search directions.
At relative-interior points the certificate is full $d$-stationarity; at smooth boundary points, \cref{ex:PDM_boundary} bounds the descent missed outside the $\eta$-FD cone by $O(\eta)$.
\end{remark}

\begin{table}[ht]
\centering
\footnotesize
\caption{Cones on which vanishing PDM certifies $h'(\bar{x};\,\cdot)\ge 0$. For \cref{alg:GRFD} with $\alpha>0$, ``small'' refers to the pointwise thresholds in \cref{cor:noised-pdm-certificate}; the last entry assumes $x_0$ is almost surely not a nonsmooth boundary point.}
\label{tab:pdm-scope}
\begin{tabular}{@{}llll@{}}
\toprule
 & \cref{alg:Oracle_GRFD} & \cref{alg:GRFD}, $\alpha=0$ & \cref{alg:GRFD}, small $\alpha>0$ \\
\midrule
$\bar{x}\in\ri{\mathcal{C}}$
 & $D_{\mathcal{C}}(\bar{x})$
 & $D_{\mathcal{C}}(\bar{x})$
 & $D_{\mathcal{C}}(\bar{x})$ \\
 & {\footnotesize \cref{def:eta_directions_oracle}, \cref{lem:eta_FD_properties}(a)}
 & {\footnotesize \cref{lem:tilde_v_spans_eta_feasible}}
 & {\footnotesize \cref{cor:noised-pdm-certificate}(1)} \\[2pt]
smooth boundary
 & $D^\eta_{\mathcal{C}}(\bar{x})$
 & $D^\eta_{\mathcal{C}}(\bar{x})$
 & $D^{\beta\eta}_{\mathcal{C}}(\bar{x})$, $\beta\in(1,1/\eta)$ \\
 & {\footnotesize \cref{def:eta_directions_oracle}}
 & {\footnotesize \cref{lem:tilde_v_spans_eta_feasible}}
 & {\footnotesize \cref{cor:noised-pdm-certificate}(2)} \\[2pt]
nonsmooth boundary
 & $D^\eta_{\mathcal{C}}(\bar{x})$
 & none
 & not reached a.s. \\
 & {\footnotesize \cref{def:eta_directions_oracle}}
 & {\footnotesize \cref{rem:nondegenerate-direction-sets}}
 & {\footnotesize \cref{lem:always_smooth_almost_surely}} \\
\bottomrule
\end{tabular}
\end{table}

Having related the PDM measure to finite-direction stationarity, we now establish a rate.
\Cref{thm:convergence_rate} shows that the iterates approach finite-direction stationarity at an $O(1/K)$ rate in the PDM measure, counted in outer iterations. Classical direct-search guarantees have the same order in a different unit, function evaluations; see the comparison after \cref{cor:eps-complexity}.
Write $s_k := |\tilde{V}^k|$ for the number of search directions at iteration $k$, as in \cref{alg:Oracle_GRFD}, and set $s_{\max}^K := \max_{0 \le i \le K-1} s_i$. When $K$ is clear from the context we write $s_{\max}$ for $s_{\max}^K$.

\begin{theorem}\label{thm:convergence_rate}
    Let $\left\{x_k\right\}_{k \ge 0}$ be a sequence generated by \cref{alg:Oracle_GRFD} or \cref{alg:GRFD}.
    Then  for any $K\geq 1$ it holds that
    \begin{equation*}
        \min\limits_{0 \le k \le K-1} \mathcal{O}^\gamma (x_k; \tilde{V}^k, r) \le \dfrac{2 s_{\max}}{\gamma K}\left[\left(f(x_0) - g(x_0)\right) - \inf\limits_{x \in \mathcal{C}} \left\{f(x) - g(x)\right\}\right].
    \end{equation*}
\end{theorem}
\begin{proof}
    Enumerate $\tilde{V}^k = \left\{\tilde{v}_1,\ldots,\tilde{v}_{s_k}\right\}$.
    Then by the definition of $t_i^k$, for all $i \in \left\{1,\ldots,s_k \right\}$,
    \begin{equation*}
        \left(f(x_k) - g(x_k)\right) - \left(f(x_k+t_i^k \tilde{v}_i) - g(x_k + t_i^k \tilde{v}_i)\right) \ge \dfrac{\gamma}{2}{t_i^k}^2 = \dfrac{\gamma}{2} (\mu^\gamma(x_k; \tilde{v}_i, r))^2.
    \end{equation*}
    Summing over $i$, 
    \begin{equation*}
        \sum\limits_{i=1}^{s_k}\left[\left(f(x_k) - g(x_k)\right) - \left(f(x_k+t_i^k \tilde{v}_i) - g(x_k + t_i^k \tilde{v}_i)\right)\right] \ge \dfrac{\gamma}{2}\mathcal{O}^\gamma(x_k;\tilde{V}^k, r).
    \end{equation*}
    Using the fact that $i^*_k$ maximizes the difference between the function values
    $$\left(f(x_k) - g(x_k)\right) - \left(f(x_k+t_i^k \tilde{v}_i) - g(x_k + t_i^k \tilde{v}_i)\right),$$ 
    the fact that $s_{max} \ge s_k$, and that $x_{k+1} = x_k + t_{i^*_k}^k \tilde{v}_{i^*_k}$, we deduce that
    \begin{equation}\label{eq:convergence_rate:optimality_measure_per_step_bound}
        s_{max} \cdot  \left[\left(f(x_k) - g(x_k)\right) - \left( f(x_{k+1}) - g(x_{k+1})\right)\right] \ge \dfrac{\gamma}{2}\mathcal{O}^\gamma(x_k;\tilde{V}^k, r).
    \end{equation}
    Summing over $k=0$ to $k=K-1$ and dividing by $K$ finally yields that 
    \begin{align*}
        \dfrac{2s_{max}}{\gamma K}\left[\left(f(x_0) - g(x_0)\right) - \inf\limits_{x \in \mathcal{C}} \left\{f(x) - g(x)\right\}\right] &\ge \dfrac{1}{K} \sum_{k=0}^{K-1} \mathcal{O}^\gamma(x_k;\tilde{V}^k, r) \\
        &\ge \min\limits_{0 \le k \le K-1}\mathcal{O}^\gamma(x_k;\tilde{V}^k, r).
    \end{align*}
\end{proof}
The rate bound requires neither directional spanning nor smoothness
of $\mathcal{C}$. Spanning determines the scope of the non-descent
certificate in
\cref{lem:optimality_measure_zero_implies_no_descent_in_span}
(see \cref{rem:pdm-vs-dstat}).
When $\tilde{V}^k$ contains only the exploration direction, the bound
concerns that direction alone; \cref{rem:nondegenerate-direction-sets}
gives conditions excluding this case.
We next bound the direction-set size $s_{\max}$ appearing in
\cref{thm:convergence_rate}.

\begin{remark}[Bounds on $s_{\max}$]\label{rem:smax}
Let $m=\dim\mathcal{C}$. In \cref{alg:GRFD}, counting the
exploration direction gives
\[
s_k\le
\begin{cases}
2m+1, & \rho(x_k)>\rho_0,\\
1, & x_k \text{ is a nonsmooth boundary point},\\
2m, & x_k \text{ is a smooth boundary point},\\
4m-1, & 0<\rho(x_k)\le\rho_0.
\end{cases}
\]
The last bound counts both the tilted and untilted directions
generated from $T\subseteq U$
(\cref{rem:interior-gap-shallow}). Consequently,
$s_{\max}\le 2m+1$ when $\rho_0=0$, and
$s_{\max}\le 4m-1$ when $\rho_0>0$.

In \cref{alg:Oracle_GRFD}, $s_k\le |V^k|+1$, so the bound
depends on the oracle implementation. For the constructions
in \cref{sec:directions_oracle_and_alg}, the polyhedral case
gives one plus the maximum number of cone generators returned
by the oracle. For smooth convex sets, the zero-noise
construction of \cref{alg:construct_directions} gives the same
bounds as above. For orthogonally structured constraints,
\eqref{eq:spanning_directions_for_orthogonal_constraint_gradients}
gives $s_k\le 2n+2-|\mathcal{I}(x_k)|\le 2n+2$.
\end{remark}

\subsection{Subsequence convergence to $d$-stationarity}\label{sec:asymptotic_convergence_result}

The PDM rate in \cref{thm:convergence_rate} concerns the directions
searched at each iteration. At boundary points, their positive span
may omit feasible directions, so this bound alone does not establish
$d$-stationarity.
\Cref{thm:main} shows that randomized exploration ensures, almost
surely, that every accumulation point is $d$-stationary. The proof
combines two observations: objective decrease forces the exploration
step lengths to vanish, while repeated sampling near any
nonstationary accumulation point would produce exploration steps
bounded away from zero.



\begin{theorem}[Subsequence convergence; adapted from {\citep[Theorem~4.1]{greenstein2026achieving}}]%
\label{thm:main}
    Let $\{x_k\}_{k=0}^\infty$ be a sequence generated by \Cref{alg:Oracle_GRFD} or \Cref{alg:GRFD}.
    Then, almost surely, every accumulation point of $\{x_k\}$ is a $d$-stationary point of \eqref{eq:feasible_directions_beck_hallak}.
\end{theorem}

\Cref{thm:main} is, in substance, \citep[Theorem~4.1(ii)]{greenstein2026achieving}. There, an arbitrary method producing feasible iterates is augmented with an exploration step whose step size is sampled at random and accepted upon sufficient decrease. \Cref{alg:Oracle_GRFD} and \Cref{alg:GRFD} have the same structure, with the non-exploration line searches in the role of the base method, except that the exploration step size is obtained from the regularized line search of \cref{rem:line-search-access}. This affects only the final step of the argument. In \citep{greenstein2026achieving}, a sampled step size bounded away from zero is accepted infinitely often, contradicting the summability of the squared exploration steps; here, such a step size is admissible in the line search, and the optimality of the line-search minimizer yields the contradiction (\cref{lem:non_descent}). The remaining ingredients are those of \citep{greenstein2026achieving}, with minor modifications. For completeness, \cref{sec:appendix_proof_thm_main} gives the full proof and indicates the source of each ingredient.

\begin{remark}[Existence of accumulation points; cf.\ {\citep[Remark~4.1]{greenstein2026achieving}}]\label{rem:accumulation-points}
A sufficient condition for the existence of an accumulation point of $\{x_k\}$ is boundedness of the level set $\{x \in \mathcal{C} : h(x) \le h(x_0)\}$.
Indeed, $t=0$ is admissible in each regularized line search, so every candidate satisfies
\begin{equation*}
    h(x_k + t_i^k \tilde{v}_i^k)
    \le h(x_k + t_i^k \tilde{v}_i^k) + \tfrac{\gamma}{2}(t_i^k)^2
    \le h(x_k),
\end{equation*}
and therefore $h(x_{k+1}) \le h(x_k)$ for all $k$, independently of which candidate is selected. Hence $h(x_k) \le h(x_0)$ and $\{x_k\}$ remains in this set. If the set is bounded, then $\{x_k\}$ is bounded in $\real^n$ and admits an accumulation point by the Bolzano--Weierstrass theorem.
\end{remark}

\subsection{Finite-time approximate directional stationarity}\label{sec:finite-time-complexity}

\Cref{thm:convergence_rate} bounds the PDM measure on the directions actually searched. We now bound the
directional derivative over all feasible directions. By \cref{prop:PDM_lower_bound}, when $g\equiv 0$ a small PDM
value controls $[-\inner{\nabla f(\bar{x})}{v}]_+$ only up to truncation by the feasible step $\min\{\tau(\bar{x}; v), r\}$, and
near $\rb{\mathcal{C}}$ that step can be arbitrarily small.

The truncation is not an artifact of the bound. Let $u$ be a unit vector and
$\mathcal{C}:=\{x:\inner{x}{u}\ge 0\}$. Let $f(x):=\inner{u}{x}$, $g\equiv 0$ and $\bar{x}:=\rho u$ with $\rho>0$. Then $h'(\bar{x};d)=\inner{d}{u}$,
and the steepest unit feasible direction is $-u$, with rate $1$. A unit $d$ with $c:=\inner{d}{u}<0$ has
$\tau(\bar{x};d)=\rho/|c|$. Every feasible step along a descent direction therefore decreases $h$ by at most
$|c|\,\tau(\bar{x};d)=\rho$.

As $\rho\downarrow 0$ the decrease available to any method vanishes while the steepest rate stays $1$. Hence
per-iteration decrease cannot certify a small directional derivative over all of $D_{\mathcal{C}}(\bar{x})$, and we
discard the directions along which one cannot travel far.

\begin{defin}[$R$-deep residual; $(\varepsilon,R)$-d-stationarity]\label{def:R-deep-residual}
Let $\mathcal{C}$ be closed and convex, let $f$ be differentiable, and let $g$ be convex.
For $\bar{x}\in\mathcal{C}$ and $R\ge 0$ let
\begin{equation}\label{eq:R-deep-residual}
    \Delta^{[R]}(\bar{x})\;:=\;\sup\bigl\{[-h'(\bar{x};d)]_+\;:\;d\in D_{\mathcal{C}}(\bar{x}),\ \norm{d}=1,\ \tau(\bar{x};d)\ge R\bigr\},
\end{equation}
with $\sup\emptyset:=0$, and write $\Delta:=\Delta^{[0]}$.
We call $\bar{x}$ $(\varepsilon,R)$-$d$-stationary if $\Delta^{[R]}(\bar{x})\le\varepsilon$.
\end{defin}

In the example above, $\Delta^{[R]}$ keeps only descent directions with $\rho/|c|\ge R$, so
$\Delta^{[R]}(\bar{x})\le \rho/R$. For fixed $R$ this is of the same order as the available decrease, whereas
$\Delta(\bar{x})=1$. More generally, when $\mathcal{C}$ is compact and $g$ is finite near $\mathcal{C}$, a direction with $\tau(\bar{x};d)<R$ offers at most $\Lambda R$ of
first-order decrease by \cref{lem:h-prime-lipschitz-span}(1). The threshold $R$ is on the feasible step, not on the depth
$\rho(\bar{x})$.

$\Delta(\bar{x})=0$ if and only if $\bar{x}$ is $d$-stationary. Since $h'(\bar{x};\,\cdot)$ is continuous, the supremum over
$D_{\mathcal{C}}(\bar{x})$ agrees with the one over $\cl{D_{\mathcal{C}}(\bar{x})}$, so $\Delta$ measures $d$-stationarity in the strong sense of
\cref{thm:main}. The map $R\mapsto\Delta^{[R]}(\bar{x})$ is nonincreasing and increases to $\Delta(\bar{x})$ as $R\downarrow 0$. If every unit
feasible direction at $\bar{x}$ admits a step of length at least $R$, then $\Delta^{[R]}(\bar{x})=\Delta(\bar{x})$. In
particular this holds at every point when $\mathcal{C}=\real^n$, and \cref{thm:finite-time-rate} is then a rate for $\Delta$.

We write $\tau_r(x;v):=\min\{\tau(x;v),r\}$ and $H:=h(x_0)-\inf_{\mathcal{C}}h$. The decrease along a single
direction is quantified by the following analogue of \cref{prop:PDM_lower_bound} for the objective of \eqref{eq:feasible_directions_beck_hallak}.

\begin{lem}\label{lem:per-direction-decrease}
Let $\mathcal{C}$ be closed and convex, let $f\in C^{1,1}_L$, let $g$ be convex, and let $\gamma\ge 0$.
Let $\bar{x}\in\mathcal{C}$, let $v\in D_{\mathcal{C}}(\bar{x})$ with $\norm{v}=1$, set $\theta:=[-h'(\bar{x};v)]_+$. Then
\begin{equation*}
    h(\bar{x})-\min_{0\le t\le \tau_r(\bar{x};v)}\Bigl\{h(\bar{x}+tv)+\tfrac{\gamma}{2}t^2\Bigr\}\;\ge\;\tfrac12\min\Bigl\{\tfrac{\theta^2}{L+\gamma},\,\theta \tau_r(\bar{x};v)\Bigr\}.
\end{equation*}
\end{lem}

The proof is in \cref{sec:appendix_proofs_analysis}. As in the proof of \cref{thm:convergence_rate}, the left-hand side is the regularized
decrease along $v$, and is comparable to $\bigl(\mu^\gamma(\bar{x};v,r)\bigr)^2$.

The first term of the minimum yields an $\varepsilon^{-2}$ rate, but it is the smaller term only when
$\tau_r(\bar{x};v)\ge\theta/(L+\gamma)$. To transfer decrease to every $R$-deep direction we therefore write each
such direction as a nonnegative combination of searched directions with controlled coefficients and controlled feasible-step lengths.

\begin{defin}[$\eta$-direction oracle with bounds on the spanning coefficients]\label{def:quantitative-eta-oracle}
Let $\mathcal{C}$ be closed and convex and let $r\in(0,\infty]$.
A mapping $x\mapsto V(x)\subseteq D_{\mathcal{C}}(x)$ with $\norm{v}=1$ for all $v\in V(x)$ has spanning-coefficient bounds $(A,B,R)$ if for every $x\in\mathcal{C}$ and every $d\in D_{\mathcal{C}}(x)$ with $\norm{d}=1$ and $\tau(x;d)\ge R$ there exist $\lambda_v\ge 0$, $v\in V(x)$, with
\begin{equation}\label{eq:quantitative-oracle}
    d=\sum_{v}\lambda_v v,\qquad\sum_v\lambda_v\le A,\qquad\sum_v\frac{\lambda_v}{\tau_r(x;v)}\le B.
\end{equation}
\end{defin}

The first constant is the usual $\ell_1$ positive-spanning constant; the second weights each coefficient by the reciprocal of the guaranteed step length along its direction, and the parameter $\eta$ enters the complexity only through $B$.

\begin{theorem}[finite-time rate]\label{thm:finite-time-rate}
Let $\mathcal{C}$ be closed and convex, let $f\in C^{1,1}_L$, let $g$ be convex, and assume $H<\infty$.
Let $\gamma>0$ and $r\in(0,\infty]$.
Let $x\mapsto V(x)$ have spanning-coefficient bounds $(A,B,R)$ and let $\{x_k\}$ be generated by any method that, at each iteration, solves the regularized line search along every direction of a finite set containing $V(x_k)$ and moves to the candidate with the smallest value of $h$.
Then for every $K\ge 1$,
\begin{equation}\label{eq:finite-time-rate}
    \min_{0\le k\le K-1}\Delta^{[R]}(x_k)\;\le\;A\sqrt{\frac{2(L+\gamma)H}{K}}\;+\;\frac{2BH}{K}.
\end{equation}
\end{theorem}

\begin{proof}
Write $q_k:=h(x_k)-h(x_{k+1})$.
Since the method sets $t_i^k\in\arg\min_{t\in[0,\tau_r(x_k;v_i)]}\{h(x_k+tv_i)+\tfrac{\gamma}{2}t^2\}$ along each searched direction $v_i$ and then picks the candidate with the smallest $h$-value,
\begin{equation*}
    q_k
    \;=\;
    \max_{i}\bigl[h(x_k)-h(x_k+t_i^k v_i)\bigr]
    \;\ge\;
    \max_{i}\Bigl[h(x_k)-\min_{t\in[0,\tau_r(x_k;v_i)]}\bigl\{h(x_k+tv_i)+\tfrac{\gamma}{2}t^2\bigr\}\Bigr],
\end{equation*}
using $\tfrac{\gamma}{2}(t_i^k)^2\ge 0$.
Combining with \cref{lem:per-direction-decrease}, for every $v\in V(x_k)$,
\begin{equation*}
    \min\Bigl\{\tfrac{\theta_v^2}{L+\gamma},\,\theta_v \tau_r(x_k;v)\Bigr\}\;\le\;2q_k,
    \qquad
    \theta_v:=[-h'(x_k;v)]_+.
\end{equation*}
Hence, examining which branch attains the minimum,
\begin{equation*}
    \theta_v
    \;\le\;
    \max\Bigl\{\sqrt{2(L+\gamma)q_k},\,\frac{2q_k}{\tau_r(x_k;v)}\Bigr\}
    \;\le\;
    \sqrt{2(L+\gamma)q_k}+\frac{2q_k}{\tau_r(x_k;v)}.
\end{equation*}
Now let $d$ be a unit feasible direction at $x_k$ with $\tau(x_k;d)\ge R$, and take the decomposition guaranteed by \eqref{eq:quantitative-oracle}.
By \cref{lem:h-prime-lipschitz-span}(2) and the bound on $\theta_v$,
\begin{equation*}
    -h'(x_k;d)
    \;\le\;
    \sum_v\lambda_v\theta_v
    \;\le\;
    \sqrt{2(L+\gamma)q_k}\sum_v\lambda_v
    +
    2q_k\sum_v\frac{\lambda_v}{\tau_r(x_k;v)}
    \;\le\;
    A\sqrt{2(L+\gamma)q_k}+2Bq_k.
\end{equation*}
Taking the supremum over such $d$ yields $\Delta^{[R]}(x_k)\le A\sqrt{2(L+\gamma)q_k}+2Bq_k$.
Finally $\sum_{k=0}^{K-1}q_k=h(x_0)-h(x_K)\le H$, so $\min_k q_k\le H/K$; evaluating the right-hand side at $k^\star=\arg\min_k q_k$ gives \eqref{eq:finite-time-rate}.
\end{proof}

\Cref{thm:finite-time-rate} controls $\Delta^{[R]}$ over all unit
feasible directions with $\tau(x;d)\ge R$. Its proof uses
spanning-coefficient bounds, is deterministic, and yields a bound
independent of $s_{\max}$. The proof also yields an early-stopping
rule. It gives $\Delta^{[R]}(x_k)\le A\sqrt{2(L+\gamma)q_k}+2Bq_k$
with $q_k:=h(x_k)-h(x_{k+1})$, so when $A$, $B$ and $L$ are known
(e.g., via \cref{prop:quantitative-spanning}) one may stop at the
first $k$ at which this bound is at most $\varepsilon$. This certifies
$\Delta^{[R]}(x_k)\le\varepsilon$, and since $\min_{k<K}q_k\le H/K$,
it occurs before any $K$ for which the right-hand side of
\eqref{eq:finite-time-rate} is at most $\varepsilon$.

We next establish the required oracle properties for the
zero-noise construction in \cref{alg:construct_directions}
under curvature bounds. We begin with rolling-ball geometry
and the resulting estimates for $\tau(\bar{x};d)$.
The geometric proofs are deferred to
\cref{sec:appendix_proofs_analysis}.

\begin{lem}[Rolling balls]\label{lem:rolling_balls}
    Let $\mathcal{C}\subseteq\mathbb{R}^n$ be a convex body with a $C^2$
    boundary whose principal curvatures are uniformly bounded in
    $[\kappa_{\min},\kappa_{\max}]$ at every boundary point, where
    $0<\kappa_{\min}\le\kappa_{\max}<\infty$.
    Write $r_{\mathrm{in}}:=1/\kappa_{\max}$ and $r_{\mathrm{out}}:=1/\kappa_{\min}$.
    Then:
    \begin{enumerate}
        \item At every $\hat{x}\in\rb{\mathcal{C}}$ with inward unit normal $u$,
              $\mathcal{B}(\hat{x}+r_{\mathrm{in}}u,\,r_{\mathrm{in}})\subseteq\mathcal{C}\subseteq\mathcal{B}(\hat{x}+r_{\mathrm{out}}u,\,r_{\mathrm{out}})$.
        \item If $x\in\mathcal{C}$ and $\rho(x)<r_{\mathrm{in}}$, then $\Pi(x)$ is a singleton.
    \end{enumerate}
\end{lem}

\begin{lem}[Curvature Feasibility Bounds]\label{lem:curvature_feasibility_bounds}
    Suppose the assumptions of \cref{lem:rolling_balls} hold.
    Let $\hat{x}\in\rb{\mathcal{C}}$ with inward unit normal $u$, let $0\le\rho\le r_{\mathrm{in}}$, set $\bar{x}:=\hat{x}+\rho u$, and let $d\in S^{n-1}$ with $c:=\inner{d}{u}$. Then:
    \begin{enumerate}
        \item $\displaystyle\tau(\bar{x};d)\;\ge\;(r_{\mathrm{in}}-\rho)c+\sqrt{(r_{\mathrm{in}}-\rho)^2c^2+2r_{\mathrm{in}}\rho-\rho^2}$.
        \begin{itemize}
            \item[(1a)] If $c\ge 0$, then $\tau(\bar{x};d)\ge\min\{2c,1\}\,r_{\mathrm{in}}$.
            \item[(1b)] If $c\ge 0$, then $\tau(\bar{x};d)\ge\sqrt{r_{\mathrm{in}}\rho}$.
            \item[(1c)] $\tau(\bar{x};u)\ge 2r_{\mathrm{in}}-\rho\ge r_{\mathrm{in}}$.
        \end{itemize}
        \item $\displaystyle\tau(\bar{x};d)\;\le\;2r_{\mathrm{out}}[c]_+ + \sqrt{2r_{\mathrm{out}}\rho}$ for every $c\in[-1,1]$.
        \item $\tau(\bar{x};d)\le\rho/|c|$ for $c<0$, and $\tau(\bar{x};-u)=\rho$.
    \end{enumerate}
\end{lem}

\begin{corollary}\label{cor:boundary_eta_tau}
    Suppose the assumptions of \cref{lem:rolling_balls} hold, and let $\eta\in(0,1)$ and $\hat{x}\in\rb{\mathcal{C}}$.
    If a unit direction $d$ is an $\eta$-FD at $\hat{x}$, then $\tau(\hat{x};d)\ge 2\eta/\kappa_{\max}$;
    if $d$ is not an $\eta$-FD, then $\tau(\hat{x};d)\le 2\eta/\kappa_{\min}$.
\end{corollary}

Under the assumptions of \cref{lem:rolling_balls} every boundary point of $\mathcal{C}$ is smooth, so the analysis of nonsmooth boundary points in
\cref{sec:noised_alg_properties} is not needed. \Cref{cor:boundary_eta_tau} relates $\Delta^{[R]}$ to the
$\eta$-FD cone. For every $R>2\eta/\kappa_{\min}$ and $\bar{x}\in\rb{\mathcal{C}}$,
\begin{equation*}
    \Delta^{[R]}(\bar{x})\le\Delta_{\eta}(\bar{x}):=\sup_{d\in D^{\eta}_{\mathcal{C}}(\bar{x}),\,\norm{d}=1}[-h'(\bar{x};d)]_+,
\end{equation*}
since by \cref{cor:boundary_eta_tau} any $d$ with $\tau(\bar{x};d)>2\eta/\kappa_{\min}$ is an $\eta$-FD (\cref{cor:eps-complexity} takes
$R=8\eta/\kappa_{\min}$). Thus $\Delta^{[R]}$ accounts both for the $O(\eta)$ residual of \cref{ex:PDM_boundary} and for the step
truncation in \cref{prop:PDM_lower_bound}.

When $\mathcal{C}$ satisfies the assumptions of \cref{lem:rolling_balls}, \cref{thm:finite-time-rate} applies to \cref{alg:GRFD}
with $\alpha=0$ and $\rho_0$ of order $r_{\mathrm{in}}$, because the non-exploration output of \cref{alg:construct_directions} implements
such an oracle. We provide the constants at $\rho_0=r_{\mathrm{in}}/2$. In this regime $\rho\le\rho_0<r_{\mathrm{in}}$, so $\Pi(x)$ is
a singleton by \cref{lem:rolling_balls}.

\begin{proposition}\label{prop:quantitative-spanning}
Let $\mathcal{C}$ satisfy the assumptions of \cref{lem:rolling_balls}, let $\eta\in(0,\tfrac12]$, and let $r\in(0,\infty]$.
Then \cref{alg:construct_directions} with $\rho_0=r_{\mathrm{in}}/2$ and $\alpha=0$ has spanning-coefficient bounds $(A,B,R)$ with
\begin{equation*}
    A=\sqrt{\dim\mathcal{C}},
    \qquad
    R=\frac{8\eta}{\kappa_{\min}},
    \qquad
    B\le\frac{\sqrt{\dim\mathcal{C}}}{r}+\frac{2\,\dim\mathcal{C}\cdot\kappa_{\max}}{\eta},
\end{equation*}
and $|V(x)|\le 4\dim\mathcal{C}$.
\end{proposition}

The proof is in \cref{sec:appendix_proofs_analysis}.

\begin{remark}[Role of the adapted basis and the tilt]\label{rem:both-modifications-necessary}
On the Euclidean unit disk, where $r_{\mathrm{in}}=r_{\mathrm{out}}=1$, dropping either the adapted basis or the tilt from \cref{alg:construct_directions} yields no finite constant $B$ in \cref{def:quantitative-eta-oracle}.
Let $\hat{x}=(-1,0)$ with inward unit normal $u=(1,0)$, let $x=\hat{x}+\rho u$ with $\rho$ small, and let $d=(\tfrac12,\tfrac{\sqrt3}{2})$, so $\inner{d}{u}=\tfrac12$ and \cref{lem:curvature_feasibility_bounds}(1a) gives $\tau(x;d)\ge 1$.
This $d$ is $R$-deep for every $R\le 1$ and every $\rho$.
\begin{itemize}
    \item \textit{Unadapted basis.} Let $w_1=\tfrac1{\sqrt2}(1,1)$ and $w_2=\tfrac1{\sqrt2}(-1,1)$, and take $V=\{\pm w_1,\pm w_2\}$.
    This set positively spans $D_{\mathcal{C}}(x)$.
    Since $\{w_1,w_2\}$ is orthonormal, $d=\inner{d}{w_1}w_1+\inner{d}{w_2}w_2$ with
    $\inner{d}{w_1}=\bigl(\tfrac12+\tfrac{\sqrt3}{2}\bigr)/\sqrt2>0$ and
    $\inner{d}{w_2}=\bigl(-\tfrac12+\tfrac{\sqrt3}{2}\bigr)/\sqrt2>0$.
    Since $\inner{w_2}{u}=-\tfrac1{\sqrt2}$, \cref{lem:curvature_feasibility_bounds}(3) gives $\tau(x;w_2)\le\sqrt{2}\,\rho$, so $\lambda_{w_2}/\tau_r(x;w_2)\ge\inner{d}{w_2}/(\sqrt{2}\,\rho)=(\sqrt3-1)/(4\rho)\to\infty$.
    Adding canceling $\pm$ pairs only inflates $\sum\lambda_v/\tau_r(x;v)$, so $B=\infty$.
    \item \textit{Adapted basis, no tilt.} $V=\{\pm u,\pm w\}$ with $w=(0,1)$ forces $\lambda_w=\tfrac{\sqrt3}{2}$, and $\tau(x;w)\le\sqrt{2\rho}\to 0$ by \cref{lem:curvature_feasibility_bounds}(2). Again $B=\infty$.
    \item \textit{Adapted basis with tilt.} $\inner{d}{u}=\tfrac12\ge\eta$, so the inward-tilt case of \cref{prop:quantitative-spanning} applies: the decomposition uses $u$ and $\bar{\eta}u+(1-\bar{\eta})w$, whose feasibility radius is at least $2\bar{\eta}$ independently of $\rho$, so the resulting bound is finite.
\end{itemize}
In this example, both modifications are therefore needed for a bound on $B$ that is uniform in $\rho$.
\end{remark}

\begin{lem}\label{lem:delta-R-upgrade}
Suppose the assumptions of \cref{lem:rolling_balls} hold. Let $f$ be continuously differentiable, let $g$ be convex and finite near $\mathcal{C}$, and let $\Lambda$ be as in \cref{lem:h-prime-lipschitz-span}(1).
Let $R\in(0,r_{\mathrm{in}}]$ and $x\in\mathcal{C}$ with $\rho:=\rho(x)$.
\begin{enumerate}
    \item If $\rho\ge r_{\mathrm{in}}$ then $\Delta^{[R]}(x)=\Delta(x)$.
    \item If $\rho=0$ then $\Delta(x)\le\Delta^{[R]}(x)+2\Lambda\kappa_{\max}R$.
    \item If $0<\rho<r_{\mathrm{in}}$, with $u$ the inward normal at the nearest boundary point, the bound of (2) holds for every unit $d\in D_{\mathcal{C}}(x)$ with $\inner{d}{u}\ge 0$.
    A unit $d$ with $\inner{d}{u}<0$ enters $\Delta^{[R]}(x)$ if $\tau(x;d)\ge R$, and is excluded from it if $\tau(x;d)<R$; in the extreme case $d=-u$, $\tau(x;-u)=\rho$.
\end{enumerate}
\end{lem}

The proof is in \cref{sec:appendix_proofs_analysis}.
A bound on $\Delta^{[R]}$ therefore controls $\Delta$ when $\rho=0$ or $\rho\ge r_{\mathrm{in}}$, up to an additive $2\Lambda\kappa_{\max}R$ on the boundary.
At a point with $0<\rho<r_{\mathrm{in}}$ it controls every unit direction with $\tau(x;d)\ge R$, but not an outward direction that leaves $\mathcal{C}$ before distance $R$.
Whether a given outward direction is excluded depends on both $R$ and $\rho$, not solely on $\rho$: a nearly tangential outward $d$ can have $\tau(x;d)$ much larger than $\rho$ and possibly at least $R$.

\begin{corollary}[$\varepsilon$-complexity]\label{cor:eps-complexity}
Suppose the assumptions of \cref{lem:rolling_balls} hold, let $f\in C^{1,1}_L$, let $g$ be convex, and assume $H<\infty$.
Let $\varepsilon>0$, $R\in(0,r_{\mathrm{in}}]$, and run \cref{alg:Oracle_GRFD} with the non-exploration output of \cref{alg:construct_directions} at $\rho_0=r_{\mathrm{in}}/2$ and $\eta=R\kappa_{\min}/8$ (then $\eta\le\kappa_{\min}/(8\kappa_{\max})\le\tfrac18$, so \cref{prop:quantitative-spanning} applies).
Then some $x_k$ with $k<K$ is $(\varepsilon,R)$-$d$-stationary as soon as
\begin{equation}\label{eq:eps-complexity}
    K\;\ge\;\frac{8\dim\mathcal{C}\,(L+\gamma)H}{\varepsilon^2}+\frac{4H}{\varepsilon}\left(\frac{\sqrt{\dim\mathcal{C}}}{r}+\frac{16\,\dim\mathcal{C}\cdot\kappa_{\max}}{R\,\kappa_{\min}}\right).
\end{equation}
If in addition the constant $\Lambda$ of \cref{lem:h-prime-lipschitz-span}(1) satisfies $\Lambda>0$ and $0<\varepsilon\le 4\Lambda$, then $R=\varepsilon/(4\Lambda\kappa_{\max})\in(0,r_{\mathrm{in}}]$ is admissible and gives
$K=O\bigl(\dim\mathcal{C}\cdot H(L+\gamma+\Lambda\kappa_{\max}^2/\kappa_{\min})\varepsilon^{-2}\bigr)$ as $\varepsilon\downarrow 0$ with $r$ held fixed.
With this $R$, if $K$ satisfies \eqref{eq:eps-complexity} with $\varepsilon$ replaced by $\varepsilon/2$, then some $x_k$ with $k<K$ is $(\varepsilon/2,R)$-$d$-stationary, and at such an iterate $\Delta(x_k)\le\varepsilon$ whenever $x_k\in\rb{\mathcal{C}}$ or $\rho(x_k)\ge r_{\mathrm{in}}$.
\end{corollary}

\begin{proof}
Substitute $A=\sqrt{\dim\mathcal{C}}$ and the bound on $B$ from \cref{prop:quantitative-spanning} into \eqref{eq:finite-time-rate} and require each of the two terms to be at most $\varepsilon/2$.
This is \eqref{eq:eps-complexity}.
For the bound on $\Delta$, apply \eqref{eq:eps-complexity} at accuracy $\varepsilon/2$ with $R=\varepsilon/(4\Lambda\kappa_{\max})$; this yields an iterate with $\Delta^{[R]}(x_k)\le\varepsilon/2$. If $\rho(x_k)\ge r_{\mathrm{in}}$, \cref{lem:delta-R-upgrade}(1) gives $\Delta(x_k)=\Delta^{[R]}(x_k)\le\varepsilon/2$. If $x_k\in\rb{\mathcal{C}}$, \cref{lem:delta-R-upgrade}(2) gives $\Delta(x_k)\le\varepsilon/2+2\Lambda\kappa_{\max}R=\varepsilon$.
\end{proof}

Under the parameter choices of \cref{cor:eps-complexity},
the guarantee requires $O(\dim\mathcal{C}\,\varepsilon^{-2})$
outer iterations and
$O((\dim\mathcal{C})^2\varepsilon^{-2})$ exact line searches,
since each iteration uses $O(\dim\mathcal{C})$ directions.
These counts are not directly comparable with the
function-evaluation bounds of \citet{gratton2019direct}:
each line search here solves an exact one-dimensional
minimization problem (\cref{rem:line-search-access}).

The choices of $\rho_0$ and $\eta$ use the curvature bounds;
choosing $R$ from $\varepsilon$ also uses $\Lambda$.
Reducing $\eta$ decreases $R=8\eta/\kappa_{\min}$, extending
the guarantee to directions with shorter feasible steps,
but increases the spanning-coefficient bound $B$.
The choice $\eta=\Theta(\varepsilon)$ yields the stated complexity.

\begin{remark}[On the parameters $\gamma$ and $r$]
\label{rem:gamma-r-tension}
The finite-time bound favors smaller $\gamma$ and larger $r$.
However, $\gamma>0$ is required by the PDM rate and the
almost-sure convergence result
(\cref{thm:convergence_rate,thm:main}).
\end{remark}

\section{Geometry for Directional Spanning and Trajectories}\label{sec:noised_alg_properties}

This section is devoted to establishing the geometric properties underlying
\cref{sec:analysis}. We first address feasibility of the
non-exploration directions and almost-sure avoidance of nonsmooth
boundary points. We then establish directional spanning guarantees
and explicit noise thresholds, leading to the PDM certificate
in \cref{cor:noised-pdm-certificate}.




\medskip

\begin{lem}\label{lem:GRFD_all_considered_directions_are_feasible}
Let $\{\tilde{V}^k \}_{k\geq 0}$ be the sequence of finite directions set generated by \cref{alg:GRFD}, and let $\left\{\tilde{v}_{\mathrm{ex}}^k\right\}_{k \geq 0}$ denote the generated exploration directions. 
Then for any $k\geq 0$, it holds that $\Tilde{V}^k \setminus \left\{\tilde{v}_{\mathrm{ex}}^k\right\} \subseteq D_\mathcal{C}(x_k).$
\end{lem}
\begin{proof} 
    If $x_k$ is a nonsmooth boundary point, then by \cref{alg:construct_directions}, $\tilde{V}^k \setminus \left\{\tilde{v}_{\mathrm{ex}}^k\right\} = \emptyset$, and the result holds trivially.

    Otherwise, the non-exploration output of \cref{alg:construct_directions} consists of the normalized images of
    $\alpha\bar{\eta}\zeta_k+(1-\alpha\bar{\eta})v$ for $v\in U$ and of
    $(1-\alpha)\bar{\eta}\,u+\alpha\bar{\eta}\zeta_k+(1-\bar{\eta})v$ for $v\in T$.
    Normalization preserves membership in $D_{\mathcal{C}}(x_k)$.

    By construction $U\subseteq \ri{D_{\mathcal{C}}(x_k)}$.
    In the first branch $T=\emptyset$ and $U$ is a signed orthonormal basis of $\aff{\mathcal{C}-x_k}=D_{\mathcal{C}}(x_k)$.
    In the third branch there are two cases.
    If $\rho(x_k)=0$, then $\hat{x}=x_k$ and $\mathcal{C}$ is smooth there, so $\cl{D_{\mathcal{C}}(x_k)}$ is a halfspace in $\aff{\mathcal{C}-x_k}$ with inward unit normal $u$; the set $T=\{\pm w_i\}$ lies in the corresponding tangent hyperplane, hence in $\cl{D_{\mathcal{C}}(x_k)}$, and $u\in\ri{D_{\mathcal{C}}(x_k)}$.
    If $0<\rho(x_k)\le\rho_0$, then $x_k\in\ri{\mathcal{C}}$ and $D_{\mathcal{C}}(x_k)=\aff{\mathcal{C}-x_k}$ by \cref{lem:feasible_directions_at_relative_interior}, so $U\cup T\subseteq\ri{D_{\mathcal{C}}(x_k)}$.
    Thus $U\cup T\subseteq\cl{D_{\mathcal{C}}(x_k)}$.
    Since $D_{\mathcal{C}}(x_k)$ is a convex cone (\cref{lem:feasible_directions_are_convex_cones}), its closure is convex, and for any $\delta\in\Delta_m$,
    \begin{equation*}
        \zeta_k = \sum_{i=1}^{m} \delta_i v_i\in \cl{D_\mathcal{C}(x_k)}.
    \end{equation*}

    For $v\in U$, the untilted vector is a convex combination of $v\in\ri{D_{\mathcal{C}}(x_k)}$ and $\zeta_k\in\cl{D_{\mathcal{C}}(x_k)}$, with weight $1-\alpha\bar{\eta}>0$ on $v$.
    \Cref{lem:foundational properties} therefore yields that the combination lies in $\ri{D_{\mathcal{C}}(x_k)}$.

    For $v\in T$, the tilted vector is a convex combination of $u\in\ri{D_{\mathcal{C}}(x_k)}$, $\zeta_k\in\cl{D_{\mathcal{C}}(x_k)}$, and $v\in\cl{D_{\mathcal{C}}(x_k)}$, with weight $(1-\alpha)\bar{\eta}>0$ on $u$.
    \Cref{lem:foundational properties} therefore yields membership in $\ri{D_{\mathcal{C}}(x_k)}$.
\end{proof}

Nonsmooth boundary points form a small subset of $\rb{\mathcal{C}}$, much as the faces of dimension at most $\dim\mathcal{C}-2$ of a polytope form a small subset of its boundary. The precise statement, obtained via an isometric flattening of $\mathcal{C}$ onto a full-dimensional convex set, is given in \cref{sec:appendix_mt_background} (\cref{prop:closed_convex_sets_are_almost_everywhere_smooth,prop:nonsmooth_boundary_borel_null}). The next theorem shows that \cref{alg:GRFD} with $\alpha>0$, started away from such points, almost surely never visits one. We give a proof sketch; the complete proof is in \cref{sec:appendix_mt_background}.
\begin{theorem}\label{lem:always_smooth_almost_surely}
    Let $\mathcal{C}_{ns}\subseteq\mathcal{C}$ be the set of nonsmooth
    boundary points of $\mathcal{C}$ (a Borel set), and let $\{x_k\}_{k=0}^{K}$ be a
    sequence generated by the first $K\ge 0$ iterations of \cref{alg:GRFD}
    with noise parameter $\alpha\in(0,1)$ and
    $\eta\in(0,1]$, initialized at
    $x_0\in\mathcal{C}$ with $\mathbb{P}\!\left(x_0\in\mathcal{C}_{ns}\right)=0$.
    Then $\mathbb{P}\!\left(x_k\in\mathcal{C}_{ns}\right)=0$ for all $k\in \left\{0,\ldots,K\right\}$.
\end{theorem}

\begin{proof}[Proof sketch]
    The argument proceeds by induction on $k$; the base case holds by
    assumption. The set $\mathcal{C}_{ns}$ is a negligible subset of $\rb{\mathcal{C}}$,
    in the sense made precise in \cref{sec:appendix_mt_background}. The
    randomized construction of $\tilde{V}^k$ in
    \cref{alg:construct_directions} ensures that the direction
    $\omega_i := \tilde{v}_i^k$ has a conditional law
    absolutely continuous with respect to surface measure $\sigma$ on
    $\mathbb{S}_L$, where $L:=\aff{\mathcal{C}-x_k}$.
    This uses only that, on every branch that returns spanning directions, the differences of the elements of $U\cup T$ span $\aff{\mathcal{C}-x_k}$ when $\dim\mathcal{C}\ge 2$ (\cref{lem:direction_density}); when $\dim\mathcal{C}=1$ there are no nonsmooth boundary points.
    Because negligible subsets of
    $\rb{\mathcal{C}}$ are mapped to $\sigma$-null sets under radial projection, $\omega_i$
    avoids the radial image of $\mathcal{C}_{ns}$ almost surely. Landing in
    $\mathcal{C}_{ns}$ would force $\omega_i$ into that null set, so each
    candidate iterate $x_k + t_i^k\tilde{v}_i^k$ misses $\mathcal{C}_{ns}$
    with probability one; a union bound over the $s_k$ directions closes the
    induction.
\end{proof}

\begin{example}[Positive noise is necessary]
\label{ex:noise-needed}
The assumption $\alpha\in(0,1)$ in \cref{lem:always_smooth_almost_surely} cannot be relaxed to $\alpha=0$. Let $\mathcal{C}=\{z\in\real^2 : z\ge 0,\ \norm{z}\le 1\}$ be the quarter disc, whose nonsmooth boundary points are the corners $0$, $(1,0)$ and $(0,1)$. We run \cref{alg:GRFD} with $\alpha=0$, $\rho_0=0$, $\gamma>0$ and $r\ge 1$, using the standard basis on the first branch of \cref{alg:construct_directions}. Write $[x_i]_j$ for the $j$th coordinate of the $i$th iterate.

\textit{Construction.} Fix $\delta\in(0,1/8]$ and $K>\gamma/(4\delta)$. Let $\psi:\real\to\real$ be of class $C^2$, nonincreasing on $[0,\infty)$, constant on $[0,\tfrac{1}{4}]$, with $\psi'\le -K$ on $[\tfrac{9}{16},1]$. Let $\bar{\chi}:\real\to[0,1]$ be of class $C^2$, nonincreasing on $[0,\infty)$, equal to $1$ on $[0,\delta^2]$ and to $0$ on $[4\delta^2,\infty)$. Both exist by $C^2$ interpolation on the gaps between the prescribed pieces: $\bar{\chi}$ joins the value $1$ on $[0,\delta^2]$ to the value $0$ on $[4\delta^2,\infty)$, and $\psi$ joins its plateau on $[0,\tfrac{1}{4}]$ to an affine piece of slope $-K$ on $[\tfrac{9}{16},\infty)$. Fix $c>2\gamma$ and $M>\max\{\psi(0)-\psi(1)+\gamma/2,\ 2\sqrt{2}\,c\delta\}$. Set $g\equiv 0$ and
\begin{equation*}
  h(z)=f(z)=\psi\bigl(\norm{z}^2\bigr)+\bar{\chi}\bigl(\norm{z}^2\bigr)\bigl(-M+c(z_1+z_2)\bigr).
\end{equation*}
Then $h$ is of class $C^2$ on $\real^2$, so its second derivatives are bounded on $\mathcal{C}$. On $\mathcal{C}\cap\mathcal{B}(0,2\delta)$ one has $z_1+z_2\le\sqrt{2}\,\norm{z}<2\sqrt{2}\,\delta$, so $-M+c(z_1+z_2)<0$ and $h(z)\ge\psi(0)-M=h(0)$, with equality only at $z=0$. On $\mathcal{C}\setminus\mathcal{B}(0,2\delta)$, $\bar{\chi}(\norm{z}^2)=0$, hence $h(z)=\psi(\norm{z}^2)\ge\psi(1)>h(0)$. Thus $0$ is the unique minimizer of $h$ on $\mathcal{C}$.

\textit{Step 1: from $x_0$ to the arc.} Fix $x_0\in\mathcal{C}$ with $[x_0]_1,[x_0]_2>2\delta$ and $\tfrac{3}{4}<\norm{x_0}<1$. Along $e_1$ the second coordinate is unchanged, so the segment stays outside $\mathcal{B}(0,2\delta)$. There $\frac{d}{dt}\norm{x_0+te_1}^2=2[x_0+te_1]_1>4\delta$ and $\norm{x_0+te_1}^2\ge\tfrac{9}{16}$. The maximal feasible step along $+e_1$ is strictly less than $1$, and $4\delta K>\gamma$, so the derivative of $t\mapsto h(x_0+te_1)+\tfrac{\gamma}{2}t^2$ is strictly less than $-4\delta K+\gamma t<0$ up to the arc. The unique minimizer of the regularized line search is therefore the maximal feasible step, which lands on the open arc. Along $-e_1$ the value of $h$ does not decrease, so the step is $0$; $\pm e_2$ are symmetric. Let $E=\{\inner{\tilde{v}_{\mathrm{ex}}}{x_0}\ge 0\}$. Then $\mathbb{P}(E)=1/2$, and on $E$ one has $\norm{x_0+t\tilde{v}_{\mathrm{ex}}}\ge\norm{x_0}>\tfrac{3}{4}>2\delta$ for every $t\ge 0$, so the exploration segment avoids $\mathcal{B}(0,2\delta)$. On $E$, every candidate has value at least $\psi(1)$, with equality only on the arc, and almost surely not at a corner. Hence, on $E$, $x_1$ lies on the open arc.

\textit{Step 2: from the arc to the corner.} On the event $E$, write $p=x_1$, so $p\ge 0$, $\norm{p}=1$ and $[p]_1+[p]_2\ge 1$. At $p$ the third branch of \cref{alg:construct_directions} with $\alpha=0$ returns the inward normal $u=-p$, and $\tau(p;u)=1\le r$. Let $s:=1-t$, $A:=[p]_1+[p]_2$, and $\varphi(t)=h(sp)+\tfrac{\gamma}{2}t^2$ on $[0,1]$, so $\varphi(1)=\psi(0)-M+\gamma/2$.
If $s\ge 2\delta$, then $\varphi(t)\ge\psi(1)>\varphi(1)$.
On $[0,2\delta]$ the map $s\mapsto h(sp)$ is $C^2$ with
\begin{equation*}
  \frac{d}{ds}h(sp)=2s\bar{\chi}'(s^2)(-M+csA)+cA\bar{\chi}(s^2)\ge 0,
\end{equation*}
since $\bar{\chi}'\le 0$ and $-M+csA<0$; thus $s\mapsto h(sp)$ is nondecreasing on $[0,2\delta]$.
If $s\in[\delta,2\delta)$, then $\varphi(t)\ge h(\delta p)+\tfrac{\gamma}{2}(1-2\delta)^2\ge\psi(0)-M+c\delta+\tfrac{\gamma}{2}(1-2\delta)^2$, which exceeds $\varphi(1)$ since $c\delta>2\gamma\delta\ge\tfrac{\gamma}{2}\bigl(1-(1-2\delta)^2\bigr)$.
If $s<\delta$, then $\varphi'(t)=-cA+\gamma t\le\gamma-c<0$.
Hence $t=1$ is the unique minimizer of $\varphi$, and the candidate along $u$ is $0$. The feasible line-search segments along the images of $T$ do not pass through $0$, and the exploration segment passes through $0$ with probability zero, so every other candidate almost surely has a strictly larger value of $h$. Therefore, conditionally on $E$, one has $x_2=0$ almost surely, and $\mathbb{P}(x_2\in\mathcal{C}_{ns})\ge\mathbb{P}(E)>0$.
\end{example}

\begin{remark}[Non-degenerate direction sets]\label{rem:nondegenerate-direction-sets}
\Cref{thm:convergence_rate} bounds the PDM measure on $\tilde{V}^k$. This is uninformative at an iterate where \cref{alg:construct_directions} takes its second branch ($\rho=0$ and $\mathcal{C}$ not smooth at $x_k$) and returns only the exploration direction.

If $\alpha\in(0,1)$ and $x_0$ is almost surely not a nonsmooth boundary point of $\mathcal{C}$, then by \cref{lem:always_smooth_almost_surely} no iterate is a nonsmooth boundary point, almost surely. Since spanning directions are constructed whenever $x_k$ is not a nonsmooth boundary point, $\tilde{V}^k \neq \{\tilde{v}_{\mathrm{ex}}^k\}$ almost surely for every $k$.

If $\mathcal{C}$ is smooth at every point of $\rb{\mathcal{C}}$, the second branch is never taken for any $\alpha\in[0,1)$, so $\tilde{V}^k \neq \{\tilde{v}_{\mathrm{ex}}^k\}$ at every iterate. When neither condition holds, \cref{ex:noise-needed} shows that the iterates can reach a nonsmooth boundary point.
\end{remark}

The following lemma shows that, at every point that is not a nonsmooth boundary point, the non-exploration output of \cref{alg:construct_directions} with $\alpha=0$ positively spans $D^\eta_{\mathcal{C}}(\bar{x})$. If $\mathcal{C}$ has no nonsmooth boundary points, that output is therefore an $\eta$-DO, and \cref{alg:GRFD} with $\alpha=0$ is an instance of \cref{alg:Oracle_GRFD}.

\begin{lem}\label{lem:tilde_v_spans_eta_feasible}
    Let $\bar{x} \in \mathcal{C}$ and $\eta\in(0,1]$, and assume that $\bar{x}$ is not a nonsmooth boundary point of $\mathcal{C}$. Denote by $\tilde{V}$ the set of directions returned by \cref{alg:construct_directions} for $\bar{x}$, with this $\eta$ and $\alpha = 0$, where $v_{\mathrm{ex}}\in \tilde{V}$ is the exploration direction. Then, 
    \begin{equation*}
        D_{\mathcal{C}}^{\eta}(\bar{x}) \subseteq \pos{\tilde{V} \setminus \left\{v_{\mathrm{ex}}\right\}}.
    \end{equation*}
\end{lem}
\begin{proof}
Let $d\in D_{\mathcal{C}}^{\eta}(\bar{x})$.
By \cref{def:eta_feasible_directions}, $d\neq 0$.
Positive homogeneity of the positive span then lets us assume without loss of generality that $\|d\|=1$.
Set $\bar{\eta}=\eta/(2\sqrt{\dim\mathcal{C}})$, as in \cref{alg:construct_directions}.
Because $\bar{x}$ is not a nonsmooth boundary point, the second branch of \cref{alg:construct_directions} is excluded, so $\tilde{V}\setminus\{\tilde{v}_{\mathrm{ex}}\}$ is produced by the $T$- and $U$-updates of \cref{alg:construct_directions}.

\textbf{Case 1: $\rho(\bar{x})>\rho_0$.} By \cref{lem:eta_FD_properties}(a), every $0 \ne d \in D_\mathcal{C}(\bar{x})$ is an $\eta$-FD, and $D_\mathcal{C}(\bar{x}) = \aff{\mathcal{C}-x_0}$.
The first branch of \cref{alg:construct_directions} sets $T=\emptyset$ and $U=\{\pm w_1,\ldots,\pm w_{\dim\mathcal{C}}\}$, an untilted signed orthonormal basis of $\aff{\mathcal{C}-\bar{x}}$.
At $\alpha=0$ the untilted output is $U$ itself, which positively spans $\aff{\mathcal{C}-x_0}=D_\mathcal{C}(\bar{x})$.

\textbf{Case 2: $\bar{x} \in \rb{\mathcal{C}}$, $\mathcal{C}$ smooth at $\bar{x}$.} By \cref{lem:eta_FD_properties}(b), $d$ is an $\eta$-FD if and only if $\langle d, u \rangle \geq \eta$. Decompose $d$ in the orthonormal basis $W = \{u, w_1, \ldots, w_{\dim\mathcal{C}-1}\}$:
\begin{equation*}
    d = c \cdot u + \sum_{i=1}^{\dim\mathcal{C}-1} a_i w_i, \qquad c = \langle d, u \rangle \geq \eta.
\end{equation*}
Since $\|d\|=1$ and $W$ is orthonormal, $c^2 + \sum a_i^2 = 1$, and by $\|\cdot\|_1 \leq \sqrt{\dim\mathcal{C}}\|\cdot\|_2$,
\begin{equation*}
    \sum_{i=1}^{\dim\mathcal{C}-1}|a_i| \leq \sqrt{\dim\mathcal{C}}.
\end{equation*}
Set $\tilde{c} = c - \frac{\bar{\eta}}{1-\bar{\eta}}\sum|a_i|$ and $\tilde{a}_i = \frac{a_i}{1-\bar{\eta}}$. Since $1 - \bar{\eta} > \frac{1}{2}$, we have $\frac{\bar{\eta}}{1-\bar{\eta}} < \frac{\eta}{\sqrt{\dim\mathcal{C}}}$, and combined with $\sum |a_i| \leq \sqrt{\dim\mathcal{C}}$ and $c \geq \eta$, this gives $\tilde{c} > c - \eta \geq 0$. We now rewrite $d$ by introducing the substitution $a_i = (1-\bar{\eta})\tilde{a}_i$ and isolating the $u$-components:
\begin{align*}
    d &= c \cdot u + \sum_{i=1}^{\dim\mathcal{C}-1} |a_i|\,\mathrm{sign}(a_i) \cdot w_i \\
      &= \tilde{c} \cdot u + \sum_{i=1}^{\dim\mathcal{C}-1} |\tilde{a}_i|
         \Bigl(\bar{\eta}\, u + (1-\bar{\eta})\,\mathrm{sign}(a_i) \cdot w_i\Bigr).
\end{align*}
At a smooth boundary point the third branch sets $U=\{u\}$ and $T=\{\pm w_1,\ldots,\pm w_{\dim\mathcal{C}-1}\}$.
At $\alpha=0$ the output (excluding exploration) consists of the untilted vector $u$ together with the normalized tilts of the directions in $T$, and therefore positively spans the same cone as
\begin{equation*}
    \bigl\{u,\;\pm(1-\bar{\eta})w_1+\bar{\eta}u,\;\ldots,\;\pm(1-\bar{\eta})w_{\dim\mathcal{C}-1}+\bar{\eta}u\bigr\}.
\end{equation*}
Each term in the decomposition of $d$ is a nonnegative multiple of an element of this set, so $d \in \pos{\tilde{V}\setminus\{v_{\mathrm{ex}}\}}$.

\textbf{Case 3: $0<\rho(\bar{x})\le\rho_0$.}
The point is interior, so $D^\eta_{\mathcal{C}}(\bar{x})=D_{\mathcal{C}}(\bar{x})=\aff{\mathcal{C}-\bar{x}}$ by \cref{lem:eta_FD_properties}(a).
The third branch of \cref{alg:construct_directions} produces $U=\{\pm u,\pm w_1,\ldots,\pm w_{\dim\mathcal{C}-1}\}$, which spans $\aff{\mathcal{C}-\bar{x}}$.
At $\alpha=0$ the untilted output includes $U$, and the claim follows.

Since $d \in D_\mathcal{C}^\eta(\bar{x})$ was arbitrary, the result follows in all cases.
\end{proof}

\Cref{lem:tilde_v_noised_spans_beta_eta_feasible} extends \cref{lem:tilde_v_spans_eta_feasible} to $\alpha>0$, the regime in which \cref{lem:always_smooth_almost_surely} applies. \Cref{ex:noise-needed} shows that nonzero noise is necessary for that almost-sure avoidance. For sufficiently small $\alpha$, the directions positively span $D_{\mathcal{C}}(\bar{x})$ at every interior point, independently of $\rho_0$, and the relaxed cone $D^{\beta\eta}_{\mathcal{C}}(\bar{x})$ at smooth boundary points.


\begin{theorem}\label{lem:tilde_v_noised_spans_beta_eta_feasible}
    Let $\bar{x} \in \mathcal{C}$ and $\eta\in(0,1]$. Denote by $\tilde{V}_{\bar{\alpha}}$ the set of directions returned by \cref{alg:construct_directions} for $\bar{x}$, with this $\eta$ and $\alpha = \bar{\alpha}$, where $v_{\mathrm{ex}}^{\bar{\alpha}}\in \tilde{V}_{\bar{\alpha}}$ is the exploration direction. 
    \begin{enumerate}
        \item If $\bar{x} \in \ri{\mathcal{C}}$, then there exists $\tilde{\alpha} > 0$ 
        such that for every $\bar{\alpha} \in (0, \tilde{\alpha})$, 
        \begin{equation*}
            D_\mathcal{C}(\bar{x}) \subseteq \pos{\tilde{V}_{\bar{\alpha}}\setminus\left\{v_{\mathrm{ex}}^{\bar{\alpha}}\right\}}.
        \end{equation*}
        \item If $\bar{x} \in \rb{\mathcal{C}}$ and $\mathcal{C}$ is smooth at $\bar{x}$, then 
        for every $\beta \in (1, 1/\eta)$ there exists $\tilde{\alpha} > 0$ 
        such that for every $\bar{\alpha} \in (0, \tilde{\alpha})$,
        \begin{equation*}
            D_\mathcal{C}^{\beta \eta}(\bar{x}) \subseteq \pos{\tilde{V}_{\bar{\alpha}} \setminus \left\{v_{\mathrm{ex}}^{\bar{\alpha}}\right\}}.
        \end{equation*}
    \end{enumerate}
\end{theorem}
To prove \cref{lem:tilde_v_noised_spans_beta_eta_feasible}, we first prove two auxiliary lemmas, \cref{lem:interior_gap} and \cref{lem:boundary_gap}, which are used to prove cases $1$ and $2$ of \cref{lem:tilde_v_noised_spans_beta_eta_feasible}, respectively.

\begin{lem}[Interior gap]\label{lem:interior_gap}
    Let $\bar{x} \in \ri{\mathcal{C}}$. Denote by $\tilde{V}_0$ the set of directions returned by \cref{alg:construct_directions} for $\bar{x}$, with $\alpha = 0$, where $v_{\mathrm{ex}}^0\in \tilde{V}$ is the exploration direction. Then there exists a compact set 
    $S \subseteq D_\mathcal{C}(\bar{x})$ and a polytope $Q \subseteq \pos{\tilde{V}_0 \setminus \left\{v_{\mathrm{ex}}^0\right\}}$ with 
    the following properties:
    \begin{enumerate}
        \item $S \subseteq \ri{Q}$ and $\dist(S, \rb{Q}) > 0$.
        \item For every $d \in D_\mathcal{C}(\bar{x})$ there exists $\lambda > 0$ such that 
        $\lambda d \in S$.
    \end{enumerate}
\end{lem}
\begin{proof}
We construct the objects $S$ and $Q$ claimed by the lemma, and verify both properties.

Define the non-exploration directions $\bar{V}_0 := \tilde{V}_0 \setminus \left\{v_{\mathrm{ex}}^0\right\}$, and let $W = \{w_1, \ldots, w_{\dim\mathcal{C}}\}$ be the orthonormal basis of $\aff{\mathcal{C} - \bar{x}}$ generated by \cref{alg:construct_directions}.
Let $U=\{w_1,-w_1,\ldots,w_{\dim\mathcal{C}},-w_{\dim\mathcal{C}}\}$ be the untilted signed basis, so that $U\subseteq\bar{V}_0$, with equality on the first branch.
Define the polytope
\[
Q := \conv\left\{ \sum_{i=1}^{\dim\mathcal{C}} \theta_i w_i \;\middle|\; \theta_i \in \{-1, 1\} \right\},
\]
and the compact set $S := D_\mathcal{C}(\bar{x}) \cap \mathcal{B}[0, 0.5]$.
We claim that $Q \subseteq \pos{U}\subseteq\pos{\bar{V}_0}$ and $S \subseteq \ri{Q}$, which together establish the setup for both properties.

For the first claim, each vertex $\sum_{i=1}^{\dim\mathcal{C}} \theta_i w_i$ of $Q$ is a nonnegative combination of elements of $U$, so $Q \subseteq \pos{U}\subseteq\pos{\bar{V}_0}$. For the second, define the linear isometry
\[
\Phi(z) := \sum_{i=1}^{\dim\mathcal{C}} z_i w_i,
\]
from $\mathbb{R}^{\dim\mathcal{C}}$ onto $\aff{\mathcal{C}-\bar{x}}$. Since $\mathcal{B}_\infty[0,1]$ is the convex hull of the $2^{\dim\mathcal{C}}$ sign vectors in $\{-1,1\}^{\dim\mathcal{C}}$, 
and $\Phi$ is linear (hence preserves convex combinations), we have
\[
\Phi(\mathcal{B}_\infty[0,1]) = \conv\left\{\Phi(\theta) \;\middle|\; \theta \in \{-1,1\}^{\dim\mathcal{C}}\right\} 
= \conv\left\{\sum_{i=1}^{\dim\mathcal{C}} \theta_i w_i \;\middle|\; \theta_i \in \{-1,1\}\right\} = Q.
\]
As $\Phi$ is an isometry, $\Phi(\mathcal{B}_2[0,1]) = \mathcal{B}[0,1] \cap \aff{\mathcal{C}-\bar{x}}$. The norm inequality $\mathcal{B}_2[0,1] \subseteq \mathcal{B}_\infty[0,1]$ then gives
\begin{equation}\label{eq:l2_ball_in_Q}
\mathcal{B}[0,1] \cap \aff{\mathcal{C} - \bar{x}} = \Phi(\mathcal{B}_2[0,1]) \subseteq \Phi(\mathcal{B}_\infty[0,1]) = Q.
\end{equation}
Since $S = \mathcal{B}[0, 0.5] \cap \aff{\mathcal{C} - \bar{x}} \subseteq \frac{1}{2}Q$, and $Q$ is a full-dimensional polytope in $\aff{\mathcal{C} - \bar{x}}$ containing the origin in its interior, we have $\frac{1}{2} Q \subseteq \ri{Q}$, and hence $S \subseteq \ri{Q}$.

Since $S$ is compact, $\rb{Q}$ is closed, and $S \subseteq \ri{Q}$ implies 
$S \cap \rb{Q} = \emptyset$, it follows that $\dist(S, \rb{Q}) > 0$ 
\cite[Section~9.5, Problem~72]{RoydenFitzpatrickRealAnalysis4e}.

For any $d \in D_\mathcal{C}(\bar{x})$ with $d \ne 0$, set $\lambda = 0.5 / \|d\|_2 > 0$. Then $\lambda d \in D_\mathcal{C}(\bar{x})$ (since $D_\mathcal{C}(\bar{x})$ is a cone) and $\|\lambda d\|_2 = 0.5$, so $\lambda d \in S$. The case $d = 0$ is immediate.
\end{proof}

\begin{remark}[Interior gap at every interior point]\label{rem:interior-gap-shallow}
Let $\bar{x}\in\ri{\mathcal{C}}$.
By \cref{lem:eta_FD_properties}(a), $D_{\mathcal{C}}(\bar{x})=\aff{\mathcal{C}-\bar{x}}$, so
$\ri{D_{\mathcal{C}}(\bar{x})}=D_{\mathcal{C}}(\bar{x})$.
On either interior branch, $U$ is a signed orthonormal basis of $\aff{\mathcal{C}-\bar{x}}$:
on the first branch $T=\emptyset$ and $U=\{w_1,-w_1,\ldots,w_{\dim\mathcal{C}},-w_{\dim\mathcal{C}}\}$;
on the third branch one of the $\dim\mathcal{C}$ basis vectors is the inward unit normal $u$, and
\begin{equation*}
    T=\{w_1,-w_1,\ldots,w_{\dim\mathcal{C}-1},-w_{\dim\mathcal{C}-1}\}\subseteq U=\{u,-u,w_1,-w_1,\ldots,w_{\dim\mathcal{C}-1},-w_{\dim\mathcal{C}-1}\}.
\end{equation*}
The sets $S$ and $Q$ in \cref{lem:interior_gap} depend only on this signed basis, so they
serve on both interior branches, and $Q\subseteq\pos{U}$.
On the third branch at $\alpha=0$ the construction also returns the tilted images of $T$, hence $\pos{U}$ is
contained in the positive span of the $\alpha=0$ non-exploration output; on the first branch the two coincide.
\end{remark}

\begin{lem}[Boundary gap]\label{lem:boundary_gap}
    Let $\bar{x} \in \rb{\mathcal{C}}$ with $\mathcal{C}$ smooth at $\bar{x}$, let $\eta \in (0,1)$, 
    and let $\beta \in (1, 1/\eta)$. Denote by $\tilde{V}_0$ the set of directions returned by \cref{alg:construct_directions} for $\bar{x}$, with $\alpha = 0$, where $v_{\mathrm{ex}}^0\in \tilde{V}_0$ is the exploration direction. Then there exist a compact set 
    $S \subseteq D_\mathcal{C}^{\beta\eta}(\bar{x})$ and a polytope $P \subseteq \pos{\tilde{V}_0 \setminus \left\{v_{\mathrm{ex}}^0\right\}}$
    with the following properties:
    \begin{enumerate}
        \item $S \subseteq \ri{P}$ and $\dist(S, \rb{P}) > 0$.
        \item For every $d \in D_\mathcal{C}^{\beta\eta}(\bar{x})$ there exists $\lambda > 0$ 
        such that $\lambda d \in S$.
    \end{enumerate}
\end{lem}
\begin{proof}
Denote the non-exploration directions generated by \cref{alg:construct_directions} by
\begin{equation*}
    \bar{V}_0 := \tilde{V}_0 \setminus \left\{v_{\mathrm{ex}}^0\right\}.
\end{equation*}

By \cref{lem:eta_FD_properties}(b), 
\begin{equation}\label{eq:eta_feasible_directions}
    D_\mathcal{C}^\eta (\bar{x}) = \left\{0 \ne d \in \real^n \; | \; \left\langle \dfrac{d}{\left\Vert d \right\Vert}, u \right\rangle \ge \eta\right\} \cap \aff{\mathcal{C}-\bar{x}},
\end{equation}
for a nonzero $u \in \aff{\mathcal{C}-\bar{x}}$ such that 
\begin{equation}\label{eq:closure_of_feasible_directions}
    \cl{D_\mathcal{C}(\bar{x})} = \left\{d \in \real^n \; | \; \langle d, u \rangle \ge 0\right\} \cap \aff{\mathcal{C}-\bar{x}}.
\end{equation}
By the same arguments,
\begin{equation}\label{eq:beta_eta_feasible_directions}
    D_\mathcal{C}^{\beta\eta} (\bar{x}) = \left\{0 \ne d \in \real^n \; | \; \left\langle \dfrac{d}{\left\Vert d \right\Vert}, u \right\rangle \ge \beta\eta\right\} \cap \aff{\mathcal{C}-\bar{x}}.
\end{equation}
Since $\beta > 1$, the containment $D_\mathcal{C}^{\beta\eta}(\bar{x}) \subseteq D_\mathcal{C}^\eta(\bar{x})$ is strict, and therefore
\begin{equation*}
D_\mathcal{C}^{\beta\eta}(\bar{x}) \subseteq \ri{D_\mathcal{C}^\eta(\bar{x})}.
\end{equation*}
By \cref{lem:tilde_v_spans_eta_feasible}, $D_\mathcal{C}^\eta(\bar{x}) \subseteq \pos{\bar{V}_0}$, and since both sets have the same affine hull,
\begin{equation}\label{eq:beta_eta_in_ri_of_V0}
D_\mathcal{C}^{\beta \eta}(\bar{x}) \subseteq \ri{\pos{\bar{V}_0}}.    
\end{equation}

Let $u \in \mathbb{R}^n$ be the unit vector from \eqref{eq:closure_of_feasible_directions}, and consider the closed halfspace
\[
H_{u,-1} := \left\{x \in \mathbb{R}^n \,\middle|\, \langle x, u \rangle \le 1 \right\}.
\]
Set $\bar{\eta}=\eta/(2\sqrt{\dim\mathcal{C}})$, as in \cref{alg:construct_directions}.
At a smooth boundary point the third branch sets $U=\{u\}$ and $T=\{\pm w_1,\ldots,\pm w_{\dim\mathcal{C}-1}\}$, where $\{u,w_1,\ldots,w_{\dim\mathcal{C}-1}\}$ is an orthonormal basis of $\aff{\mathcal{C}-\bar{x}}$.
At $\alpha=0$, each vector in $\bar{V}_0$ is obtained by normalizing $u$ or $\bar{\eta}\,u+(1-\bar{\eta})v$ for some $v\in T$.
In particular $\langle v, u \rangle > 0$ for every $v \in \bar{V}_0$. We show that the set $\pos{\bar{V}_0} \cap H_{u,-1}$ is a polytope with vertex description
\begin{equation}\label{eq:truncated_spanning_set_is_polytope}
    P := \pos{\bar{V}_0} \cap H_{u,-1} = \conv{\left(\left\{0\right\} \cup \left\{ \frac{1}{\langle v, u \rangle} v \,\middle|\, v \in \bar{V}_0 \right\}\right)}.
\end{equation}
To see this, we consider $\ext{(P)}$, the extreme points of $P$. Let $w \in P$. Since $w \in \pos{\bar{V}_0}$, there exist $\lambda_1,\ldots,\lambda_s \in \real_{+}$, $v_1,\ldots,v_s \in \bar{V}_0$ such that $w = \sum_{i=1}^s \lambda_i v_i$. If $v_i=u$, then $\langle v_i, u \rangle = 1$. If $v_i=\mathrm{normalize}\bigl(\bar{\eta}\,u+(1-\bar{\eta})v\bigr)$ for some $v\in T$, then
\begin{equation*}
    \langle v_i, u \rangle = \frac{\bar{\eta}}{\sqrt{\bar{\eta}^2+(1-\bar{\eta})^2}} \in (\bar{\eta},1).
\end{equation*}
In both cases $v_i \in H_{u,-1}$, so $v_i \in P$ for all $i$. Therefore, if $w \notin \real_+ v_i$ for some $i$, then $w$ is not an extreme point of $P$. If $w \in \real_+ v_i$ for some $i$ but $w \notin \left\{0, \frac{1}{\langle v_i, u \rangle} v_i\right\}$, then $w$ lies in the convex hull of these two points and is again not extreme. Thus,
\begin{equation*}
    \ext{(P)} \subseteq \left\{0\right\} \cup \left\{ \frac{1}{\langle v, u \rangle} v \,\middle|\, v \in \bar{V}_0 \right\}.
\end{equation*}
To apply the Krein--Milman theorem, we verify that $P$ is nonempty, compact, and convex. 
Nonemptiness is trivial since $0 \in P$. Convexity holds as $P$ is an intersection of two 
convex sets. Closedness follows as $P$ is an intersection of two closed sets, where 
closedness of $\pos{\bar{V}_0}$ follows from \cite[Theorem 19.1]{rockafellar1997convex}. 
For boundedness, suppose for contradiction that there exists a sequence $\{w_i\} \subseteq P$ 
with $\|w_i\| \to \infty$. Writing $w_i = \sum_j \lambda_j^i v_j$ with $\lambda_j^i \geq 0$, 
we get $\sum_j \lambda_j^i \to \infty$, and hence
\begin{equation*}
    \langle w_i, u \rangle = \left\langle \sum_{j} \lambda_j^i v_j, u \right\rangle \ge 
    \left(\min_{j \in [s]} \langle v_j, u \rangle\right) \sum_{j} \lambda_j^i \rightarrow \infty,
\end{equation*}
contradicting $w_i \in H_{u,-1}$. Thus $P$ is nonempty, compact, and convex, so by \cite[The Krein-Milman Lemma]{RoydenFitzpatrickRealAnalysis4e}, $\ext{(P)}$ is 
nonempty, and by \cite[The Krein-Milman Theorem]{RoydenFitzpatrickRealAnalysis4e}, 
$P = \conv(\ext{(P)})$. It is possible that $\left\{0\right\} \cup \left\{ 
\frac{1}{\langle v, u \rangle} v \,\middle|\, v \in \bar{V}_0 \right\} \setminus \ext{(P)}$ 
is nonempty, but since every point in this set belongs to $P$, adding them does not change 
the convex hull. Therefore, \eqref{eq:truncated_spanning_set_is_polytope} holds.

Define
\[
S := D_\mathcal{C}^{\beta \eta}(\bar{x}) \cap \left\{ d \in \mathbb{R}^n \,\middle|\, 0.25 \le \|d\| \le 0.5 \right\}.
\]
We first show $S \subseteq \ri{P}$. Since for any $d$ with $\|d\| \le 0.5$ we have $\langle d, u \rangle \le 0.5 < 1$, it follows that $\{d : \|d\| \le 0.5\} \subseteq \ri{H_{u,-1}}$, and therefore by \eqref{eq:beta_eta_in_ri_of_V0},
\[
S \subseteq \ri{\pos{\bar{V}_0}} \cap \ri{H_{u,-1}} = \ri{P}.
\]
It remains to show $\dist(S, \rb{P}) > 0$. The set $S$ is compact: it is bounded by definition, and closed since $f(d) = \langle d/\|d\|, u \rangle$ is continuous away from $0$, making $D_\mathcal{C}^{\beta\eta}(\bar{x})$ closed on $\{\|d\| \ge 0.25\}$. Since $S$ is compact, $\rb{P}$ is closed, and $S \cap \rb{P} = \emptyset$, it follows that $\dist(S, \rb{P}) > 0$ \cite[Section~9.5, Problem~72]{RoydenFitzpatrickRealAnalysis4e}.

Let $d \in D_\mathcal{C}^{\beta\eta}(\bar{x})$. By \cref{def:eta_feasible_directions}, $d \ne 0$. Set $\lambda := 0.5 / \|d\| > 0$. Since $D_\mathcal{C}^{\beta\eta}(\bar{x})$ is a cone (for any $\mu > 0$, the direction $\mu d$ satisfies $\langle \mu d / \|\mu d\|, u \rangle = \langle d / \|d\|, u \rangle \ge \beta\eta$ and $\mu d \in \aff{\mathcal{C} - \bar{x}}$, so $\mu d \in D_\mathcal{C}^{\beta\eta}(\bar{x})$), we have $\lambda d \in D_\mathcal{C}^{\beta\eta}(\bar{x})$. Moreover, $\|\lambda d\| = 0.5 \in [0.25, 0.5]$, so $\lambda d \in S$.
\end{proof}

The proof of \cref{lem:tilde_v_noised_spans_beta_eta_feasible} combines the gap
lemmas above with two elementary facts: a quantitative Hausdorff bound for convex
hulls indexed by a common set, and a criterion under which Hausdorff proximity
forces containment. We record them here, as both are applied with constants that
do not depend on the shared random direction $\zeta$.

\begin{lem}[Hausdorff bound for matched convex hulls]\label{lem:hausdorff_matched_convex_hulls}
    Let $N \in \mathbb{N}$ and let $\{a_j\}_{j \in [N]}, \{b_j\}_{j \in [N]} \subseteq \real^n$. Then
    \begin{equation*}
        d_H\!\left(\conv\{a_j\}_{j \in [N]},\ \conv\{b_j\}_{j \in [N]}\right)
        \le \max_{j \in [N]} \norm{a_j - b_j}.
    \end{equation*}
\end{lem}
\begin{proof}
    Set $\varepsilon := \max_{j \in [N]} \norm{a_j - b_j}$. Let $x \in \conv\{a_j\}_{j \in [N]}$,
    so $x = \sum_{j \in [N]} \lambda_j a_j$ for some $\lambda \in \Delta_N$. The point
    $y := \sum_{j \in [N]} \lambda_j b_j$ lies in $\conv\{b_j\}_{j \in [N]}$ and satisfies
    \begin{equation*}
        \norm{x - y} = \norm{\sum_{j \in [N]} \lambda_j (a_j - b_j)}
        \le \sum_{j \in [N]} \lambda_j \norm{a_j - b_j} \le \varepsilon.
    \end{equation*}
    Hence $\sup_{x \in \conv\{a_j\}} \inf_{y \in \conv\{b_j\}} \norm{x - y} \le \varepsilon$.
    Exchanging the roles of $\{a_j\}$ and $\{b_j\}$ yields the reverse bound, so
    $d_H \le \varepsilon$.
\end{proof}

\begin{lem}[Containment under Hausdorff proximity]\label{lem:containment_under_hausdorff}
    Let $K \subseteq \real^n$ be a nonempty compact convex set, and let $S \subseteq \ri{K}$
    satisfy $\rho := \dist(S, \rb{K}) > 0$. Then for every compact convex set
    $K' \subseteq \aff{K}$ with $d_H(K, K') < \rho$, it holds that $S \subseteq K'$.
\end{lem}
\begin{proof}
    Fix $s \in S$ and suppose, for contradiction, that $s \notin K'$. Since $K'$ is a
    nonempty closed convex subset of $\aff{K}$ and $s \in \aff{K}$, separation within
    $\aff{K}$ yields a unit vector $n$ in the direction subspace of $\aff{K}$ with
    \begin{equation*}
        \inner{n}{s} > h_{K'}(n) := \sup_{y \in K'} \inner{n}{y}.
    \end{equation*}
    Because $\dist(s, \rb{K}) \ge \rho$, the closed ball of radius $\rho$ around $s$ within
    $\aff{K}$ is contained in $K$; in particular $s + \rho n \in K$, whence
    $\inner{n}{s} + \rho \le h_K(n)$. For unit $n$ one has
    $|h_K(n) - h_{K'}(n)| \le d_H(K, K')$, so
    \begin{equation*}
        \inner{n}{s} + \rho \le h_K(n) \le h_{K'}(n) + d_H(K, K')
        < h_{K'}(n) + \rho,
    \end{equation*}
    giving $\inner{n}{s} < h_{K'}(n)$, a contradiction. Hence $s \in K'$, and since
    $s \in S$ was arbitrary, $S \subseteq K'$.
\end{proof}

\begin{proof}[Proof of \cref{lem:tilde_v_noised_spans_beta_eta_feasible}]

Denote the non-exploration directions generated by \cref{alg:construct_directions} by
\begin{equation*}
    \bar{V}_{\bar{\alpha}} := \tilde{V}_{\bar{\alpha}} \setminus \left\{v_{\mathrm{ex}}^{\bar{\alpha}}\right\}.
\end{equation*}

\textbf{Case 1: $\bar{x} \in \ri{\mathcal{C}}$.}
Let $U=\{w_1,-w_1,\ldots,w_{\dim\mathcal{C}},-w_{\dim\mathcal{C}}\}$ be the untilted signed orthonormal basis of
$\aff{\mathcal{C}-\bar{x}}$ produced by \cref{alg:construct_directions} at $\bar{x}$
(in the third branch, one of the $\dim\mathcal{C}$ basis vectors is the inward unit normal $u$).
By \cref{lem:interior_gap,rem:interior-gap-shallow}, there exist a compact set
$S\subseteq D_\mathcal{C}(\bar{x})$ and a polytope $Q\subseteq\pos{U}$ such that
$S\subseteq\ri{Q}$, $\dist(S,\rb{Q})>0$, and every direction in $D_\mathcal{C}(\bar{x})$
can be scaled into $S$.

We now construct a family of polytopes that track the noise level~$\bar{\alpha}$.
By the proof of \cref{lem:interior_gap}, the polytope~$Q$ is the convex hull
of the $2^{\dim\mathcal{C}}$ vertices
\begin{equation*}
    w_\theta = \sum_{i=1}^{\dim\mathcal{C}} \theta_i w_i, \qquad \theta \in \{-1,1\}^{\dim\mathcal{C}},
\end{equation*}
corresponding to all sign combinations over the orthonormal basis.

Let $\{\tilde{v}_1^{\bar{\alpha}}, \ldots,
\tilde{v}_{2\dim\mathcal{C}}^{\bar{\alpha}}\}$ be the unnormalized untilted combinations of
\cref{alg:construct_directions}, enumerated so that
$\tilde{v}_{2i-1}^{\bar{\alpha}}$ and $\tilde{v}_{2i}^{\bar{\alpha}}$ are the perturbations of
$w_i$ and $-w_i$:
\begin{equation*}
    \tilde{v}_j^{\bar{\alpha}}
    = \bar{\alpha}\bar{\eta}\zeta
      + (1 - \bar{\alpha}\bar{\eta})\, v_j,
    \qquad j \in [2\dim\mathcal{C}],
\end{equation*}
where $v_j$ denotes the corresponding element of $U$ and $\zeta$ is the shared random direction.
For each $j$, \cref{alg:construct_directions} includes $\mathrm{normalize}(\tilde{v}_j^{\bar{\alpha}})$ in $\bar{V}_{\bar{\alpha}}$, so $\tilde{v}_j^{\bar{\alpha}}\in\pos{\bar{V}_{\bar{\alpha}}}$.
If $T$ is nonempty the algorithm also returns normalized tilts of $T$; those are not used below.
Note that $\tilde{v}_j^{\bar{\alpha}} \to v_j$ as $\bar{\alpha} \to 0^+$ for each $j \in [2\dim\mathcal{C}]$.

For each vertex $w_\theta$ of~$Q$ with sign vector~$\theta$, define the 
perturbed vertex
\[
    w_{\theta}^{\bar{\alpha}} 
    := \sum_{i=1}^{\dim\mathcal{C}} \begin{cases}
        \tilde{v}_{2i-1}^{\bar{\alpha}} & \text{if } \theta_i = 1,\\[4pt]
        \tilde{v}_{2i}^{\bar{\alpha}} & \text{if } \theta_i = -1.
    \end{cases}
\]
Each summand lies in $\pos{\bar{V}_{\bar{\alpha}}}$, 
so $w_{\theta}^{\bar{\alpha}} \in \pos{\bar{V}_{\bar{\alpha}}}$. Define
\[
    Q_{\bar{\alpha}} 
    := \conv\!\left\{ w_{\theta}^{\bar{\alpha}} \;\middle|\; 
       \theta \in \{-1,1\}^{\dim\mathcal{C}} \right\}.
\]
Then $Q_{\bar{\alpha}} \subseteq \pos{\bar{V}_{\bar{\alpha}}}$, since
$\pos{\bar{V}_{\bar{\alpha}}}$ is convex. We bound the vertex displacement uniformly in $\zeta$. Let $Z := \conv(U\cup T)$
be the compact set in which the shared direction $\zeta$ ranges.
The vector $\zeta$ is a convex combination of $U\cup T$, so $\zeta\in Z$.
Every vector in $U\cup T$ has unit norm and $\mathcal{B}[0,1]$ is convex, hence $Z\subseteq\mathcal{B}[0,1]$ and
$\norm{\zeta}\le 1$ for every $\zeta\in Z$. The set $Z$ does not depend on $\bar{\alpha}$.
For each $j \in [2\dim\mathcal{C}]$, writing $v_j \in U$ for the corresponding noiseless
direction,
\begin{equation*}
    \tilde{v}_j^{\bar{\alpha}} - v_j = \bar{\alpha}\bar{\eta}\,(\zeta - v_j),
    \qquad \text{so} \qquad
    \norm{\tilde{v}_j^{\bar{\alpha}} - v_j} = \bar{\alpha}\bar{\eta}\,\norm{\zeta - v_j}
    \le 2\bar{\alpha}\bar{\eta},
\end{equation*}
where the last inequality uses $\norm{\zeta} \le 1$ and $\norm{v_j} = 1$ (so
$\norm{\zeta - v_j} \le 2$), and holds for every $\zeta \in Z$. Each vertex
$w_\theta^{\bar{\alpha}}$ is a sum of $\dim\mathcal{C}$ such directions and $w_\theta$ the sum of the
matching $v_j$, so by the triangle inequality
\begin{equation*}
    \norm{w_\theta^{\bar{\alpha}} - w_\theta} \le 2\dim\mathcal{C}\bar{\alpha}\bar{\eta}
    \qquad \text{for every } \theta \in \{-1,1\}^{\dim\mathcal{C}} \text{ and every } \zeta \in Z.
\end{equation*}
As $Q$ and $Q_{\bar{\alpha}}$ are the convex hulls of the matched vertex families
$\{w_\theta\}_\theta$ and $\{w_\theta^{\bar{\alpha}}\}_\theta$,
\cref{lem:hausdorff_matched_convex_hulls} gives
\begin{equation}\label{eq:interior_uniform_hausdorff}
    d_H(Q_{\bar{\alpha}}, Q) \le 2\dim\mathcal{C}\bar{\alpha}\bar{\eta}
    \qquad \text{for every } \zeta \in Z.
\end{equation}
Each perturbed vertex $w_\theta^{\bar{\alpha}}$ is a sum of vectors
$\tilde{v}_j^{\bar{\alpha}} = \bar{\alpha}\bar{\eta}\zeta + (1-\bar{\alpha}\bar{\eta})v_j$,
which are affine combinations of $\zeta, v_j \in \aff{\mathcal{C}-\bar{x}}$ and hence
lie in the subspace $\aff{\mathcal{C}-\bar{x}} = \aff{Q}$; therefore
$Q_{\bar{\alpha}} \subseteq \aff{Q}$, so \cref{lem:containment_under_hausdorff} applies.

Set $\tilde{\alpha} := \min\!\left\{ 1,\ \dist(S, \rb{Q}) / (2\dim\mathcal{C}\bar{\eta}) \right\} > 0$. This
depends only on $\dim\mathcal{C}$, $\bar{\eta}$, and the noiseless gap $\dist(S, \rb{Q})$
supplied by \cref{lem:interior_gap}, and is therefore independent of $\zeta$ (the cap
by $1$ keeps $\bar{\alpha}$ within the admissible noise range). For every
$\bar{\alpha} \in (0, \tilde{\alpha})$ and every $\zeta \in Z$,
\eqref{eq:interior_uniform_hausdorff} yields $d_H(Q_{\bar{\alpha}}, Q) < \dist(S, \rb{Q})$,
so $S \subseteq Q_{\bar{\alpha}}$ by \cref{lem:containment_under_hausdorff}. For any
$d \in D_\mathcal{C}(\bar{x})$, there exists $\lambda > 0$ with
$\lambda d \in S \subseteq Q_{\bar{\alpha}} \subseteq \pos{\bar{V}_{\bar{\alpha}}}$,
so $d \in \pos{\bar{V}_{\bar{\alpha}}}$ since $\pos{\bar{V}_{\bar{\alpha}}}$ is a cone.
Because $\tilde{\alpha}$ does not depend on $\zeta$, the conclusion holds for every
realization of the shared direction.

\medskip
\textbf{Case 2: $\bar{x} \in \rb{\mathcal{C}}$, $\mathcal{C}$ smooth at $\bar{x}$.}
By \cref{lem:boundary_gap}, there exist a compact set 
$S \subseteq D_\mathcal{C}^{\beta\eta}(\bar{x})$ and a polytope 
$P \subseteq \pos{\bar{V}_0}$ such that $S \subseteq \ri{P}$, 
$\dist(S, \rb{P}) > 0$, and every direction in $D_\mathcal{C}^{\beta\eta}(\bar{x})$ 
can be scaled into~$S$.

Let $u \in \aff{\mathcal{C} - \bar{x}}$ be the unit vector from 
\eqref{eq:closure_of_feasible_directions}. By the proof of \cref{lem:boundary_gap} 
(specifically \eqref{eq:truncated_spanning_set_is_polytope}), the polytope~$P$ 
satisfies
\[
    P = \conv\!\left( \{0\} \cup \left\{ \frac{v}{\langle v, u \rangle} 
    \;\middle|\; v \in \bar{V}_0 \right\} \right).
\]
For any $a$ with $\inner{a}{u}>0$, write
$R_u(a):=a/\inner{a}{u}$. This map depends only on the positive ray through $a$:
$R_u(\lambda a)=R_u(a)$ for every $\lambda>0$.
At a smooth boundary point the third branch of \cref{alg:construct_directions} sets $U=\{u\}$ and $T=\{\pm w_1,\ldots,\pm w_{\dim\mathcal{C}-1}\}$, where $\{u,w_1,\ldots,w_{\dim\mathcal{C}-1}\}$ is an orthonormal basis of $\aff{\mathcal{C}-\bar{x}}$.
Write $m:=|U\cup T|$ and $V:=U\cup T=\{u, w_1, -w_1, \ldots, w_{\dim\mathcal{C}-1}, -w_{\dim\mathcal{C}-1}\}$.
For each $\hat{v}_i\in V$, write $v_i$ for the corresponding unnormalized noiseless generator: $u$ itself if $\hat{v}_i\in U$, and $\bar{\eta}\,u+(1-\bar{\eta})\,\hat{v}_i$ if $\hat{v}_i\in T$.
These differ from the unit vectors in $\bar{V}_0$ only by positive scaling, so
\[
    P=\conv\!\left(\{0\}\cup\{R_u(v_i)\mid i\in[m]\}\right).
\]
Moreover, $\norm{v_i}\le 1$ and $\inner{v_i}{u}\ge\bar{\eta}$ for every $i$
(the latter inner product equals $1$ when $\hat{v}_i\in U$ and $\bar{\eta}$
when $\hat{v}_i\in T$).

For the noisy construction we work with the unnormalized combinations of \cref{alg:construct_directions}, enumerated so that $\tilde{v}_i^{\bar{\alpha}}$ is the perturbation of~$v_i$:
\begin{itemize}
    \item For $\hat{v}_i \in U$: \; 
    $\tilde{v}_i^{\bar{\alpha}} 
    = \bar{\alpha}\bar{\eta}\zeta + (1 - \bar{\alpha}\bar{\eta})\,\hat{v}_i$.
    \item For $\hat{v}_i \in T$: \; 
    $\tilde{v}_i^{\bar{\alpha}} 
    = (1-\bar{\alpha})\bar{\eta}\,u + \bar{\alpha}\bar{\eta}\zeta 
      + (1 - \bar{\eta})\,\hat{v}_i$.
\end{itemize}
Let $Z := \conv(V)$ be the compact set in which the shared direction $\zeta$
ranges. Since every element of $V$ has unit norm, $Z\subseteq\mathcal{B}[0,1]$;
in particular, $\norm{\zeta}\le 1$ uniformly over $\zeta\in Z$, and $Z$ does
not depend on $\bar{\alpha}$. In both branches above the
displacement has the form $\tilde{v}_i^{\bar{\alpha}} - v_i = \bar{\alpha}\bar{\eta}\,(\zeta - c_i)$
with $c_i \in \{\hat{v}_i, u\}$ satisfying $\norm{c_i}=1$, hence
\begin{equation}\label{eq:boundary_uniform_displacement}
    \norm{\tilde{v}_i^{\bar{\alpha}} - v_i} = \bar{\alpha}\bar{\eta}\,\norm{\zeta - c_i}
    \le 2\bar{\alpha}\bar{\eta}
    \qquad \text{for every } \zeta \in Z.
\end{equation}
Combining this estimate with $\inner{v_i}{u}\ge\bar{\eta}$ and $\norm{u}=1$ gives
\begin{equation}\label{eq:boundary_uniform_denominator}
    \inner{\tilde{v}_i^{\bar{\alpha}}}{u}
    \ge \inner{v_i}{u} - \norm{\tilde{v}_i^{\bar{\alpha}} - v_i}
    \ge \bar{\eta} - 2\bar{\alpha}\bar{\eta}
    = \bar{\eta}(1 - 2\bar{\alpha})
    \ge \tfrac{\bar{\eta}}{2}
\end{equation}
for every $\bar{\alpha} \le \tfrac{1}{4}$ and every $\zeta \in Z$; in particular
$\inner{\tilde{v}_i^{\bar{\alpha}}}{u} > 0$. For such $\bar{\alpha}$, define the
perturbed polytope
\[
    P_{\bar{\alpha}} 
    := \conv\!\left( \{0\} \cup \left\{ 
       \frac{\tilde{v}_i^{\bar{\alpha}}}%
            {\langle \tilde{v}_i^{\bar{\alpha}}, u \rangle} 
       \;\middle|\; i \in [m] \right\} \right).
\]
For each $i$, \cref{alg:construct_directions} includes
$\mathrm{normalize}(\tilde{v}_i^{\bar{\alpha}})$ in
$\bar{V}_{\bar{\alpha}}$. Since $R_u(\tilde{v}_i^{\bar{\alpha}})$ lies on the
same positive ray, convexity gives
$P_{\bar{\alpha}}\subseteq\pos{\bar{V}_{\bar{\alpha}}}$.

We now bound the vertex displacement uniformly in $\zeta \in Z$. For any $a, b$ with
$\inner{a}{u}, \inner{b}{u} > 0$,
\begin{equation*}
    \frac{a}{\inner{a}{u}} - \frac{b}{\inner{b}{u}}
    = \frac{a\,\inner{b - a}{u} + (a - b)\,\inner{a}{u}}
           {\inner{a}{u} \inner{b}{u}},
    \qquad \text{so} \qquad
    \norm{\frac{a}{\inner{a}{u}} - \frac{b}{\inner{b}{u}}}
    \le \frac{2\norm{a}}{\inner{a}{u} \inner{b}{u}}\,\norm{a - b},
\end{equation*}
using $|\inner{\cdot}{u}| \le \norm{\cdot}$. Apply this with
$a = \tilde{v}_i^{\bar{\alpha}}$ and $b = v_i$. Since $\norm{v_i}\le 1$, for
$\bar{\alpha} \le \tfrac{1}{4}$,
$\norm{\tilde{v}_i^{\bar{\alpha}}} \le \norm{v_i} + 2\bar{\alpha}\bar{\eta} \le 2$, and
\eqref{eq:boundary_uniform_denominator} gives $\inner{\tilde{v}_i^{\bar{\alpha}}}{u}
\ge \tfrac{\bar{\eta}}{2}$, while $\inner{v_i}{u} \ge \bar{\eta}$; combined with
\eqref{eq:boundary_uniform_displacement},
\begin{equation*}
    \norm{\frac{\tilde{v}_i^{\bar{\alpha}}}{\inner{\tilde{v}_i^{\bar{\alpha}}}{u}}
            - \frac{v_i}{\inner{v_i}{u}}}
    \le \frac{2 \cdot 2}{(\bar{\eta}/2)\,\bar{\eta}} \cdot 2\bar{\alpha}\bar{\eta}
    = \frac{16\,\bar{\alpha}}{\bar{\eta}}
    \qquad \text{for every } \zeta \in Z.
\end{equation*}
The vertex $0$ is common to $P$ and $P_{\bar{\alpha}}$, so all matched vertices are
within $16\bar{\alpha}/\bar{\eta}$ of one another. By
\cref{lem:hausdorff_matched_convex_hulls},
\begin{equation}\label{eq:boundary_uniform_hausdorff}
    d_H(P_{\bar{\alpha}}, P) \le \frac{16\,\bar{\alpha}}{\bar{\eta}}
    \qquad \text{for every } \zeta \in Z.
\end{equation}
Each vertex of $P_{\bar{\alpha}}$ is either $0$ or a positive multiple of
$\tilde{v}_i^{\bar{\alpha}}$, and $\tilde{v}_i^{\bar{\alpha}}$ is an affine combination of
$\zeta, \hat{v}_i, u \in \aff{\mathcal{C}-\bar{x}}$; hence
$P_{\bar{\alpha}} \subseteq \aff{\mathcal{C}-\bar{x}} = \aff{P}$, so
\cref{lem:containment_under_hausdorff} applies.

Set $\tilde{\alpha} := \min\!\left\{ \tfrac{1}{4},\ \dfrac{\bar{\eta}\,\dist(S, \rb{P})}{16} \right\} > 0$.
This depends only on $\bar{\eta}$ and the noiseless gap $\dist(S, \rb{P})$ supplied by
\cref{lem:boundary_gap}, and is therefore independent of $\zeta$. For every
$\bar{\alpha} \in (0, \tilde{\alpha})$ and every $\zeta \in Z$,
\eqref{eq:boundary_uniform_hausdorff} yields $d_H(P_{\bar{\alpha}}, P) < \dist(S, \rb{P})$,
so $S \subseteq P_{\bar{\alpha}}$ by \cref{lem:containment_under_hausdorff}. For any
$d \in D_\mathcal{C}^{\beta\eta}(\bar{x})$, there exists $\lambda > 0$ with
$\lambda d \in S \subseteq P_{\bar{\alpha}} \subseteq \pos{\bar{V}_{\bar{\alpha}}}$,
so $d \in \pos{\bar{V}_{\bar{\alpha}}}$ since $\pos{\bar{V}_{\bar{\alpha}}}$ is a cone.
Because $\tilde{\alpha}$ does not depend on $\zeta$, the conclusion holds for every
realization of the shared direction.
\end{proof}

\cref{lem:tilde_v_noised_spans_beta_eta_feasible} asserts that a positive
noise threshold exists. We now compute the two gap distances on which it depends; this gives explicit thresholds in terms of $\eta$, $\beta$ and $\dim\mathcal{C}$. We first record an
elementary fact used in the boundary computation.

\begin{lem}[inradius of a polyhedron]\label{lem:inradius_polyhedron}
    Let $\nu_1,\dots,\nu_N\in\real^m\setminus\{0\}$, $b_1,\dots,b_N\in\real$, and
    $C=\bigcap_{i=1}^N\{x\in\real^m:\inner{\nu_i}{x}\le b_i\}$. If $C$ is closed
    with nonempty interior, then for every $x\in\interior{C}$,
    \begin{equation*}
        \dist\bigl(x,\rb{C}\bigr)=\min_{i\in[N]}\frac{b_i-\inner{\nu_i}{x}}{\norm{\nu_i}}.
    \end{equation*}
\end{lem}
\begin{proof}
    For $x\in\interior{C}$ the distance to the boundary equals the largest
    inscribed-ball radius, $\dist(x,\rb{C})=\sup\{r\ge0:\mathcal{B}[x,r]\subseteq C\}$.
    For each $i$, $\mathcal{B}[x,r]\subseteq\{\inner{\nu_i}{\cdot}\le b_i\}$ iff
    $\inner{\nu_i}{x}+r\norm{\nu_i}\le b_i$, i.e.\ iff
    $r\le(b_i-\inner{\nu_i}{x})/\norm{\nu_i}$. Intersecting over $i$ and taking the
    supremum over $r$ gives the claim.
\end{proof}

\begin{proposition}[interior gap value]\label{prop:interior_gap_value}
    In the setting of \cref{lem:interior_gap}, $\dist\bigl(S,\rb{Q}\bigr)=\tfrac12$.
\end{proposition}
\begin{proof}
    Let $\Phi:\real^{\dim\mathcal{C}}\to L$, $\Phi(z)=\sum_{i=1}^{\dim\mathcal{C}} z_i w_i$; since $\{w_i\}$ is
    orthonormal, $\Phi$ is a linear isometry onto $L:=\aff{\mathcal{C}-\bar x}$. The cube
    $\mathcal{B}_\infty[0,1]=[-1,1]^{\dim\mathcal{C}}$ is the convex hull of $\{-1,1\}^{\dim\mathcal{C}}$, and $\Phi$
    preserves convex hulls, so $Q=\Phi(\mathcal{B}_\infty[0,1])$; as $\Phi$ is a
    homeomorphism onto $L$, $\rb{Q}=\Phi(\{z:\norm{z}_\infty=1\})$. Since
    $\bar x\in\ri{\mathcal{C}}$, \cref{lem:feasible_directions_at_relative_interior} gives $D_\mathcal{C}(\bar x)=L$,
    hence $S=\mathcal{B}[0,\tfrac12]\cap L=\Phi(\mathcal{B}_2[0,\tfrac12])$.

    For $s=\Phi(z_s)\in S$ and $v=\Phi(z_v)\in\rb{Q}$ we have $\norm{z_s}_2\le\tfrac12$ and
    $\norm{z_v}_2\ge\norm{z_v}_\infty=1$, so $\norm{s-v}=\norm{z_s-z_v}_2\ge1-\tfrac12=\tfrac12$;
    thus $\dist(S,\rb{Q})\ge\tfrac12$. Equality holds at $s=\Phi(\tfrac12 e_1)$, $v=\Phi(e_1)$.
\end{proof}

\begin{proposition}[boundary gap value]\label{prop:boundary_gap_value}
    In the setting of \cref{lem:boundary_gap}, set $\vartheta:=\tfrac{1-\bar\eta}{\bar\eta}$.
    \begin{enumerate}
        \item If $\dim\mathcal{C}\ge 2$, then
        $\displaystyle \dist\bigl(S,\rb{P}\bigr)
        =\frac{\vartheta\beta\eta-\sqrt{(\dim\mathcal{C}-1)(1-\beta^2\eta^2)}}{4\sqrt{\vartheta^2+\dim\mathcal{C}-1}}$.
        \item If $\dim\mathcal{C}=1$, then $\dist\bigl(S,\rb{P}\bigr)=\tfrac14$.
        \item In all cases, $\displaystyle \dist\bigl(S,\rb{P}\bigr)\ge\frac{(\beta-1)\,\eta}{4\sqrt{\dim\mathcal{C}}}$.
    \end{enumerate}
\end{proposition}
\begin{proof}
    We work throughout in the setting of \cref{lem:boundary_gap}: $\eta\in(0,1)$ and
    $\beta\in(1,1/\eta)$, hence $\beta\eta\in(\eta,1)$ and $1-\beta^2\eta^2>0$. Write
    $L:=\aff{\mathcal{C}-\bar{x}}$ and let $\{u,w_1,\dots,w_{\dim\mathcal{C}-1}\}$ be the
    orthonormal basis of $L$ from the third branch of \cref{alg:construct_directions}, with
    $u$ the unit vector of \eqref{eq:closure_of_feasible_directions}. The two objects furnished
    by \cref{lem:boundary_gap} are the compact set
    $S=D_\mathcal{C}^{\beta\eta}(\bar{x})\cap\{d:0.25\le\norm{d}\le 0.5\}$ and the truncated
    polytope $P=\pos{\bar{V}_0}\cap H_{u,-1}$. As in \cref{alg:construct_directions},
    $\bar\eta=\eta/(2\sqrt{\dim\mathcal{C}})$. Combined with the standing definition
    $\vartheta=(1-\bar\eta)/\bar\eta$, this gives
    \begin{equation*}
        \vartheta=\frac{1-\bar\eta}{\bar\eta}=\frac{2\sqrt{\dim\mathcal{C}}-\eta}{\eta},
        \qquad \vartheta\ge 1.
    \end{equation*}

    Coordinatize $L$ by the linear isometry $\Psi:\real\times\real^{\dim\mathcal{C}-1}\to L$,
    $\Psi(c,a)=c\,u+\sum_{j=1}^{\dim\mathcal{C}-1}a_j w_j$. Distances are preserved, so we compute in
    $\real^{\dim\mathcal{C}}$, writing points as $(c,a)$ with $c=\inner{d}{u}$ and
    $\norm{(c,a)}^2=c^2+\norm{a}_2^2$.

    Recall from the $\alpha=0$ construction of \cref{alg:construct_directions} at a smooth
    boundary point (cf.\ \cref{lem:tilde_v_spans_eta_feasible}) that the non-exploration
    output is obtained by applying $\mathrm{normalize}$ to each vector in
    $\{u\}\cup\{\pm(1-\bar\eta)w_j+\bar\eta u:j\in[\dim\mathcal{C}-1]\}$.
    Since $R_u(v)=v/\inner{v}{u}$ depends only on the positive ray through $v$, the truncated
    polytope $P$ of \cref{lem:boundary_gap} is unchanged if each generator is replaced by a
    positive multiple, and we compute with this unnormalized set.
    Scaling each direction by its $u$-component, for $v=u$ we have $\inner{v}{u}=1$, so
    $R_u(v)=u=(1,\mathbf{0})$. For $v=\pm(1-\bar\eta)w_j+\bar\eta u$, orthonormality
    gives $\inner{v}{u}=\bar\eta$, hence
    $R_u(v)=u\pm\tfrac{1-\bar\eta}{\bar\eta}w_j=(1,\pm\vartheta e_j)$. Thus the vertices
    $R_u(v)$ are $(1,\mathbf{0})$ and $(1,\pm\vartheta e_j)$, $j\in[\dim\mathcal{C}-1]$.
    Writing $K:=\{(c,a):\norm{a}_1\le\vartheta c\}$ for the untruncated cone, we have
    $K\subseteq\{c\ge 0\}$ and therefore
    \begin{equation*}
        P=K\cap\{c\le 1\}=\{(c,a):0\le c\le 1,\ \norm{a}_1\le\vartheta c\}.
    \end{equation*}
    The inclusion $\subseteq$ holds as the right-hand set is convex and contains every vertex;
    for $\supseteq$, a point with $c>0$ and $\norm{a}_1\le\vartheta c$ has
    $a/(\vartheta c)\in\conv\{\pm e_j\}$, so $(1,a/c)\in\conv\{(1,\pm\vartheta e_j)\}$ and
    $(c,a)=c\,(1,a/c)+(1-c)(0,\mathbf{0})$, while $c=0$ forces the apex.

    We next write $S$ in these coordinates. For $d\ne 0$, since $c=\inner{d}{u}$ and
    $\norm{d}^2=c^2+\norm{a}_2^2$, the defining inequality of $D_\mathcal{C}^{\beta\eta}(\bar{x})$
    reads $c/\sqrt{c^2+\norm{a}_2^2}\ge\beta\eta$. As $\beta\eta\in(\eta,1)\subseteq(0,1)$ while
    the denominator is positive, the inequality forces $c>0$. Set
    \begin{equation*}
        \mu:=\frac{\beta\eta}{\sqrt{1-\beta^2\eta^2}}.
    \end{equation*}
    Both sides being positive, squaring is an equivalence, and $1-\beta^2\eta^2>0$ because
    $\beta<1/\eta$, so
    \begin{align*}
        \inner{d/\norm{d}}{u}\ge\beta\eta
        &\iff \frac{c}{\sqrt{c^2+\norm{a}_2^2}}\ge\beta\eta \\
        &\iff c^2(1-\beta^2\eta^2)\ge\beta^2\eta^2\norm{a}_2^2 \\
        &\iff c\ge\mu\norm{a}_2.
    \end{align*}
    Every step reverses on $\{d\ne 0\}$. The definition of $\mu$ also gives
    \begin{equation*}
        1+\mu^2=\frac{1}{1-\beta^2\eta^2},\qquad
        \frac{\mu}{\sqrt{1+\mu^2}}=\beta\eta,\qquad
        \frac{1}{\sqrt{1+\mu^2}}=\sqrt{1-\beta^2\eta^2}.
    \end{equation*}
    Thus
    \begin{equation*}
        S=\bigl\{(c,a):c\ge\mu\norm{a}_2,\ \frac14\le\sqrt{c^2+\norm{a}_2^2}\le\frac12\bigr\}.
    \end{equation*}

    The cone $K$ is a polyhedron. Writing
    $\norm{a}_1=\max_{\epsilon\in\{-1,1\}^{\dim\mathcal{C}-1}}\inner{\epsilon}{a}$ expresses $K$ as
    the halfspace intersection
    \begin{equation*}
        K=\bigcap_{\epsilon\in\{-1,1\}^{\dim\mathcal{C}-1}}\{\inner{\epsilon}{a}-\vartheta c\le 0\},
    \end{equation*}
    with facet normals $\nu_\epsilon=(-\vartheta,\epsilon)$ and right-hand sides $0$. Each normal
    has $\norm{\nu_\epsilon}=\sqrt{\vartheta^2+\norm{\epsilon}^2}$, hence
    $\norm{\nu_\epsilon}=\sqrt{\vartheta^2+\dim\mathcal{C}-1}$, and
    $(1,\mathbf{0})\in\interior{K}$ since $0<\vartheta\cdot 1$. In particular $K$ and $P$ are
    full-dimensional in these coordinates, so $\rb{K}=\bdr{K}$ and $\rb{P}=\bdr{P}$ throughout.
    By \cref{lem:boundary_gap}(1), $S\subseteq\ri{P}$ and $\aff{P}=L$, hence in $\Psi$-coordinates
    $S\subseteq\{0<c<1,\ \norm{a}_1<\vartheta c\}\subseteq\interior{K}$. Applying
    \cref{lem:inradius_polyhedron} at every $(c,a)\in\interior{K}$ and minimizing the per-facet
    residual over $\epsilon$,
    \begin{equation}\label{eq:dist_to_K}
        \dist\bigl((c,a),\rb{K}\bigr)
        =\min_{\epsilon}\frac{\vartheta c-\inner{\epsilon}{a}}{\sqrt{\vartheta^2+\dim\mathcal{C}-1}}
        =\frac{\vartheta c-\max_{\epsilon}\inner{\epsilon}{a}}{\sqrt{\vartheta^2+\dim\mathcal{C}-1}}
        =\frac{\vartheta c-\norm{a}_1}{\sqrt{\vartheta^2+\dim\mathcal{C}-1}}.
    \end{equation}

    We first evaluate $\inf_{s\in S}\dist(s,\rb{K})$, which \eqref{eq:dist_to_K} reduces to a
    scaled infimum of $\vartheta c-\norm{a}_1$ over $S$. Write $r=\norm{(c,a)}\in[\tfrac14,\tfrac12]$.
    Since $r^2=c^2+\norm{a}_2^2$ and $c\ge\mu\norm{a}_2\ge 0$, we have $c/r\in[0,1]$. Set
    $\psi:=\arccos(c/r)$, so $\psi\in[0,\tfrac\pi2]$ and $c=r\cos\psi$. Then $\sin\psi\ge 0$ and
    \begin{equation*}
        \sin\psi=\sqrt{1-\cos^2\psi}=\sqrt{1-\frac{c^2}{r^2}}=\frac{\norm{a}_2}{r},
    \end{equation*}
    hence $\norm{a}_2=r\sin\psi$. The cone constraint $c\ge\mu\norm{a}_2$ becomes $\cot\psi\ge\mu$,
    i.e.\ $\psi\le\psi_0:=\arctan(1/\mu)$. Bounding $\norm{a}_1\le\sqrt{\dim\mathcal{C}-1}\,\norm{a}_2$
    by Cauchy--Schwarz,
    \begin{equation*}
        \vartheta c-\norm{a}_1\ge r\bigl(\vartheta\cos\psi-\sqrt{\dim\mathcal{C}-1}\sin\psi\bigr)=:r\,\phi(\psi).
    \end{equation*}
    Since $\phi'(\psi)=-\vartheta\sin\psi-\sqrt{\dim\mathcal{C}-1}\cos\psi\le 0$ on $[0,\tfrac\pi2]$,
    $\phi$ is nonincreasing and is minimized over $[0,\psi_0]$ at $\psi_0$. From
    $\tan\psi_0=1/\mu$ and $\sin^2\psi_0+\cos^2\psi_0=1$ we have
    \begin{equation*}
        \cos\psi_0=\frac{\mu}{\sqrt{1+\mu^2}},\qquad
        \sin\psi_0=\frac{1}{\sqrt{1+\mu^2}},
    \end{equation*}
    so
    \begin{equation*}
        \phi(\psi_0)=\frac{\vartheta\mu-\sqrt{\dim\mathcal{C}-1}}{\sqrt{1+\mu^2}}.
    \end{equation*}

    The numerator is nonnegative: substituting the identity for $\vartheta$ from the opening,
    \begin{equation*}
        \vartheta\mu=\frac{\beta\bigl(2\sqrt{\dim\mathcal{C}}-\eta\bigr)}{\sqrt{1-\beta^2\eta^2}}
        \ge\beta\bigl(2\sqrt{\dim\mathcal{C}}-\eta\bigr)
        >2\sqrt{\dim\mathcal{C}}-1
        \ge\sqrt{\dim\mathcal{C}}
        \ge\sqrt{\dim\mathcal{C}-1},
    \end{equation*}
    the first inequality as $\sqrt{1-\beta^2\eta^2}\le 1$, the second as $\beta>1$ and $\eta<1$,
    and the last two as $\dim\mathcal{C}\ge 1$. The same chain yields
    \begin{equation}\label{eq:varkappa_beta_eta}
        \vartheta\beta\eta=\beta\bigl(2\sqrt{\dim\mathcal{C}}-\eta\bigr)
        >\sqrt{\dim\mathcal{C}}
        >\sqrt{\dim\mathcal{C}-1},
    \end{equation}
    used for the lower bound in~(3) below. Thus $\phi(\psi_0)\ge 0$.

    Because $r\ge\tfrac14$ and $\phi(\psi_0)\ge 0$, we get
    $\vartheta c-\norm{a}_1\ge\tfrac14\phi(\psi_0)$. For $\dim\mathcal{C}\ge 2$ this bound is
    attained, and that is the only place the assumption $\dim\mathcal{C}\ge 2$ is used in the evaluation of
    $\inf_{s\in S}\dist(s,\rb{K})$. Explicitly, take $r=\tfrac14$, $\psi=\psi_0$, and $a$ with
    equal-magnitude entries: the resulting point sits on the cone boundary $c=\mu\norm{a}_2$
    (as $\psi=\psi_0$), saturates $\norm{a}_1=\sqrt{\dim\mathcal{C}-1}\,\norm{a}_2$ (equal
    magnitudes attain Cauchy--Schwarz), and lies in $S$ since $r=\tfrac14\in[\tfrac14,\tfrac12]$
    and $c=\mu\norm{a}_2\ge 0$ places it in the closed cone $D_\mathcal{C}^{\beta\eta}(\bar{x})$.
    Therefore
    \begin{equation}\label{eq:inf_dist_rbK}
        \inf_{s\in S}\dist(s,\rb{K})
        =\frac{\phi(\psi_0)}{4\sqrt{\vartheta^2+\dim\mathcal{C}-1}}\le\frac14.
    \end{equation}
    For the upper bound, dropping the nonnegative term gives
    $\phi(\psi_0)\le\vartheta\mu/\sqrt{1+\mu^2}$, and
    \begin{equation*}
        \sqrt{(1+\mu^2)(\vartheta^2+\dim\mathcal{C}-1)}\ge\vartheta\mu
    \end{equation*}
    (square both sides: $\vartheta^2+(\dim\mathcal{C}-1)(1+\mu^2)\ge 0$), so the ratio is at most
    $\tfrac14$.

    It remains to identify $\dist(S,\rb{P})$ with $\inf_{s\in S}\dist(s,\rb{K})$ when
    $\dim\mathcal{C}\ge 2$. Since $P=K\cap\{c\le 1\}$, the elementary inclusion
    $\bdr{A\cap B}\subseteq(\bdr{A}\cap B)\cup(A\cap\bdr{B})$ gives
    $\rb{P}\subseteq(\rb{K}\cap\{0\le c\le 1\})\cup\{c=1\}$. Conversely, a point of $\rb{K}$ with
    first coordinate in $[0,1]$ lies on $\rb{P}$: if $p\in\rb{K}$ and $c_p\in[0,1]$, then
    $p\in K$ and $c_p\le 1$, so $p\in P$. Every neighbourhood of $p$ contains a point outside
    $K$, and that point lies outside $P$ as well, since $P\subseteq K$. Thus $p\in\rb{P}$, and
    $\rb{K}\cap\{0\le c\le 1\}\subseteq\rb{P}$.

    Fix $s=(c_s,a_s)\in S$, and let $p\in\rb{K}$ be a nearest point to $s$. We check that the
    first coordinate of $p$ lies in $[0,1]$, so that $p\in\rb{P}$. Since $c_s\le\norm{s}\le\tfrac12$
    and $p$ realizes the inradius distance \eqref{eq:dist_to_K}, writing
    $\delta:=\vartheta c_s-\norm{a_s}_1>0$ (as $S\subseteq\interior{K}$),
    \begin{equation*}
        \norm{p-s}=\frac{\delta}{\sqrt{\vartheta^2+\dim\mathcal{C}-1}}\le c_s,
    \end{equation*}
    using $\delta\le\vartheta c_s$ and $\sqrt{\vartheta^2+\dim\mathcal{C}-1}\ge\vartheta$. Hence
    the first coordinate of $p$ obeys $|c_p-c_s|\le\norm{p-s}\le c_s$, i.e.\
    $c_p\in[0,2c_s]\subseteq[0,1]$ (as $c_s\le\tfrac12$), so
    $p\in\rb{K}\cap\{0\le c\le 1\}\subseteq\rb{P}$ and $\dist(s,\rb{P})\le\dist(s,\rb{K})$.
    Conversely, the distance from $s$ to the cap $\{c=1\}$ is $1-c_s\ge\tfrac12$, and
    $\rb{K}\cap\{0\le c\le 1\}\subseteq\rb{K}$ only increases distance, so
    \begin{equation*}
        \dist(s,\rb{P})\ge\min\bigl\{\dist(s,\rb{K}),\,1-c_s\bigr\}
        \ge\min\bigl\{\dist(s,\rb{K}),\,\frac12\bigr\}.
    \end{equation*}
    Taking infima and using the pointwise bound $\dist(s,\rb{P})\le\dist(s,\rb{K})$ together with
    $\inf_{s\in S}\dist(s,\rb{K})\le\tfrac14$ \eqref{eq:inf_dist_rbK},
    \begin{align*}
        \inf_{s\in S}\dist(s,\rb{P})
        &\le\inf_{s\in S}\dist(s,\rb{K}),\\
        \inf_{s\in S}\dist(s,\rb{P})
        &\ge\min\bigl\{\inf_{s\in S}\dist(s,\rb{K}),\,\frac12\bigr\}
        =\inf_{s\in S}\dist(s,\rb{K}).
    \end{align*}
    Therefore
    \begin{equation}\label{eq:rbP_eq_K}
        \dist(S,\rb{P})=\inf_{s\in S}\dist(s,\rb{K})
        =\frac{1}{\sqrt{\vartheta^2+\dim\mathcal{C}-1}}\,\inf_{s\in S}\bigl(\vartheta c_s-\norm{a_s}_1\bigr).
    \end{equation}

    Plugging the attained value into \eqref{eq:rbP_eq_K},
    \begin{equation*}
        \dist(S,\rb{P})=\frac{1}{\sqrt{\vartheta^2+\dim\mathcal{C}-1}}\cdot\frac14\cdot
        \frac{\vartheta\mu-\sqrt{\dim\mathcal{C}-1}}{\sqrt{1+\mu^2}},
    \end{equation*}
    and using $\tfrac{\mu}{\sqrt{1+\mu^2}}=\beta\eta$,
    $\tfrac{1}{\sqrt{1+\mu^2}}=\sqrt{1-\beta^2\eta^2}$ to clear $\mu$ from the numerator,
    \begin{equation*}
        \frac{\vartheta\mu-\sqrt{\dim\mathcal{C}-1}}{\sqrt{1+\mu^2}}
        =\vartheta\beta\eta-\sqrt{\dim\mathcal{C}-1}\,\sqrt{1-\beta^2\eta^2}
        =\vartheta\beta\eta-\sqrt{(\dim\mathcal{C}-1)(1-\beta^2\eta^2)},
    \end{equation*}
    yields the formula in~(1).

    The case $\dim\mathcal{C}=1$ is computed directly and does not use the identification
    \eqref{eq:rbP_eq_K}. Here $K=\{c\ge 0\}$, $P=[0,1]\cdot u$, $\rb{P}=\{0,u\}$, and
    $S=\{c\,u:\tfrac14\le c\le\tfrac12\}$, so $\dist(S,\rb{P})=\tfrac14$, establishing~(2).
    The case split is necessary: with no transverse coordinate $a$, the minimizing angle
    $\psi_0$ is unavailable, and the expression in~(1) evaluates to $\beta\eta/4$, which is
    strictly smaller than $\tfrac14$.

    For the lower bound in~(3), first take $\dim\mathcal{C}=1$. Then
    $(\beta-1)\eta=\beta\eta-\eta<1-\eta<1$ (using $\beta<1/\eta$), so
    $\tfrac14\ge\tfrac{(\beta-1)\eta}{4}$.

    Now take $\dim\mathcal{C}\ge 2$ and bound formula~(1). The numerator satisfies
    \begin{equation*}
        \vartheta\beta\eta-\sqrt{(\dim\mathcal{C}-1)(1-\beta^2\eta^2)}
        \ge\vartheta\beta\eta-\sqrt{\dim\mathcal{C}-1},
    \end{equation*}
    as $1-\beta^2\eta^2\le 1$; the right-hand side is positive by the estimate
    $\vartheta\beta\eta>\sqrt{\dim\mathcal{C}-1}$ \eqref{eq:varkappa_beta_eta}. For the denominator,
    $\vartheta\ge 1$ as already noted, hence
    \begin{equation*}
        \sqrt{\vartheta^2+\dim\mathcal{C}-1}
        \le\sqrt{\vartheta^2+(\dim\mathcal{C}-1)\vartheta^2}
        =\vartheta\sqrt{\dim\mathcal{C}}.
    \end{equation*}
    Combining,
    \begin{equation*}
        \dist(S,\rb{P})\ge\frac{\vartheta\beta\eta-\sqrt{\dim\mathcal{C}-1}}{4\vartheta\sqrt{\dim\mathcal{C}}}
        =\frac{1}{4\sqrt{\dim\mathcal{C}}}\Bigl(\beta\eta-\frac{\sqrt{\dim\mathcal{C}-1}}{\vartheta}\Bigr).
    \end{equation*}
    Finally
    \begin{equation*}
        \frac{1}{\vartheta}=\frac{\bar\eta}{1-\bar\eta}\le 2\bar\eta=\frac{\eta}{\sqrt{\dim\mathcal{C}}},
    \end{equation*}
    from $1-\bar\eta\ge\tfrac12$, so
    \begin{equation*}
        \frac{\sqrt{\dim\mathcal{C}-1}}{\vartheta}
        \le\eta\sqrt{\frac{\dim\mathcal{C}-1}{\dim\mathcal{C}}}\le\eta.
    \end{equation*}
    Therefore $\beta\eta-\tfrac{\sqrt{\dim\mathcal{C}-1}}{\vartheta}\ge(\beta-1)\eta$ and
    $\dist(S,\rb{P})\ge\tfrac{(\beta-1)\eta}{4\sqrt{\dim\mathcal{C}}}$.

    In both cases the argument gives $0<\dist(S,\rb{P})\le\tfrac14$.
\end{proof}

\begin{corollary}[explicit admissible noise levels]\label{cor:explicit_alpha}
    Recall $\bar\eta=\eta/(2\sqrt{\dim\mathcal{C}})$ as in \cref{alg:construct_directions}.
    \begin{enumerate}
        \item If $\bar x\in\ri{\mathcal{C}}$, the conclusion of
        \cref{lem:tilde_v_noised_spans_beta_eta_feasible} holds whenever
        \begin{equation*}
            0<\bar\alpha<\min\Bigl\{1,\ \tfrac{1}{4\dim\mathcal{C}\bar\eta}\Bigr\}
            =\min\Bigl\{1,\ \tfrac{1}{2\sqrt{\dim\mathcal{C}}\,\eta}\Bigr\}.
        \end{equation*}
        \item If $\bar x\in\rb{\mathcal{C}}$ and $\mathcal{C}$ is smooth at $\bar x$, then for
        every $\beta\in(1,1/\eta)$ the conclusion holds whenever
        $0<\bar\alpha<\min\{\tfrac14,\ \bar\eta\,\dist(S,\rb{P})/16\}$, and in particular whenever
        \begin{equation*}
            0<\bar\alpha<\min\Bigl\{\tfrac14,\ \tfrac{(\beta-1)\,\eta^2}{128\dim\mathcal{C}}\Bigr\}.
        \end{equation*}
    \end{enumerate}
\end{corollary}
\begin{proof}
    (1): the Case~1 bound $d_H(Q_{\bar\alpha},Q)\le 2\dim\mathcal{C}\bar\alpha\bar\eta$ and
    \cref{prop:interior_gap_value} give $d_H(Q_{\bar\alpha},Q)<\tfrac12=\dist(S,\rb{Q})$ for
    $\bar\alpha<\tfrac{1}{4\dim\mathcal{C}\bar\eta}$, whence $S\subseteq Q_{\bar\alpha}$ by
    \cref{lem:containment_under_hausdorff}; the cap at $1$ is $\alpha\in[0,1)$.
    (2): the first range is the Case~2 threshold; the second uses
    \cref{prop:boundary_gap_value}(3), since
    $\bar\eta\,\dist(S,\rb{P})/16\ge\bar\eta(\beta-1)\eta/(64\sqrt{\dim\mathcal{C}})
    =(\beta-1)\eta^2/(128\dim\mathcal{C})$.
\end{proof}

\Cref{lem:optimality_measure_zero_implies_no_descent_in_span} turns a zero value of the PDM measure into the first-order condition $h'(x_k;\,\cdot)\ge 0$ on any cone contained in the positive span of the direction set. \Cref{lem:tilde_v_noised_spans_beta_eta_feasible} shows that this cone may be taken to be $D_{\mathcal{C}}(x_k)$ at a relative interior point and $D_{\mathcal{C}}^{\beta\eta}(x_k)$ at a smooth boundary point. By \cref{cor:explicit_alpha} the thresholds depend only on $\eta$, $\dim\mathcal{C}$ and $\beta$, so a single admissible $\alpha$ serves at every interior point and every smooth boundary point.

\begin{corollary}[PDM certificate for the noised method]\label{cor:noised-pdm-certificate}
    Let $\{x_k\}_{k\ge 0}$ be generated by \cref{alg:GRFD} with $\eta\in(0,1]$ and
    $\alpha\in(0,1)$, fix $k\ge 0$, and suppose
    $\mathcal{O}^\gamma(x_k;\tilde{V}^k,r)=0$.
    \begin{enumerate}
        \item If $x_k\in\ri{\mathcal{C}}$ and $\alpha\in(0,\tilde{\alpha})$ for the threshold
        $\tilde{\alpha}$ of part~1 of \cref{lem:tilde_v_noised_spans_beta_eta_feasible} at $x_k$,
        then
        \begin{equation*}
            h'(x_k;\,d)\ge 0 \qquad \forall d\in D_{\mathcal{C}}(x_k).
        \end{equation*}
        \item If $x_k\in\rb{\mathcal{C}}$, $\mathcal{C}$ is smooth at $x_k$,
        $\beta\in(1,1/\eta)$, and $\alpha\in(0,\tilde{\alpha})$ for the threshold
        $\tilde{\alpha}$ of part~2 of \cref{lem:tilde_v_noised_spans_beta_eta_feasible} at $x_k$,
        then
        \begin{equation*}
            h'(x_k;\,d)\ge 0 \qquad \forall d\in D_{\mathcal{C}}^{\beta\eta}(x_k).
        \end{equation*}
    \end{enumerate}
\end{corollary}
\begin{proof}
    Write $\bar{V}^k := \tilde{V}^k \setminus \{\tilde{v}_{\mathrm{ex}}^k\}$.
    By \cref{lem:GRFD_all_considered_directions_are_feasible},
    $\bar{V}^k \subseteq D_{\mathcal{C}}(x_k)$.
    The measure \eqref{eq:optimality_measure} is a sum of squares, so
    $\mathcal{O}^\gamma(x_k;\tilde{V}^k,r)=0$ implies
    $\mathcal{O}^\gamma(x_k;\bar{V}^k,r)=0$.
    The set $\tilde{V}^k$ is the output of \cref{alg:construct_directions} at $x_k$ with this
    $\eta$ and $\alpha$, so \cref{lem:tilde_v_noised_spans_beta_eta_feasible} applies.

    If $x_k\in\ri{\mathcal{C}}$ and $\alpha$ is as in part~1, then
    $D_{\mathcal{C}}(x_k)\subseteq\pos{\bar{V}^k}$.
    \Cref{lem:optimality_measure_zero_implies_no_descent_in_span} with
    $\mathcal{K}=D_{\mathcal{C}}(x_k)$ yields
    $\langle\nabla f(x_k),d\rangle - g'(x_k;\,d)\ge 0$ for every
    $d\in D_{\mathcal{C}}(x_k)$, which is $h'(x_k;\,d)\ge 0$.

    If $x_k\in\rb{\mathcal{C}}$ with $\mathcal{C}$ smooth at $x_k$ and $\alpha$ is as in part~2,
    then $D_{\mathcal{C}}^{\beta\eta}(x_k)\subseteq\pos{\bar{V}^k}$.
    The same lemma with $\mathcal{K}=D_{\mathcal{C}}^{\beta\eta}(x_k)$ gives the claim.
\end{proof}

\section{Conclusion}
In this paper we developed a feasible-directions framework for \eqref{eq:feasible_directions_beck_hallak} over a general closed convex set, based on $\eta$-feasible directions: feasible directions that remain uniformly separated from the boundary of the feasible direction cone. We proposed two methods that search along finitely many feasible directions whose positive span contains the $\eta$-feasible directions, together with a randomized exploration direction. The first assumes access to an $\eta$-direction oracle and extends the linearly constrained method of \citet{BECK2022517}. The second constructs spanning directions by a noise-augmented procedure and avoids nonsmooth boundary points almost surely. For both methods we established an $O(1/K)$ rate in the proximal directional minimization measure of \citet{BECK2022517}. In addition, by an adaptation of the exploration result of \citet{greenstein2026achieving} to the regularized line search, every accumulation point of either method is almost surely $d$-stationary. Under curvature bounds on a smooth constraint boundary, \cref{alg:Oracle_GRFD} driven by the depth-aware construction of \cref{alg:construct_directions} reaches $(\varepsilon,R)$-$d$-stationarity in $O(\dim\mathcal{C}\cdot\varepsilon^{-2})$ iterations. As a byproduct, we obtained a characterization of $\eta$-feasible directions at interior and smooth boundary points, and directional spanning properties under noise.

\bibliography{sntl_bib}

@article{greenstein2026achieving,
  title={Achieving Directional-Stationarity from a Single Random Direction Step},
  author={Greenstein, Dan and Hallak, Nadav},
  journal={arXiv preprint arXiv:2605.22045},
  year={2026}
}

@article{fan2001variable,
  title={Variable selection via nonconcave penalized likelihood and its oracle properties},
  author={Fan, Jianqing and Li, Runze},
  journal={Journal of the American statistical Association},
  volume={96},
  number={456},
  pages={1348--1360},
  year={2001},
  publisher={Taylor \& Francis}
}

@book{RW98,
  author    = {R. Tyrrell Rockafellar and Roger J.-B. Wets},
  title     = {Variational Analysis},
  series    = {Grundlehren der mathematischen Wissenschaften},
  volume    = {317},
  edition   = {Corrected 3rd printing},
  publisher = {Springer},
  address   = {Berlin, Heidelberg},
  year      = {2009},
  isbn      = {978-3-540-62772-2},
  doi       = {10.1007/978-3-642-02431-3}
}

@book{boyd2004convex,
  title={Convex optimization},
  author={Boyd, Stephen and Vandenberghe, Lieven},
  year={2004},
  publisher={Cambridge university press}
}

@article{Frank1956,
  title={An algorithm for quadratic programming},
  author={Frank, Marguerite and Wolfe, Philip},
  journal={Naval Research Logistics Quarterly},
  volume={3},
  number={1-2},
  pages={95--110},
  year={1956}
}

@inproceedings{Jaggi2013,
  title={Revisiting {F}rank-{W}olfe: Projection-free sparse convex optimization},
  author={Jaggi, Martin},
  booktitle={Proceedings of the 30th International Conference on Machine Learning},
  pages={427--435},
  year={2013},
  organization={PMLR}
}

@article{feasible_directions_2020,
author = {Beck, Amir and Hallak, Nadav},
title = {On the Convergence to Stationary Points of Deterministic and Randomized Feasible Descent Directions Methods},
journal = {SIAM Journal on Optimization},
volume = {30},
number = {1},
pages = {56-79},
year = {2020},
doi = {10.1137/18M1217760},

URL = { 
    
        https://doi.org/10.1137/18M1217760
    
    

},
eprint = { 
    
        https://doi.org/10.1137/18M1217760
    
    

}
}

@article{BECK2022517,
title = {The regularized feasible directions method for nonconvex optimization},
journal = {Operations Research Letters},
volume = {50},
number = {5},
pages = {517-523},
year = {2022},
issn = {0167-6377},
doi = {https://doi.org/10.1016/j.orl.2022.07.005},
url = {https://www.sciencedirect.com/science/article/pii/S0167637722000906},
author = {Amir Beck and Nadav Hallak}
}

@misc{fukuda2016lecture,
  author       = {Komei Fukuda},
  title        = {Lecture: Polyhedral Computation, Spring 2016},
  howpublished = {Lecture notes, Department of Mathematics and Institute of Theoretical Computer Science, {ETH} Z{\"u}rich},
  year         = {2016},
  note         = {\url{https://people.inf.ethz.ch/fukudak/lect/pclect/notes2016/PolyComp2016.pdf}}
}

@Article{PRA17,
  author          = {J.~S.~Pang and M.~Razaviyayn and A.~Alvarado},
  title           = {{Computing B-Stationary Points of Nonsmooth DC Programs}},
  journal         = {Mathematics of Operations Research},
  year            = {2017},
  volume          = {42},
  number          = {1},
  pages           = {95--118},
  month           = {jan},
  issn            = {0364-765X},
 }

@article{brilli2025interior,
  title={An interior point method for nonlinear constrained derivative-free optimization},
  author={Brilli, Andrea and Liuzzi, Giampaolo and Lucidi, Stefano},
  journal={Optimization Methods and Software},
  pages={1--39},
  year={2025},
  publisher={Taylor \& Francis}
}

@article{Fasano2014Linesearch,
author = {Fasano, G. and Liuzzi, G. and Lucidi, S. and Rinaldi, F.},
title = {A Linesearch-Based Derivative-Free Approach for Nonsmooth Constrained Optimization},
journal = {SIAM Journal on Optimization},
volume = {24},
number = {3},
pages = {959-992},
year = {2014},
doi = {10.1137/130940037},
URL = {https://doi.org/10.1137/130940037},
eprint = {https://doi.org/10.1137/130940037}
}

@article{larson2019derivative,
  title={Derivative-free optimization methods},
  author={Larson, Jeffrey and Menickelly, Matt and Wild, Stefan M},
  journal={Acta Numerica},
  volume={28},
  pages={287--404},
  year={2019},
  publisher={Cambridge University Press}
}

@book{conn2009introduction,
  title={Introduction to derivative-free optimization},
  author={Conn, Andrew R and Scheinberg, Katya and Vicente, Luis N},
  year={2009},
  publisher={SIAM}
}

@article{kolda2003optimization,
  title={Optimization by direct search: New perspectives on some classical and modern methods},
  author={Kolda, Tamara G and Lewis, Robert Michael and Torczon, Virginia},
  journal={SIAM review},
  volume={45},
  number={3},
  pages={385--482},
  year={2003},
  publisher={SIAM}
}

@article{Kolda2006,
  title={Stationarity results for generating set search for linearly constrained optimization},
  author={Kolda, Tamara G and Lewis, Robert Michael and Torczon, Virginia},
  journal={SIAM Journal on Optimization},
  volume={17},
  number={4},
  pages={943--968},
  year={2006},
  publisher={SIAM}
}

@article{Lewis1999,
  title={Pattern search algorithms for bound constrained minimization},
  author={Lewis, Robert Michael and Torczon, Virginia},
  journal={SIAM Journal on Optimization},
  volume={9},
  number={4},
  pages={1082--1099},
  year={1999},
  publisher={SIAM}
}

@article{Lewis2000,
  title={Pattern search methods for linearly constrained minimization},
  author={Lewis, Robert Michael and Torczon, Virginia},
  journal={SIAM Journal on Optimization},
  volume={10},
  number={3},
  pages={917--941},
  year={2000},
  publisher={SIAM}
}

@article{phan2024difference,
  title={Difference-of-convex algorithm with extrapolation for nonconvex, nonsmooth optimization problems},
  author={Phan, Duy Nhat and Le Thi, Hoai An},
  journal={Mathematics of Operations Research},
  volume={49},
  number={3},
  pages={1973--1985},
  year={2024},
  publisher={INFORMS}
}

@article{candes2008enhancing,
  title={Enhancing sparsity by reweighted l1 minimization},
  author={Candes, Emmanuel J and Wakin, Michael B and Boyd, Stephen P},
  journal={Journal of Fourier analysis and applications},
  volume={14},
  pages={877--905},
  year={2008},
  publisher={Springer}
}

@article{abbaszadehpeivasti2024rate,
  title={On the rate of convergence of the difference-of-convex algorithm (DCA)},
  author={Abbaszadehpeivasti, Hadi and de Klerk, Etienne and Zamani, Moslem},
  journal={Journal of Optimization Theory and Applications},
  volume={202},
  number={1},
  pages={475--496},
  year={2024},
  publisher={Springer}
}

@Article{Ahn2017,
  author        = {M.~Ahn and J.~S.~Pang and J.~Xin},
  title         = {Difference-of-Convex Learning: Directional Stationarity, Optimality, and Sparsity},
  journal       = {{SIAM} Journal on Optimization},
  year          = {2017},
  volume        = {27},
  number        = {3},
  pages         = {1637--1665},
  month         = {jan},
  __markedentry = {[Nadav Hallak:6]},
  doi           = {10.1137/16m1084754},
  owner         = {Nadav Hallak},
  publisher     = {Society for Industrial {\&} Applied Mathematics ({SIAM})},
  timestamp     = {2018.04.12},
}

@book{ben2001lectures,
  title={Lectures on modern convex optimization: analysis, algorithms, and engineering applications},
  author={Ben-Tal, Aharon and Nemirovski, Arkadi},
  year={2001},
  publisher={SIAM}
}

@article{vandenberghe1996semidefinite,
  title={Semidefinite programming},
  author={Vandenberghe, Lieven and Boyd, Stephen},
  journal={SIAM review},
  volume={38},
  number={1},
  pages={49--95},
  year={1996},
  publisher={SIAM}
}

@article{goemans1995improved,
  title={Improved approximation algorithms for maximum cut and satisfiability problems using semidefinite programming},
  author={Goemans, Michel X and Williamson, David P},
  journal={Journal of the ACM (JACM)},
  volume={42},
  number={6},
  pages={1115--1145},
  year={1995},
  publisher={ACM}
}

@article{parrilo2003semidefinite,
  title={Semidefinite programming relaxations for semialgebraic problems},
  author={Parrilo, Pablo A},
  journal={Mathematical Programming},
  volume={96},
  pages={293--320},
  year={2003},
  publisher={Springer}
}

@book{boyd1994linear,
  title={Linear matrix inequalities in system and control theory},
  author={Boyd, Stephen and El Ghaoui, Laurent and Feron, Eric and Balakrishnan, Venkataramanan},
  volume={15},
  year={1994},
  publisher={SIAM}
}

@book{skrzypczyk2023semidefinite,
  title={Semidefinite programming in quantum information science},
  author={Skrzypczyk, Paul and Cavalcanti, Daniel},
  year={2023},
  publisher={IOP Publishing}
}

@book{hiriart1996convex,
  author    = {Jean-Baptiste Hiriart-Urruty and Claude Lemar{\'e}chal},
  title     = {Convex Analysis and Minimization Algorithms~I: Fundamentals},
  series    = {Grundlehren der mathematischen Wissenschaften},
  volume    = {305},
  publisher = {Springer},
  address   = {Berlin, Heidelberg},
  year      = {1993},
  isbn      = {978-3-540-56850-6},
  doi       = {10.1007/978-3-662-02796-7}
}

@book{schneider1993convex,
  author    = {Rolf Schneider},
  title     = {Convex Bodies: The {Brunn}--{Minkowski} Theory},
  series    = {Encyclopedia of Mathematics and its Applications},
  volume    = {44},
  publisher = {Cambridge University Press},
  address   = {Cambridge},
  year      = {1993},
  isbn      = {0521352207},
  doi       = {10.1017/CBO9780511526282}
}

@article{alberti1994structure,
  title={On the structure of singular sets of convex functions},
  author={Alberti, Giovanni},
  journal={Calculus of Variations and Partial Differential Equations},
  volume={2},
  pages={17--27},
  year={1994},
  publisher={Springer}
}

@book{RoydenFitzpatrickRealAnalysis4e,
  author       = {H. L. Royden and P. M. Fitzpatrick},
  title        = {Real Analysis},
  edition      = {4},
  publisher    = {Pearson},
  year         = {2010},
  address      = {Upper Saddle River, NJ},
  isbn         = {978-0-13-143747-0},
}

@book{rockafellar1997convex,
  title={Convex analysis},
  author={Rockafellar, R Tyrrell},
  volume={28},
  year={1997},
  publisher={Princeton university press}
}

@article{audet2006mesh,
  title={Mesh adaptive direct search algorithms for constrained optimization},
  author={Audet, Charles and Dennis Jr, John E},
  journal={SIAM Journal on optimization},
  volume={17},
  number={1},
  pages={188--217},
  year={2006},
  publisher={SIAM}
}

@article{gratton2019direct,
  title={Direct search based on probabilistic feasible descent for bound and linearly constrained problems},
  author={Gratton, Serge and Royer, Cl{\'e}ment W and Vicente, Lu{\'\i}s Nunes and Zhang, Zaikun},
  journal={Computational Optimization and Applications},
  volume={72},
  number={3},
  pages={525--559},
  year={2019},
  publisher={Springer}
}

@article{liuzzi2010sequential,
  title={Sequential penalty derivative-free methods for nonlinear constrained optimization},
  author={Liuzzi, Giampaolo and Lucidi, Stefano and Sciandrone, Marco},
  journal={SIAM Journal on Optimization},
  volume={20},
  number={5},
  pages={2614--2635},
  year={2010},
  publisher={SIAM}
}

@article{Lucidi2002,
  title={Objective-derivative-free methods for constrained optimization},
  author={Lucidi, Stefano and Sciandrone, Marco and Tseng, Paul},
  journal={Mathematical Programming},
  volume={92},
  number={1},
  pages={37--59},
  year={2002},
  publisher={Springer}
}

@article{tao1998dc,
  title={A DC optimization algorithm for solving the trust-region subproblem},
  author={Tao, Pham Dinh and An, Le Thi Hoai},
  journal={SIAM Journal on Optimization},
  volume={8},
  number={2},
  pages={476--505},
  year={1998},
  publisher={SIAM}
}

@book{shiryaev2016probability1,
  author    = {Shiryaev, Albert N.},
  title     = {Probability-1},
  series    = {Graduate Texts in Mathematics},
  volume    = {95},
  edition   = {3rd},
  publisher = {Springer},
  year      = {2016},
  doi       = {10.1007/978-0-387-72206-1}
}

@book{shiryaev2019probability2,
  author    = {Albert N. Shiryaev},
  title     = {Probability-2},
  publisher = {Springer},
  series    = {Graduate Texts in Mathematics},
  volume    = {95},
  year      = {2019},
  edition   = {3rd},
  doi       = {10.1007/978-0-387-72208-5}
}

@book{EvansGariepy2015,
  author    = {Lawrence C. Evans and Ronald F. Gariepy},
  title     = {Measure Theory and Fine Properties of Functions, Revised Edition},
  edition   = {1},
  publisher = {Chapman and Hall/CRC},
  year      = {2015},
  doi       = {10.1201/b18333},
}

@article{martinez1994local,
  author = {Mart{\'i}nez, Jos{\'e} Mario},
  title = {Local Minimizers of Quadratic Functions on Euclidean Balls and Spheres},
  journal = {SIAM Journal on Optimization},
  volume = {4},
  number = {1},
  pages = {159-176},
  year = {1994},
  doi = {10.1137/0804009},
  URL = {https://doi.org/10.1137/0804009},
  eprint = {https://doi.org/10.1137/0804009}
}

\appendix
\section{Proofs of Section \ref{sec:preliminaries}}\label{sec:appendix_proofs_of_preliminaries}
\begin{proof}[Proof of {\cref{lem:feasible_directions_are_convex_cones}}]
    We note that since $0 \in D_\mathcal{C}(\bar{x})$, $D_\mathcal{C}(\bar{x}) \ne \emptyset$ for every nonempty set $\mathcal{C}$ and element $\bar{x} \in \mathcal{C}$.

    First, we prove that $D_{\mathcal{C}}\left(\bar{x}\right)$ is a cone.
    Let $y \in D_{\mathcal{C}}\left(\bar{x}\right)$. By definition, $\exists \lambda > 0$ s.t. $\bar{x} + \lambda y \in \mathcal{C}$. For $\omega = 0$, $\omega y = 0 \in D_\mathcal{C}(\bar{x})$. For $\omega > 0$, we choose $\lambda' = \dfrac{\lambda}{\omega}$, and derive
    \begin{equation*}
        \bar{x} + \lambda' \omega y = \bar{x} + \lambda y \in \mathcal{C}.
    \end{equation*}
    Therefore, $\omega y \in D_{\mathcal{C}}\left(\bar{x}\right)$ for all $\omega \ge 0$, and $D_{\mathcal{C}}(\bar{x})$ is a cone.

    Next, we prove that $D_{\mathcal{C}}(\bar{x})$ is convex.
    Let $y, z \in D_{\mathcal{C}}\left(\bar{x}\right)$. By definition, 
    \begin{equation*}
        \exists \eta, \omega > 0 \text{  such that  } \bar{x} + \eta y, \bar{x} + \omega z \in \mathcal{C}.
    \end{equation*}
    Assume without loss of generality that $\eta \le \omega$. By the convexity of $\mathcal{C}$,
    \begin{equation*}
        \bar{x} + \eta z \in \mathcal{C}.
    \end{equation*}
    Let $\lambda \in \left(0,1\right)$. By the convexity of $\mathcal{C}$,
    \begin{equation*}
        \lambda \left( \bar{x} + \eta y\right) + \left(1-\lambda\right)\left(\bar{x} + \eta z\right) = \bar{x} + \eta \left(\lambda y + \left(1-\lambda\right)z\right) \in \mathcal{C}.
    \end{equation*}
    Therefore, 
    \begin{equation*}
        \lambda y + \left(1-\lambda\right)z \in D_{\mathcal{C}}(\bar{x}).
    \end{equation*}
\end{proof}

\begin{proof}[Proof of \Cref{lem:foundational properties}]
    We first show that for  $S \subseteq \real^n$, 
    it holds that $\aff{S} = \aff{\cl{S}}$.

    If $S = \emptyset$, then $\cl{S} = \emptyset$, and the affine set of an empty set is $\aff{\emptyset} = \emptyset$. Hence, the claim is trivially true. 

    Otherwise, assume that $S \ne \emptyset$. 
    The inclusion $\aff{S} \subseteq \aff{\cl{S}}$ is trivial since $S \subseteq \cl{S}$. 

    In the opposite direction, let $y \in \aff{\cl{S}}$. Therefore, there exist $\bar{x}_1, \ldots,\bar{x}_k \in \cl{S}$, and scalars $a_1, \ldots,a_k$ such that $\sum\limits_{i=1}^k a_i = 1$ and
    \begin{equation*}
        y = \sum\limits_{i=1}^k a_i \bar{x}_i.
    \end{equation*}
    Since $\bar{x}_i \in \cl{S}$, for each $i$ there exists  a sequence $\left\{x_i^j\right\}_{j \ge 0} \subseteq S$ such that $x_i^j\xrightarrow{j \to \infty} \bar{x}_i$. Let $y^j = \sum\limits_{i=1}^k a_i x_i^j$. By its construction, $y^j \in \aff{S}$ for all $j$. Since $\aff{S}$ is a closed set, and $y^j \xrightarrow{j\to\infty} y$, it follows that $y \in \aff{S}$.
    
    \medskip

    We now move to prove that for $\bar{x}, \bar{y} \in S$ it holds that $\aff{S - \bar{x}} = \aff{S - \bar{y}}.$
    Let $z \in \aff{S - \bar{x}}$. 
    Our goal is to show that $z\in \aff{S - \bar{y}}$, that is, to construct a set of vectors $u_1,\ldots,u_l \in S - \bar{y}$ and coefficients $b_1, \ldots, b_l \in \real$ such that $\sum\limits_{i=1}^l b_i =1$ and $z = \sum\limits_{i=1}^l b_i u_i$.
    
    Since $z \in \aff{S - \bar{x}}$, there exist $v_1,\ldots, v_k \in S - \bar{x}$ and $a_1, \ldots, a_k \in \real$, such that $\sum\limits_{i=1}^k a_i = 1$, and
    \begin{equation*}
        z = \sum\limits_{i=1}^k a_i v_i. 
    \end{equation*}
    By definition of $S - \bar{x}$, the fact that $v_1, \ldots, v_k \in S - \bar{x}$ implies $v_1 + \bar{x},\ldots, v_k + \bar{x} \in S$. 
    Hence, by definition of $S - \bar{y}$, we have that $0, \bar{x} - \bar{y} \in S - \bar{y}$, and $v_1 + \bar{x} - \bar{y},\ldots, v_k + \bar{x} - \bar{y} \in S - \bar{y}$. 
    Therefore,
    \begin{equation*}
        z = \sum\limits_{i=1}^k a_i v_i = 1 \cdot 0 - 1\cdot \left(\bar{x} - \bar{y}\right) + \sum\limits_{i=1}^k a_i \left(v_i + \bar{x} - \bar{y}\right).
    \end{equation*}
    Since the coefficients $1, -1, a_1, \ldots, a_k$ sum to $1$, it follows that $z \in \aff{S - \bar{y}}$, and consequently $\aff{S - \bar{x}} \subseteq \aff{S - \bar{y}}$. 
    The reverse inclusion $\aff{S-\bar{y}}\subseteq\aff{S-\bar{x}}$ is proved in the same way. 

    \medskip

    Now assume $S$ is convex. Nonemptiness of $\ri{S}$ and the line segment principle are \citep[Proposition 2.40]{RW98}. We show $\ri{S} = \ri{\cl{S}}$.
    Let $y \in \ri{S}$. 
    Then there exists $\epsilon > 0$ such that 
    \begin{equation*}
        \mathcal{B}(y, \epsilon) \cap \aff{S} \subseteq S.
    \end{equation*}
    Since $\aff{S} = \aff{\cl{S}}$ and $S \subseteq \cl{S}$, it follows that
    \begin{equation*}
        \mathcal{B}(y,\epsilon) \cap \aff{\cl{S}} \subseteq \cl{S}. 
    \end{equation*}
    Since $y \in S \subseteq \cl{S}$, it follows that $y \in \ri{\cl{S}}$.

    Let $x \in \ri{\cl{S}}$. By definition, there exists $\epsilon > 0$ such that 
    \begin{equation*}
        \mathcal{B}(x,\epsilon) \cap \aff{\cl{S}} = \mathcal{B}(x,\epsilon) \cap \aff{S} \subseteq \cl{S}.
    \end{equation*}
    Let $y \in \ri{S}$. If $x = y$, then $x \in \ri{S}$ and we are done. Otherwise, $\left\Vert x - y \right\Vert > 0$. Let $\delta = \dfrac{\epsilon}{2 \left\Vert x - y \right\Vert}$, and define
    \begin{equation*}
        z = \left(1 + \delta\right)x - \delta y.
    \end{equation*}
    Since $x, y \in \cl{S}$, it clearly follows that $z \in \aff{\cl{S}} = \aff{S}$. 
    Additionally, $\left\Vert z - x \right\Vert = \left\Vert \left(1 + \delta\right)x - \delta y - x\right\Vert = \delta \left\Vert y - x \right\Vert \le 0.5\epsilon$, and therefore, $z \in \mathcal{B}(x,\epsilon)$. 
    
    Now, $z \in \mathcal{B}(x,\epsilon) \cap \aff{\cl{S}}$, and since $x \in \ri{\cl{S}}$, it follows that $\mathcal{B}(x,\epsilon) \cap \aff{\cl{S}} \subseteq \cl{S}$. Combining both results, we deduce that $z \in \cl{S}$. Together with the equation
    \begin{equation*}
        x = \dfrac{1}{1 + \delta}z + \dfrac{\delta}{1+\delta}y,
    \end{equation*}
    and the fact that $y \in \ri{S}$, this in turn implies $x \in \ri{S}$ by the line segment principle.

\end{proof}

\section{Proofs of Section \ref{sec:analysis}}\label{sec:appendix_proofs_analysis}
\begin{proof}[Proof of \cref{lem:optimality_measure_zero_implies_no_descent_in_span}]
Write $h'(\bar{x};d)=\inner{\nabla f(\bar{x})}{d}-g'(\bar{x};d)$.
The assumption $\mathcal{O}^\gamma(\bar{x};V,r)=0$ means $\mu^\gamma(\bar{x};v,r)=0$ for every $v\in V$, so
\begin{equation*}
    h(\bar{x})\;\le\;h(\bar{x}+tv)+\tfrac{\gamma}{2}t^2
    \qquad\text{for all }t\in\bigl[0,\min\{\tau(\bar{x};v),r\}\bigr].
\end{equation*}
Since $V\subseteq D_{\mathcal{C}}(\bar{x})$, each $v\in V$ satisfies $\tau(\bar{x};v)>0$ by \cref{def:feasible_directions}. Combined with $r>0$, the feasible interval has positive length. For every $t\in\bigl(0,\min\{\tau(\bar{x};v),r\}\bigr]$,
\begin{equation*}
    \frac{h(\bar{x}+tv)-h(\bar{x})}{t}+\tfrac{\gamma}{2}t\;\ge\;0.
\end{equation*}
Sending $t\to 0^+$ yields $h'(\bar{x};v)\ge 0$; the directional derivative exists and is finite because $f$ is continuously differentiable and $g:\real^n\to\real$ is convex \citep[Theorem~23.1]{rockafellar1997convex}.

Now let $d\in\mathcal{K}$. Then $d=\sum_{v\in V}\lambda_v v$ for some $\lambda_v\ge 0$. The map $g'(\bar{x};\cdot)$ is sublinear, hence
\begin{align*}
    h'(\bar{x};d)
    &=\sum_{v}\lambda_v\inner{\nabla f(\bar{x})}{v}-g'\Bigl(\bar{x};\sum_{v}\lambda_v v\Bigr) \\
    &\ge
    \sum_{v}\lambda_v\bigl(\inner{\nabla f(\bar{x})}{v}-g'(\bar{x};v)\bigr)
    =\sum_{v}\lambda_v h'(\bar{x};v)
    \ge 0.
\end{align*}
\end{proof}

\begin{proof}[Proof of \cref{lem:h-prime-lipschitz-span}]
Since $\mathcal{C}$ is compact and $\nabla f$ is continuous, $\sup_{x\in\mathcal{C}}\norm{\nabla f(x)}<\infty$.
The set $\partial g$ is locally bounded on the interior of $\mathrm{dom}\,g$, which contains a neighborhood of $\mathcal{C}$, hence $\partial g$ is uniformly bounded on $\mathcal{C}$ and $\Lambda<\infty$.

Write $h'(\bar{x};d)=\inner{\nabla f(\bar{x})}{d}-g'(\bar{x};d)$.
The map $g'(\bar{x};\cdot)$ is the support function of $\partial g(\bar{x})$, so
\begin{equation*}
    \bigl\lvert g'(\bar{x};d)-g'(\bar{x};d')\bigr\rvert
    \;=\;
    \bigl\lvert \sigma_{\partial g(\bar{x})}(d)-\sigma_{\partial g(\bar{x})}(d')\bigr\rvert
    \;\le\;
    \max_{s\in\partial g(\bar{x})}\norm{s}\,\norm{d-d'}.
\end{equation*}
The linear map $d\mapsto\inner{\nabla f(\bar{x})}{d}$ is Lipschitz with modulus $\norm{\nabla f(\bar{x})}$.
Adding the two moduli,
\begin{equation*}
    \bigl\lvert h'(\bar{x};d)-h'(\bar{x};d')\bigr\rvert
    \;\le\;
    \Bigl(\norm{\nabla f(\bar{x})}+\max_{s\in\partial g(\bar{x})}\norm{s}\Bigr)\norm{d-d'}
    \;\le\;
    \Lambda\norm{d-d'}.
\end{equation*}
Finally $h'(\bar{x};0)=0$, so $\bigl\lvert h'(\bar{x};d)\bigr\rvert=\bigl\lvert h'(\bar{x};d)-h'(\bar{x};0)\bigr\rvert\le\Lambda\norm{d}$. This concludes part~1.

For part~2, $\inner{\nabla f(\bar{x})}{\cdot}$ is linear and $g'(\bar{x};\cdot)$ is sublinear, so
\begin{equation*}
    h'(\bar{x};d)
    =\sum_v\lambda_v\inner{\nabla f(\bar{x})}{v}-g'\Bigl(\bar{x};\sum_v\lambda_v v\Bigr)
    \;\ge\;
    \sum_v\lambda_v\bigl(\inner{\nabla f(\bar{x})}{v}-g'(\bar{x};v)\bigr)
    =\sum_v\lambda_v h'(\bar{x};v).
\end{equation*}
Negate and use $-h'\le[-h']_+$.
\end{proof}

\begin{proof}[Proof of \cref{prop:PDM_lower_bound}]
The bound is per-direction, so it is enough to treat a single $v\in V$ and sum the squares.
Write $\phi(t):=f(\bar{x}+tv)+\tfrac{\gamma}{2}t^2$ for $t\ge 0$, so $\mu^\gamma(\bar{x};v,r)$ is the smallest minimizer of $\phi$ on $[0,\tau_r(\bar{x};v)]$, and $\phi'(t)=\inner{\nabla f(\bar{x}+tv)}{v}+\gamma t$.

If $\inner{\nabla f(\bar{x})}{v}\ge 0$, this $v$ contributes $0$ to the right-hand side of \eqref{eq:PDM_lower_bound}, and $\mu^\gamma(\bar{x};v,r)\ge 0$ is immediate.

Now assume $\inner{\nabla f(\bar{x})}{v}<0$, and set $t^\star:=-\inner{\nabla f(\bar{x})}{v}/(L+\gamma)>0$.
For $0\le t<t^\star$,
\begin{align*}
    \phi'(t)
    \;&\le\;
    \lvert\phi'(t)-\phi'(0)\rvert+\phi'(0) \\
    \;&\le\;
    \lvert\inner{\nabla f(\bar{x}+tv)}{v}-\inner{\nabla f(\bar{x})}{v}\rvert+\gamma t+\phi'(0) \\
    \;&\le\;
    (L+\gamma)t+\phi'(0)
    \;<\;0,
\end{align*}
where the third inequality uses $f\in C^{1,1}_L$ and $\norm{v}=1$, so $\lvert\inner{\nabla f(\bar{x}+tv)-\nabla f(\bar{x})}{v}\rvert\le Lt$.
Thus $\phi$ is strictly decreasing on $[0,t^\star]$, so every minimizer of $\phi$ on $[0,\tau_r(\bar{x};v)]$ is at least $\min\{t^\star,\tau_r(\bar{x};v)\}$. Therefore
\begin{equation*}
    \mu^\gamma(\bar{x};v,r)
    \;\ge\;
    \min\bigl\{t^\star,\,\tau_r(\bar{x};v)\bigr\}
    \;=\;
    \frac{1}{L+\gamma}\min\bigl\{-\inner{\nabla f(\bar{x})}{v},\,(L+\gamma)\tau_r(\bar{x};v)\bigr\}.
\end{equation*}
Combining with the nonnegative inner-product case,
\begin{equation*}
    \mu^\gamma(\bar{x};v,r)
    \;\ge\;
    \frac{1}{L+\gamma}\min\bigl\{[-\inner{\nabla f(\bar{x})}{v}]_+,\,(L+\gamma)\tau_r(\bar{x};v)\bigr\}
\end{equation*}
for every $v\in V$. Squaring and summing over $V$ gives \eqref{eq:PDM_lower_bound}.
\end{proof}

\begin{proof}[Proof of \cref{lem:per-direction-decrease}]
If $\theta=0$ there is nothing to prove; assume $\theta>0$, so $h'(\bar{x};v)=-\theta$.
For $t\in[0,\tau_r(\bar{x};v)]$ the segment $[\bar{x},\bar{x}+tv]\subseteq\mathcal{C}$.
Since $f\in C^{1,1}_L$ and $\norm{v}=1$,
\begin{equation*}
    f(\bar{x}+tv)\;\le\;f(\bar{x})+t\inner{\nabla f(\bar{x})}{v}+\tfrac{L}{2}t^2.
\end{equation*}
Convexity of $g$ gives $g(\bar{x}+tv)\ge g(\bar{x})+t\,g'(\bar{x};v)$ for $t\ge 0$.
Subtracting and adding the regularizer,
\begin{align*}
    h(\bar{x}+tv)+\tfrac{\gamma}{2}t^2-h(\bar{x})
    \;&\le\;
    t\bigl(\inner{\nabla f(\bar{x})}{v}-g'(\bar{x};v)\bigr)+\tfrac{L+\gamma}{2}t^2 \\
    \;&=\;
    t\,h'(\bar{x};v)+\tfrac{L+\gamma}{2}t^2
    \;=\;
    -\theta t+\tfrac{L+\gamma}{2}t^2
    \;=:\;\psi(t).
\end{align*}
Hence
\begin{equation*}
    h(\bar{x})-\min_{0\le t\le \tau_r(\bar{x};v)}\Bigl\{h(\bar{x}+tv)+\tfrac{\gamma}{2}t^2\Bigr\}
    \;\ge\;
    -\min_{0\le t\le \tau_r(\bar{x};v)}\psi(t).
\end{equation*}
The derivative $\psi'(t)=-\theta+(L+\gamma)t$ vanishes at $t^\star=\theta/(L+\gamma)$, and
\begin{equation*}
    \psi(t^\star)
    \;=\;
    -\frac{\theta^2}{L+\gamma}+\frac{L+\gamma}{2}\cdot\frac{\theta^2}{(L+\gamma)^2}
    \;=\;
    -\frac{\theta^2}{2(L+\gamma)}.
\end{equation*}
If $t^\star\le \tau_r(\bar{x};v)$ then $\min_{[0,\tau_r(\bar{x};v)]}\psi=\psi(t^\star)$, and the claim follows.
If $t^\star>\tau_r(\bar{x};v)$ then $\psi$ is decreasing on $[0,\tau_r(\bar{x};v)]$, so the minimum is $\psi(\tau_r(\bar{x};v))$.
Since $\tau_r(\bar{x};v)<t^\star$ we have $\tau_r(\bar{x};v)^2<\tau_r(\bar{x};v)\,t^\star$, and therefore
\begin{align*}
    \psi(\tau_r(\bar{x};v))
    &=
    -\theta \tau_r(\bar{x};v)+\tfrac{L+\gamma}{2}\tau_r(\bar{x};v)^2 \\
    &<
    -\theta \tau_r(\bar{x};v)+\tfrac{L+\gamma}{2}\tau_r(\bar{x};v)\, t^\star \\
    &=
    -\theta \tau_r(\bar{x};v)+\tfrac12\theta \tau_r(\bar{x};v) \\
    &=
    -\tfrac12\theta \tau_r(\bar{x};v).
\end{align*}
Thus the left-hand side is greater than $\tfrac12\theta \tau_r(\bar{x};v)$.
\end{proof}

\begin{proof}[Proof of \cref{lem:rolling_balls}]
Write $r_{\mathrm{in}}=1/\kappa_{\max}$, $r_{\mathrm{out}}=1/\kappa_{\min}$.
Let $\hat{x}\in\rb{\mathcal{C}}$ with inward unit normal $u$, and let $\nu=-u$ be the outer unit normal at $\hat{x}$, unique since $\mathcal{C}$ has a $C^2$
boundary. Recall that the support function of a nonempty $\mathcal{X}\subseteq\mathbb{R}^n$ is
$h_{\mathcal{X}}(w)=\sup_{x\in\mathcal{X}}\inner{x}{w}$.
Let $\mathcal{B}_{\mathrm{in}}:=\mathcal{B}(\hat{x}+r_{\mathrm{in}}u,\,r_{\mathrm{in}})$.
Since all principal curvatures of $\mathcal{B}_{\mathrm{in}}$ equal $\kappa_{\max}$, they upper
bound those of $\mathcal{C}$ at every pair of boundary points sharing an outer unit normal. By
\citep[Theorem~2.5.4]{schneider1993convex}, $h_{\mathcal{C}}-h_{\mathcal{B}_{\mathrm{in}}}$ is a support function, so by
\citep[Theorem~1.7.1]{schneider1993convex} there is a unique convex body $M$ with $h_M=h_{\mathcal{C}}-h_{\mathcal{B}_{\mathrm{in}}}$.
By \citep[Theorem~1.7.5]{schneider1993convex} and uniqueness, $M+\mathcal{B}_{\mathrm{in}}=\mathcal{C}$.
We claim $0\in M$; then
$\mathcal{B}_{\mathrm{in}}=\{0\}+\mathcal{B}_{\mathrm{in}}\subseteq M+\mathcal{B}_{\mathrm{in}}=\mathcal{C}$.
Indeed there are $\bar{m}\in M$ and
$\bar{b}\in\mathcal{B}_{\mathrm{in}}$ with $\bar{m}+\bar{b}=\hat{x}$, and
\begin{equation*}
h_{\mathcal{C}}(\nu)=\inner{\hat{x}}{\nu}=\inner{\bar{b}}{\nu}+\inner{\bar{m}}{\nu}\le h_{\mathcal{B}_{\mathrm{in}}}(\nu)+h_M(\nu)=h_{\mathcal{C}}(\nu),
\end{equation*}
so equality holds throughout. In particular $\inner{\bar{b}}{\nu}=h_{\mathcal{B}_{\mathrm{in}}}(\nu)$, so $\bar{b}$ maximizes $\inner{b}{\nu}$ over $b\in\mathcal{B}_{\mathrm{in}}$.
This maximizer is unique and equal to $\hat{x}$: the center of $\mathcal{B}_{\mathrm{in}}$ is $\hat{x}+r_{\mathrm{in}}u=\hat{x}-r_{\mathrm{in}}\nu$, so the unique farthest point of the ball in the direction $\nu$ is $(\hat{x}-r_{\mathrm{in}}\nu)+r_{\mathrm{in}}\nu=\hat{x}$.
Thus $\bar{b}=\hat{x}$, and therefore $\bar{m}=0$.

Now set $\mathcal{B}_{\mathrm{out}}:=\mathcal{B}(\hat{x}+r_{\mathrm{out}}u,\,r_{\mathrm{out}})$.
The principal curvatures of $\mathcal{B}_{\mathrm{out}}$ are constantly $\kappa_{\min}$, hence they lower
bound those of $\mathcal{C}$ at every pair of boundary points sharing an outer unit normal.
By the same comparison \citep[Theorem~2.5.4]{schneider1993convex}, $h_{\mathcal{B}_{\mathrm{out}}}-h_{\mathcal{C}}$ is a support function, so by
\citep[Theorem~1.7.1]{schneider1993convex} there is a unique convex body $M'$ with $h_{M'}=h_{\mathcal{B}_{\mathrm{out}}}-h_{\mathcal{C}}$, and as above $M'+\mathcal{C}=\mathcal{B}_{\mathrm{out}}$.
We claim $0\in M'$; then $\mathcal{C}=\{0\}+\mathcal{C}\subseteq M'+\mathcal{C}=\mathcal{B}_{\mathrm{out}}$.
Since $\hat{x}\in\mathcal{B}_{\mathrm{out}}$, there are $\bar{m}'\in M'$ and $\bar{c}\in\mathcal{C}$ with $\bar{m}'+\bar{c}=\hat{x}$, and
\begin{equation*}
h_{\mathcal{B}_{\mathrm{out}}}(\nu)=\inner{\hat{x}}{\nu}=\inner{\bar{c}}{\nu}+\inner{\bar{m}'}{\nu}\le h_{\mathcal{C}}(\nu)+h_{M'}(\nu)=h_{\mathcal{B}_{\mathrm{out}}}(\nu),
\end{equation*}
so equality holds throughout. In particular $\inner{\bar{c}}{\nu}=h_{\mathcal{C}}(\nu)$.
But $\hat{x}$ itself attains $h_{\mathcal{C}}(\nu)$, because $\nu$ is an outer unit normal to $\mathcal{C}$ at $\hat{x}$.
This maximizer is unique: $\kappa_{\min}>0$ forces $\mathcal{C}$ to be strictly convex, so a supporting hyperplane meets $\mathcal{C}$ at exactly one point.
Thus $\bar{c}=\hat{x}$ and $\bar{m}'=0$. We conclude
\begin{equation*}
\mathcal{B}_{\mathrm{in}}\;\subseteq\;\mathcal{C}\;\subseteq\;\mathcal{B}_{\mathrm{out}},
\end{equation*}
which is part~1. The outer inclusion implies that $\mathcal{C}$ is bounded, hence compact.

For part~2, if $\rho(x)=0$ then $x\in\rb{\mathcal{C}}$ and $\Pi(x)=\{x\}$.
Assume $0<\rho(x)<r_{\mathrm{in}}$.
The set $\rb{\mathcal{C}}$ is compact, so $\Pi(x)$ is nonempty.
Let $y\in\Pi(x)$ and let $u$ be the inward unit normal at $y$.
Write $\rho:=\rho(x)$, so $x=y+\rho u$, and set $z:=y+r_{\mathrm{in}}u$.
Part~1 applied at $y$ gives $\mathcal{B}(z,r_{\mathrm{in}})\subseteq\mathcal{C}$.
The ball of radius $\rho$ about $x$ is contained in the ball of radius $r_{\mathrm{in}}$ about $z$: for every $w$ with $\norm{w-x}\le\rho$,
\begin{equation*}
\norm{w-z}\;\le\;\norm{w-x}+\norm{x-z}\;\le\;\rho+(r_{\mathrm{in}}-\rho)\;=\;r_{\mathrm{in}}.
\end{equation*}
Equality holds throughout if and only if $\norm{w-x}=\rho$ and $w$, $x$, $z$ are collinear with $x$ between $z$ and $w$, which forces $w=x-\rho u=y$.
Thus the two balls share exactly one boundary point, namely $y$.
Since $\mathcal{B}(x,\rho)\subseteq\mathcal{B}(z,r_{\mathrm{in}})\subseteq\mathcal{C}$, every other point of $\mathcal{B}(x,\rho)$ lies in the open ball of radius $r_{\mathrm{in}}$ about $z$, hence in $\ri{\mathcal{C}}$.
So $\mathcal{B}(x,\rho)$ meets $\rb{\mathcal{C}}$ only at $y$.
Every point of $\Pi(x)$ lies on $\rb{\mathcal{C}}$ and on the sphere of radius $\rho$ about $x$, so $\Pi(x)=\{y\}$.
\end{proof}

\begin{proof}[Proof of \cref{lem:curvature_feasibility_bounds}]
Write $r_{\mathrm{in}}=1/\kappa_{\max}$ and $r_{\mathrm{out}}=1/\kappa_{\min}$,
and let $\mathcal{B}_{\mathrm{in}}:=\mathcal{B}(\hat{x}+r_{\mathrm{in}}u,\,r_{\mathrm{in}})$ and
$\mathcal{B}_{\mathrm{out}}:=\mathcal{B}(\hat{x}+r_{\mathrm{out}}u,\,r_{\mathrm{out}})$.
By \cref{lem:rolling_balls}, $\mathcal{B}_{\mathrm{in}}\subseteq\mathcal{C}\subseteq\mathcal{B}_{\mathrm{out}}$.

For part~1, let $z:=\hat{x}+r_{\mathrm{in}}u$, so $\bar{x}-z=-(r_{\mathrm{in}}-\rho)u$.
Then $\bar{x}+td\in\mathcal{B}_{\mathrm{in}}$ iff $\norm{td-(r_{\mathrm{in}}-\rho)u}^2\le r_{\mathrm{in}}^2$,
i.e.
\begin{equation*}
t^2-2t(r_{\mathrm{in}}-\rho)c\;\le\;r_{\mathrm{in}}^2-(r_{\mathrm{in}}-\rho)^2=2r_{\mathrm{in}}\rho-\rho^2,
\end{equation*}
whose largest admissible $t$ is the largest root of the quadratic
$t^2-2t(r_{\mathrm{in}}-\rho)c-(2r_{\mathrm{in}}\rho-\rho^2)$, namely
$(r_{\mathrm{in}}-\rho)c+\sqrt{(r_{\mathrm{in}}-\rho)^2c^2+2r_{\mathrm{in}}\rho-\rho^2}$; since
$\mathcal{B}_{\mathrm{in}}\subseteq\mathcal{C}$ this lower bounds $\tau(\bar{x};d)$.

For (1a), first take $0\le c\le\tfrac12$. The inequality
\begin{equation*}
(r_{\mathrm{in}}-\rho)c+\sqrt{(r_{\mathrm{in}}-\rho)^2c^2+2r_{\mathrm{in}}\rho-\rho^2}\;\ge\;2r_{\mathrm{in}}c
\end{equation*}
is equivalent, after isolating the square root and squaring (both sides being nonnegative),
to $(r_{\mathrm{in}}-\rho)^2c^2+2r_{\mathrm{in}}\rho-\rho^2\ge(r_{\mathrm{in}}+\rho)^2c^2$,
i.e.\ to $2r_{\mathrm{in}}\rho-\rho^2\ge 4r_{\mathrm{in}}\rho c^2$, i.e.\ to
$c^2\le(2r_{\mathrm{in}}-\rho)/(4r_{\mathrm{in}})$, which holds since $\rho\le r_{\mathrm{in}}$
gives $(2r_{\mathrm{in}}-\rho)/(4r_{\mathrm{in}})\ge\tfrac14$. The stated expression is nondecreasing in $c$,
so for $c>\tfrac12$ it is at least its value at $c=\tfrac12$, hence $\tau(\bar{x};d)\ge r_{\mathrm{in}}$.
Thus (1a) holds in all cases.

For (1b), drop the nonnegative $c$-terms and use $2r_{\mathrm{in}}\rho-\rho^2\ge r_{\mathrm{in}}\rho$ for $\rho\le r_{\mathrm{in}}$.
For (1c), substituting $c=1$ gives
$(r_{\mathrm{in}}-\rho)+\sqrt{(r_{\mathrm{in}}-\rho)^2+2r_{\mathrm{in}}\rho-\rho^2}=(r_{\mathrm{in}}-\rho)+r_{\mathrm{in}}=2r_{\mathrm{in}}-\rho\ge r_{\mathrm{in}}$,
since $(r_{\mathrm{in}}-\rho)^2+2r_{\mathrm{in}}\rho-\rho^2=r_{\mathrm{in}}^2$.

For part~2, the same computation with $\mathcal{C}\subseteq\mathcal{B}_{\mathrm{out}}$
gives $\tau(\bar{x};d)\le(r_{\mathrm{out}}-\rho)c+\sqrt{(r_{\mathrm{out}}-\rho)^2c^2+X}$ with
$X:=2r_{\mathrm{out}}\rho-\rho^2\ge 0$. Put $a:=(r_{\mathrm{out}}-\rho)|c|\ge 0$ and use
$\sqrt{a^2+X}\le a+\sqrt{X}$. For $c\ge 0$ this gives
$\tau(\bar{x};d)\le 2(r_{\mathrm{out}}-\rho)c+\sqrt{X}\le 2r_{\mathrm{out}}c+\sqrt{2r_{\mathrm{out}}\rho}$;
for $c<0$ the leading term is $-a$, so $\tau(\bar{x};d)\le\sqrt{X}\le\sqrt{2r_{\mathrm{out}}\rho}$.
Combining the two cases, $\tau(\bar{x};d)\le 2r_{\mathrm{out}}[c]_+ + \sqrt{2r_{\mathrm{out}}\rho}$.

For part~3, since $\tau(\bar{x};d)=\sup\{t\ge 0:\bar{x}+td\in\mathcal{C}\}$, any condition
necessary for $\bar{x}+td\in\mathcal{C}$ gives an upper bound on $\tau(\bar{x};d)$.
As $u$ is the inward unit normal at $\hat{x}\in\rb{\mathcal{C}}$, the hyperplane through
$\hat{x}$ orthogonal to $u$ supports $\mathcal{C}$, so
$\mathcal{C}\subseteq H:=\{y:\inner{y-\hat{x}}{u}\ge 0\}$. Since $\bar{x}=\hat{x}+\rho u$ with
$\norm{u}=1$ we have $\inner{\bar{x}-\hat{x}}{u}=\rho$, hence
\begin{equation*}
\inner{(\bar{x}+td)-\hat{x}}{u}=\rho+t\inner{d}{u}=\rho+tc,
\end{equation*}
and $\bar{x}+td\in\mathcal{C}\subseteq H$ forces $\rho+tc\ge 0$. For $c<0$ this is equivalent to
$t\le\rho/|c|$, so $\tau(\bar{x};d)\le\rho/|c|$.

Taking $d=-u$ (so $c=-1$) gives $\tau(\bar{x};-u)\le\rho$. The reverse inequality follows from
convexity: $\hat{x},\bar{x}\in\mathcal{C}$, so the segment
$[\bar{x},\hat{x}]=\{\bar{x}+t(-u):0\le t\le\rho\}$ lies in $\mathcal{C}$. Hence
$\tau(\bar{x};-u)=\rho$.
\end{proof}

\begin{proof}[Proof of \cref{cor:boundary_eta_tau}]
This is \cref{lem:curvature_feasibility_bounds} at $\rho=0$. By \cref{lem:eta_FD_properties}(b) the $\eta$-FD
condition on a unit $d$ reads $c=\inner{d}{u}\ge\eta$, and part~1 at $\rho=0$ gives
$\tau(\hat{x};d)\ge r_{\mathrm{in}}c+\sqrt{r_{\mathrm{in}}^2 c^2}=2r_{\mathrm{in}}c\ge 2\eta/\kappa_{\max}$. If $d$ is not an $\eta$-FD then
$c<\eta$; if $c<0$ then $\tau(\hat{x};d)=0$ by part~3, and if $0\le c<\eta$ then part~2 at
$\rho=0$ gives $\tau(\hat{x};d)\le 2r_{\mathrm{out}}c\le 2\eta/\kappa_{\min}$.
\end{proof}

\begin{proof}[Proof of \cref{prop:quantitative-spanning}]
Set $\bar{\eta}:=\eta/(2\sqrt{\dim\mathcal{C}})$.
Fix $x\in\mathcal{C}$ and a unit vector $d\in D_{\mathcal{C}}(x)$ with $\tau(x;d)\ge R$.
Let $V(x)$ denote the non-exploration output of \cref{alg:construct_directions} at $\alpha=0$ and $\rho_0=r_{\mathrm{in}}/2$.
Under the assumptions of \cref{lem:rolling_balls} the boundary is $C^2$, so the second branch of \cref{alg:construct_directions} is never taken.
If $\rho>\rho_0$ then $T=\emptyset$ and $|U|=2\dim\mathcal{C}$, and Step~3 with $\alpha=0$ returns $U$, hence $|V(x)|=2\dim\mathcal{C}$.
If $\rho\le\rho_0$ then $|T|=2(\dim\mathcal{C}-1)$ and $|U|\le 2\dim\mathcal{C}$; Step~3 with $\alpha=0$ returns one direction for each element of $T$ and of $U$, hence $|V(x)|=|T|+|U|\le 4\dim\mathcal{C}-2$.
In all cases $|V(x)|\le 4\dim\mathcal{C}$.

\textbf{Deep interior} ($\rho>\rho_0=r_{\mathrm{in}}/2$).
Write $d=\sum_i a_i w_i$ and take $\lambda_{\mathrm{sgn}(a_i)w_i}=|a_i|$.
Then $\sum\lambda=\norm{a}_1\le\sqrt{\dim\mathcal{C}}$ by Cauchy--Schwarz, hence $A=\sqrt{\dim\mathcal{C}}$, and $\tau(x;v)\ge\rho>r_{\mathrm{in}}/2$ for every unit $v$, so $\tau_r(x;v)\ge\min\{r_{\mathrm{in}}/2,\,r\}$ and
\begin{equation*}
    \sum_v\frac{\lambda_v}{\tau_r(x;v)}\;\le\;\sqrt{\dim\mathcal{C}}\Bigl(\frac1r+2\kappa_{\max}\Bigr).
\end{equation*}
This bound is independent of $\eta$.

\textbf{Boundary or shallow interior} ($\rho\le\rho_0$).
Decompose $d=c\,u+\sum_{i=1}^{\dim\mathcal{C}-1}a_i w_i$, $c=\inner{d}{u}$, $c^2+\sum a_i^2=1$, hence $|c|+\sum_i|a_i|\le\sqrt{\dim\mathcal{C}}$ by Cauchy--Schwarz.

\textit{Inward tilt: $c\ge\eta$, use the tilted set.}
Let $\bar{v}_i$ be the normalization of $\bar{\eta}u+(1-\bar{\eta})\,\mathrm{sgn}(a_i)w_i$, and write $\varrho:=\norm{\bar{\eta}u+(1-\bar{\eta})\,\mathrm{sgn}(a_i)w_i}\in[\tfrac1{\sqrt2},1]$, so $\varrho\bar{v}_i=\bar{\eta}u+(1-\bar{\eta})\,\mathrm{sgn}(a_i)w_i$. Set
\begin{equation*}
    \tilde{c} := c-\frac{\bar{\eta}}{1-\bar{\eta}}\sum_j|a_j|,
    \qquad
    \tilde{a}_i := \frac{a_i}{1-\bar{\eta}}.
\end{equation*}
Then
\begin{align*}
    \sum_i|\tilde{a}_i|\varrho\bar{v}_i
    &=
    \sum_i\frac{|a_i|}{1-\bar{\eta}}\bigl(\bar{\eta}u+(1-\bar{\eta})\,\mathrm{sgn}(a_i)w_i\bigr)\\
    &=
    \frac{\bar{\eta}}{1-\bar{\eta}}\Bigl(\sum_i|a_i|\Bigr)u
    +\sum_i a_i w_i,
\end{align*}
and therefore $d=\tilde{c}\,u+\sum_i|\tilde{a}_i|\varrho\bar{v}_i$.
Since $\bar{\eta}\le\tfrac12$ and $\sum_j|a_j|\le\sqrt{\dim\mathcal{C}}$,
\begin{equation*}
    \frac{\bar{\eta}}{1-\bar{\eta}}\sum_j|a_j| \;\le\; 2\bar{\eta}\sqrt{\dim\mathcal{C}} \;=\;\eta\;\le\; c,
\end{equation*}
so $\tilde{c}\ge 0$.
Using $\lambda_u=\tilde{c}$ and $\lambda_{\bar{v}_i}=|\tilde{a}_i|\varrho$,
\begin{equation*}
    \sum\lambda \;=\; \tilde{c}+\sum_i|\tilde{a}_i|\varrho \;\le\; \tilde{c}+\sum_i|\tilde{a}_i| \;=\; c-\tfrac{\bar{\eta}}{1-\bar{\eta}}\textstyle\sum|a_j|+\tfrac{1}{1-\bar{\eta}}\textstyle\sum|a_i| \;=\; c+\sum_i|a_i|\;\le\;\sqrt{\dim\mathcal{C}},
\end{equation*}
since $\varrho\le 1$ and $\tfrac{1}{1-\bar{\eta}}-\tfrac{\bar{\eta}}{1-\bar{\eta}}=1$.
Hence $A=\sqrt{\dim\mathcal{C}}$, independently of $\eta$.
For the $B$-sum, write $c_i:=\inner{\bar{v}_i}{u}=\bar{\eta}/\varrho$.
Since $\varrho\le 1$ we have $c_i\ge\bar{\eta}\ge 0$, so \cref{lem:curvature_feasibility_bounds}(1a) gives $\tau(x;\bar{v}_i)\ge\min\{2c_i,1\}\,r_{\mathrm{in}}$.
Because $c_i\ge\bar{\eta}$ and $\bar{\eta}\le\tfrac12$, we have $\min\{2c_i,1\}\ge 2\bar{\eta}$, hence $\tau(x;\bar{v}_i)\ge 2\bar{\eta}\,r_{\mathrm{in}}=2\bar{\eta}/\kappa_{\max}=\eta/(\sqrt{\dim\mathcal{C}}\,\kappa_{\max})$.
Also $\tau(x;u)\ge r_{\mathrm{in}}$ by (1c).
Hence
\begin{equation*}
    \sum\frac{\lambda_v}{\tau_r(x;v)}\;\le\;\frac{\sqrt{\dim\mathcal{C}}}{r}+\sqrt{\dim\mathcal{C}}\cdot\frac{\sqrt{\dim\mathcal{C}}\,\kappa_{\max}}{\eta}\;=\;\frac{\sqrt{\dim\mathcal{C}}}{r}+\frac{\dim\mathcal{C}\kappa_{\max}}{\eta}.
\end{equation*}

\textit{Untilted: $c<\eta$, use the untilted set $U$.}
First $\rho$ is bounded below.
If $c\ge 0$, \cref{lem:curvature_feasibility_bounds}(2) gives $R\le\tau(x;d)\le 2r_{\mathrm{out}}c+\sqrt{2r_{\mathrm{out}}\rho}<2r_{\mathrm{out}}\eta+\sqrt{2r_{\mathrm{out}}\rho}$; since $R=8r_{\mathrm{out}}\eta$ this forces $\sqrt{2r_{\mathrm{out}}\rho}>6r_{\mathrm{out}}\eta$, i.e.
\begin{equation*}
    \rho \;>\; 18\,r_{\mathrm{out}}\eta^2.
\end{equation*}
If $c<0$, the same lemma gives $\tau(x;d)\le\sqrt{2r_{\mathrm{out}}\rho}$, so $R\le\sqrt{2r_{\mathrm{out}}\rho}$ and therefore
\begin{equation*}
    \rho \;\ge\; 32\,r_{\mathrm{out}}\eta^2.
\end{equation*}
In either subcase $\rho>18\,r_{\mathrm{out}}\eta^2$, so $U$ is the full set of $2\dim\mathcal{C}$ signed basis vectors and lies in $D_{\mathcal{C}}(x)$.
Decompose $d$ with $\lambda_{\mathrm{sgn}(c)u}=|c|$, $\lambda_{\mathrm{sgn}(a_i)w_i}=|a_i|$; $\sum\lambda=|c|+\sum|a_i|\le\sqrt{\dim\mathcal{C}}$, hence $A=\sqrt{\dim\mathcal{C}}$.
For the $B$-sum:
\begin{itemize}
    \item $\inner{\pm w_i}{u}=0$, so \cref{lem:curvature_feasibility_bounds}(1b) gives $\tau(x;\pm w_i)\ge\sqrt{r_{\mathrm{in}}\rho}>\eta\sqrt{18\,r_{\mathrm{in}}r_{\mathrm{out}}}$.
    Hence $\sum_i|a_i|/\tau(x;\pm w_i)\le\sqrt{\dim\mathcal{C}}\,\sqrt{\kappa_{\max}\kappa_{\min}}/(\sqrt{18}\,\eta)\le\sqrt{\dim\mathcal{C}}\,\kappa_{\max}/(4\eta)$.
    \item If $c\ge 0$: $\tau(x;u)\ge r_{\mathrm{in}}$ so $c/\tau(x;u)\le\eta\kappa_{\max}\le\kappa_{\max}$.
    \item If $c<0$: $\tau(x;-u)=\rho$ and \cref{lem:curvature_feasibility_bounds}(3) gives $|c|\le\rho/\tau(x;d)\le\rho/R$, so $|c|/\tau(x;-u)\le 1/R=\kappa_{\min}/(8\eta)$.
\end{itemize}
The last two subcases are mutually exclusive, so their contributions enter as a maximum.
Adding the $1/r$ truncation $\sum\lambda/r\le\sqrt{\dim\mathcal{C}}/r$,
\begin{align*}
    \sum\frac{\lambda_v}{\tau_r(x;v)}
    \;&\le\;
    \frac{\sqrt{\dim\mathcal{C}}}{r}+\frac{\sqrt{\dim\mathcal{C}}\,\kappa_{\max}}{4\eta}+\max\Bigl\{\kappa_{\max},\,\frac{\kappa_{\min}}{8\eta}\Bigr\} \\
    \;&\le\;
    \frac{\sqrt{\dim\mathcal{C}}}{r}+\frac{(\dim\mathcal{C}+1)\kappa_{\max}}{\eta}
    \;\le\;
    \frac{\sqrt{\dim\mathcal{C}}}{r}+\frac{2\dim\mathcal{C}\kappa_{\max}}{\eta}.
\end{align*}
The second inequality uses $\sqrt{\dim\mathcal{C}}/4\le \dim\mathcal{C}$ and $\max\{\kappa_{\max},\,\kappa_{\min}/(8\eta)\}\le\kappa_{\max}/\eta$, the latter since $\eta\le\tfrac12$ and $\kappa_{\min}\le\kappa_{\max}$; the third uses $\dim\mathcal{C}+1\le 2\dim\mathcal{C}$.

Thus $A=\sqrt{\dim\mathcal{C}}$ in each case.
The first-branch bound is $B\le\sqrt{\dim\mathcal{C}}(1/r+2\kappa_{\max})$, independent of $\eta$;
the inward-tilt bound is $B\le\sqrt{\dim\mathcal{C}}/r+\dim\mathcal{C}\kappa_{\max}/\eta$;
the untilted bound is $B\le\sqrt{\dim\mathcal{C}}/r+2\dim\mathcal{C}\kappa_{\max}/\eta$.
The last of these is the largest, and is the bound of the proposition.
\end{proof}

\begin{proof}[Proof of \cref{lem:delta-R-upgrade}]
(1) Every unit $d$ has $\tau(x;d)\ge\rho\ge r_{\mathrm{in}}\ge R$, so no direction is excluded.

(2)--(3), inward case.
Let $d$ be a unit feasible direction with $c:=\inner{d}{u}\ge 0$, put $s:=R\kappa_{\max}\in(0,1]$ and $w:=(d+su)/\norm{d+su}$.
Since $\norm{d+su}^2=1+2sc+s^2$ and $c\ge 0$, we have $\norm{d+su}\in[1,1+s]$, hence
\begin{equation*}
    \inner{w}{u}\;=\;\frac{c+s}{\norm{d+su}}\;\ge\;\frac{s}{1+s}\;\ge\;\frac{s}{2}\;=\;\frac{R\kappa_{\max}}{2}.
\end{equation*}
By \cref{lem:curvature_feasibility_bounds}(1a), $\tau(x;w)\ge\min\{2\inner{w}{u},1\}\,r_{\mathrm{in}}$.
If $\inner{w}{u}\ge\tfrac12$ this is at least $r_{\mathrm{in}}\ge R$; if $\inner{w}{u}<\tfrac12$ it is at least $2\inner{w}{u}\,r_{\mathrm{in}}\ge R$.
Thus $w$ is admissible in $\Delta^{[R]}$.
Also $\norm{w-d}\le\bigl\lvert 1-\norm{d+su}\bigr\rvert+s\le 2s=2R\kappa_{\max}$, so by \cref{lem:h-prime-lipschitz-span}(1)
\begin{equation*}
    -h'(x;d)\;\le\;-h'(x;w)+\Lambda\norm{d-w}\;\le\;\Delta^{[R]}(x)+2\Lambda\kappa_{\max}R.
\end{equation*}
At a boundary point $\cl{D_{\mathcal{C}}(x)}=\{d:\inner{d}{u}\ge 0\}$, so this covers every feasible direction and gives (2).

(3), outward case: $\mathcal{C}$ lies in the supporting halfspace at the nearest boundary point, so $\rho+t\inner{d}{u}\ge 0$ is necessary for feasibility, which is \cref{lem:curvature_feasibility_bounds}(3).
If $\tau(x;d)\ge R$ then $d$ is included in $\Delta^{[R]}(x)$; otherwise $\tau(x;d)<R$, so $d$ is excluded from $\Delta^{[R]}(x)$.
\end{proof}

\subsection{Proof of \texorpdfstring{\cref{thm:main}}{Theorem 6.2}}\label{sec:appendix_proof_thm_main}
The proof of \cref{thm:main} follows that of \citep[Theorem~4.1]{greenstein2026achieving} and combines four ingredients. First, the sufficient decrease of the line searches forces the exploration step sizes $t_{\mathrm{ex}}^k$ to vanish. Second, almost surely every direction in the relative interior of the feasible-direction cone at an accumulation point is approximated by exploration directions along some subsequence (\cref{lem:direction_exists}). Third, along such a subsequence the approximating directions admit a uniform lower bound on the feasibility radius (\cref{lem:uniform_step}). Fourth, vanishing steps together with that bound force non-descent in the limit (\cref{lem:non_descent}). This gives $h'(x^*;\,v)\ge 0$ on $\ri{D_{\mathcal{C}}(x^*)}$, which \cref{cor:extension} extends to the whole cone.
From this point onward, let
\begin{equation*}
L := \aff{\mathcal{C}-\bar{x}},
\end{equation*}
for any fixed $\bar{x}\in\mathcal{C}$, so $L$ is the linear subspace parallel to $\aff{\mathcal{C}}$.
This definition is independent of the choice of $\bar{x}$, because
$\aff{\mathcal{C}-\bar{x}}=\aff{\mathcal{C}-\bar{y}}$ for all $\bar{x},\bar{y}\in\mathcal{C}$
(\Cref{lem:foundational properties}).
In particular, $\aff{\mathcal{C}}-x^k=L$ for every iterate $x^k\in\mathcal{C}$.
We write $\mathbb{S}_L$ for the unit sphere of $L$.


\begin{lem}[Vanishing line-search steps; cf.\ the proof of {\citep[Theorem~4.1(ii)]{greenstein2026achieving}}]\label{lem:non_descent}
    Let $\{x_j\}_{j \in \mathcal{J}}$ be a sequence in $\mathcal{C}$ converging to $x^* \in \mathcal{C}$,
    and let $\{w^j\}_{j \in \mathcal{J}}$ be a sequence of unit vectors in~$L$
    converging to some $w^* \in \real^n$.
    Suppose that:
    \begin{enumerate}[label=(\roman*),nosep]
        \item there exists $c > 0$ such that
              $\tau(x_j;\, w^j) \geq c$ for all $j \in \mathcal{J}$;
        \item $t^j \to 0$ as $j \to \infty$ over $\mathcal{J}$, where for each $j \in \mathcal{J}$,
        \[
            t^j \in \argmin_{t \in [0,\, \min\{\tau(x_j;\, w^j),\, r\}]}
            h(x_j + tw^j) + \tfrac{\gamma}{2}\,t^2,
        \]
    \end{enumerate}
    where $\gamma>0$. 
    Then $h'(x^*;\, w^*) \geq 0$.
\end{lem}

\begin{proof}
    Assume for contradiction that $h'(x^*; w^*) = -\epsilon < 0$.
    By the definition of the directional derivative,
    there exists $\tau_0 > 0$ such that
    \begin{equation}\label{eq:sufficient_descent}
        h(x^* + tw^*) - h(x^*) < -\tfrac{\epsilon}{2}\, t
        \qquad \forall\, t \in (0, \tau_0].
    \end{equation}
    By local Lipschitz continuity of $h$ at $x^*$, there exist $R_h>0$ and $L_h>0$
    such that $h$ is $L_h$-Lipschitz on $B(x^*,R_h)$.

    Define
    \[
        t^* = \min\Bigl\{c,\; \tau_0,\;
        \frac{\epsilon}{4\gamma},\; r,\; \frac{R_h}{4}\Bigr\}.
    \]
    By Assumption (i), for all $j$ sufficiently large,
    $t^* \leq c \leq \tau(x_j; w^j)$ and $t^* \leq r$, so 
    $t^* \in [0,\, \min\{\tau(x_j; w^j),\, r\}]$.

    Since $x_j \to x^*$ and $w^j \to w^*$,
    there exists $J_1 \in \mathbb{N}$ such that for $j > J_1$, $j \in \mathcal{J}$:
    \[
        \norm{x_j-x^*}<\min\!\left\{\frac{R_h}{4},\frac{\epsilon t^*}{16L_h}\right\},
        \qquad
        \norm{w^j-w^*}<\frac{\epsilon}{16L_h}.
    \]
    For such $j$, we have
    \[
        \norm{(x^*+t^*w^*)-x^*}=t^*\le\frac{R_h}{4}<R_h,
    \]
    and, since $\norm{w^j}=1$,
    \[
        \norm{(x_j+t^*w^j)-x^*}
        \le\norm{x_j-x^*}+t^*\norm{w^j}
        <\frac{R_h}{4}+\frac{R_h}{4}
        =\frac{R_h}{2}<R_h.
    \]
    Hence both points $x_j+t^*w^j$ and $x^*+t^*w^*$ lie in $B(x^*,R_h)$, where
    the $L_h$-Lipschitz bound is valid, and therefore
    \[
        |h(x_j+t^*w^j)-h(x^*+t^*w^*)|
        \leq L_h\norm{(x_j-x^*)+t^*(w^j-w^*)},
    \]
    and, by the triangle inequality together with the bounds above,
    \[
        |h(x_j + t^* w^j) - h(x^* + t^* w^*)| \leq
        \frac{\epsilon t^*}{8},
        \qquad
        |h(x^*) - h(x_j)| \leq \frac{\epsilon t^*}{8}.
    \]
    For every $t\in[0,\epsilon/(4\gamma)]$, we have
    \begin{equation*}
        \frac{\gamma t^2}{2}\leq\frac{\epsilon t}{8}.
    \end{equation*}
    Combining with~\eqref{eq:sufficient_descent} at $t = t^*$
    (valid since $t^* \leq \tau_0$):
    \begin{align*}
        h(x_j + t^* w^j) - h(x_j) + \frac{\gamma(t^*)^2}{2}
        &= \bigl[h(x_j + t^* w^j) - h(x^* + t^* w^*)\bigr]
         + \bigl[h(x^* + t^* w^*) - h(x^*)\bigr] \\
        &\quad + \bigl[h(x^*) - h(x_j)\bigr] + \frac{\gamma(t^*)^2}{2} \\
        &\leq \frac{\epsilon t^*}{8} - \frac{\epsilon t^*}{2}
           + \frac{\epsilon t^*}{8} + \frac{\epsilon t^*}{8}
        \leq -\frac{\epsilon t^*}{8}.
    \end{align*}

    Since $t^j \to 0$ by Assumption (ii), there exists $J_2\in\mathbb{N}$ such that for
    $j>J_2$, $j\in\mathcal{J}$,
    \[
        t^j<\min\!\left\{\frac{R_h}{4},\frac{\epsilon t^*}{8L_h},\frac{\epsilon}{4\gamma}\right\}.
    \]
    For such $j$ (and $j>J_1$), we have $\norm{x_j-x^*}<R_h/4$ and $\norm{w^j}=1$, so
    \[
        \norm{(x_j+t^jw^j)-x^*}
        \le \norm{x_j-x^*}+t^j\norm{w^j}
        < \frac{R_h}{4}+\frac{R_h}{4}
        = \frac{R_h}{2}<R_h.
    \]
    Thus both $x_j+t^jw^j$ and $x_j$ belong to $B(x^*,R_h)$, and the Lipschitz bound gives
    \[
        h(x_j+t^jw^j)-h(x_j)\ge -L_h t^j,
    \]
    Since $\frac{\gamma(t^j)^2}{2}\ge0$, we obtain
    \[
        h(x_j+t^jw^j)-h(x_j)+\frac{\gamma(t^j)^2}{2}
        \ge h(x_j+t^jw^j)-h(x_j)
        \ge -L_h t^j
        > -\frac{\epsilon t^*}{8}.
    \]
    Since $t^* \in [0, \min\{\tau(x_j; w^j), r\}]$, optimality of $t^j$ requires
    $h(x_j + t^j w^j) + \tfrac{\gamma}{2}(t^j)^2
    \leq h(x_j + t^* w^j) + \tfrac{\gamma}{2}(t^*)^2$.
    But the estimates above give
    $h(x_j + t^j w^j) + \tfrac{\gamma}{2}(t^j)^2
    > h(x_j) - \epsilon t^*/8
    \geq h(x_j + t^* w^j) + \tfrac{\gamma}{2}(t^*)^2$,
    a contradiction.
\end{proof}
Applying \Cref{lem:non_descent} requires a uniform lower bound on the feasibility
radius $\tau(x_j;\, w^j)$ along the sequence of directions.
The following lemma provides such a bound whenever the limit direction $w^*$ lies in the relative interior of $D_{\mathcal{C}}(x^*)$.


\begin{lem}[Uniform step size bound; cf.\ {\citep[Lemma~4.3]{greenstein2026achieving}}]\label{lem:uniform_step}
    Let $\{x_j\}_{j \in \mathcal{J}} \subseteq \mathcal{C}$ with $x_j \to x^*$, and let
    $\{w^j\}_{j \in \mathcal{J}} \subseteq L$ with
    $w^j \to w^*$, where $w^* \in \ri{D_\mathcal{C}(x^*)}$.
    Then there exists $\bar{\epsilon} > 0$ and $J_0 \in \mathbb{N}$ such that
    \[
        x_j + \bar{\epsilon}\, w^j \in \mathcal{C}
        \qquad \forall\, j \in \mathcal{J},\; j \geq J_0.
    \]
    In particular, $w^j \in D_\mathcal{C}(x_j)$ and
    $\tau(x_j; w^j) \geq \bar{\epsilon}$ for all
    $j \in \mathcal{J}$ sufficiently large.
\end{lem}

\begin{proof}
    \noindent Since $w^* \in \ri{D_\mathcal{C}(x^*)}$ and $\aff{D_\mathcal{C}(x^*)} = L$, there exists $r > 0$
    such that $\mathcal{B}(w^*, r) \cap L \subseteq D_\mathcal{C}(x^*)$.
    Choose $\dim L+1$ directions $d_0, \ldots, d_{\dim L} \in \mathcal{\mathcal{B}}(w^*, r) \cap L$ such that
    $\conv\{d_0, \ldots, d_{\dim L}\}$ is a $\dim L$-dimensional simplex in $L$ with
    $w^* \in \ri{\conv\{d_0, \ldots, d_{\dim L}\}}$.
    Since each $d_i \in D_\mathcal{C}(x^*)$, by definition, there exist $\epsilon_i > 0$ with
    $x^* + \epsilon_i d_i \in \mathcal{C}$.
    Set $\bar{\epsilon} = \min_i \epsilon_i > 0$.
    Since $\bar{\epsilon} \leq \epsilon_i$ and $\mathcal{C}$ is convex (with $x^* \in \mathcal{C}$),
    each $p_i := x^* + \bar{\epsilon}\, d_i \in \mathcal{C}$.

    For any $d \in \conv\{d_0, \ldots, d_{\dim L}\}$, write
    $d = \sum_{i=0}^{\dim L} \lambda_i d_i$ with $\lambda_i \geq 0$,
    $\sum_{i=0}^{\dim L} \lambda_i = 1$.
    Then $x^* + \bar{\epsilon}\, d = \sum_{i=0}^{\dim L} \lambda_i p_i \in \mathcal{C}$
    by convexity.

    Since $w^* \in \ri{\conv\{d_0, \ldots, d_{\dim L}\}}$, there exists $\rho > 0$
    with $\mathcal{B}(w^*, \rho) \cap L \subseteq \conv\{d_0, \ldots, d_{\dim L}\}$.
    In particular, $x^* + \bar{\epsilon}\, d \in \mathcal{C}$ for every
    $d \in \mathcal{B}(w^*, \rho) \cap L$, so
    $\mathcal{B}(x^* + \bar{\epsilon}\, w^*,\; \bar{\epsilon}\rho) \cap \aff{\mathcal{C}}
    \subseteq \mathcal{C}$,
    i.e., $x^* + \bar{\epsilon}\, w^* \in \ri{\mathcal{C}}$.

    Since $w^j \in L = \aff{\mathcal{C} - x^*}$,
    we have $x_j + \bar{\epsilon}\, w^j \in \aff{\mathcal{C}}$.
    Moreover,
    \[
        \|(x_j + \bar{\epsilon}\, w^j) - (x^* + \bar{\epsilon}\, w^*)\|
        \leq \|x_j - x^*\| + \bar{\epsilon}\,\|w^j - w^*\|
        < \bar{\epsilon}\rho
    \]
    for all $j$ sufficiently large.
    Therefore $x_j + \bar{\epsilon}\, w^j \in \mathcal{C}$ for all $j \geq J_0$
    for some $J_0 \in \mathbb{N}$.
\end{proof}
\begin{remark}\label{rem:why_ri}
    The condition $w^* \in \ri{D_\mathcal{C}(x^*)}$ cannot be relaxed to
    $w^* \in D_\mathcal{C}(x^*)$ in general.
    For directions on the relative boundary of $D_\mathcal{C}(x^*)$,
    $\tau(x_j; w^j)$ may tend to zero.
    When $x^* \in \ri{\mathcal{C}}$, this issue does not arise:
    $D_\mathcal{C}(x^*) = \aff{\mathcal{C} - x^*}$ is a subspace, so
    $\ri{D_\mathcal{C}(x^*)} = D_\mathcal{C}(x^*)$ and $w^*\in\ri{D_\mathcal{C}(x^*)}$ holds
    for every limit direction.
\end{remark}
Fix a countable dense set $\mathcal{Q}\subseteq\mathbb{S}_L$.
The following lemma establishes that, almost surely, every direction from $\mathcal{Q}$
that is feasible at an accumulation point can be realized as a limit of exploration directions.
\begin{lem}[Exploration directions approximate all feasible directions; adapted from {\citep[Lemma~4.2]{greenstein2026achieving}}]\label{lem:direction_exists}
    Almost surely, for every accumulation point~$x^*$ of~$\{x_k\}$
    and every $v \in \mathcal{Q}\cap \ri{D_\mathcal{C}(x^*)}$,
    there exists a subsequence $\mathcal{J} \subseteq \mathbb{N}$ such that
    $x_j \to x^*$ and $v_{\mathrm{ex}}^{j} \to v$ as
    $j \to \infty$, $j \in \mathcal{J}$.
\end{lem}
\begin{proof}
    Let $\mathcal{F}_{k-1}$ denote the $\sigma$-algebra generated by
    $\{x^0,v_{\mathrm{ex}}^0,\dots,x^{k-1},v_{\mathrm{ex}}^{k-1},x^k\}$.
    Then $x^k$ is $\mathcal{F}_{k-1}$-measurable and $v_{\mathrm{ex}}^k$ is
    independent of $\mathcal{F}_{k-1}$.

    Fix $y\in\mathcal{C}$, $\delta>0$, and
    $V:=\mathcal{B}[w,\epsilon]\cap\mathbb{S}_L$ with
    $w\in\mathbb{S}_L$, $\epsilon>0$.
    Since the sampling law of $v_{\mathrm{ex}}^k$ has full support on
    $\mathbb{S}_L$, the probability
    $p:=P(v_{\mathrm{ex}}^k\in V)$ satisfies $p>0$ and is independent of $k$.
    Define
    \begin{equation*}
        A_k:=\{v_{\mathrm{ex}}^k\in V\}\cap\{\norm{x^k-y}<\delta\}.
    \end{equation*}
    Then
    \begin{equation*}
        P(A_k\mid\mathcal{F}_{k-1})
        =\mathbf{1}_{\{\norm{x^k-y}<\delta\}}\cdot p.
    \end{equation*}
    By the conditional Borel--Cantelli--L\'evy lemma
    \citep[Chapter~VII, \S5, Corollary~2]{shiryaev2019probability2},
    \begin{equation}\label{eq:borel_cantelli_levy_result}
        \Bigl\{\sum_{k=0}^{\infty}P(A_k\mid\mathcal{F}_{k-1})=\infty\Bigr\}
        =\{A_k\text{ i.o.}\}\quad\text{a.s.}
    \end{equation}
    Hence
    \begin{equation*}
        \Bigl\{\norm{x^k-y}<\delta\ \text{i.o.}\Bigr\}
        \Longrightarrow
        \Bigl\{v_{\mathrm{ex}}^k\in V,\ \norm{x^k-y}<\delta\ \text{i.o.}\Bigr\}
        \quad\text{a.s.}
    \end{equation*}

    Let $\{y_n\}_{n\in\N}$ be a countable dense subset of $\mathcal{C}$, let
    $\mathcal{Q}=\{q_i\}_{i\in\N}$, and let $\delta_m,\epsilon_l$ range over positive rationals.
    For each tuple $u=(n,i,m,l)$, define the logical implication event
    \begin{equation*}
        E_u:=
        \Bigl\{\norm{x^k-y_n}<\delta_m\ \text{i.o.}\Longrightarrow
        \bigl(v_{\mathrm{ex}}^k\in\mathcal{B}[q_i,\epsilon_l],\ \norm{x^k-y_n}<\delta_m\bigr)\ \text{i.o.}\Bigr\}.
    \end{equation*}
    By \eqref{eq:borel_cantelli_levy_result}, $P(E_u)=1$ for each $u$.
    Set
    \begin{equation*}
        \Omega_0:=\bigcap_u E_u.
    \end{equation*}
    Since the index set is countable, $P(\Omega_0)=1$.

    Fix $\omega\in\Omega_0$. If $\{x^k(\omega)\}$ has no accumulation point,
    the lemma is vacuous on this path. Otherwise, let $x^*$ be any accumulation point
    of $\{x^k(\omega)\}$ and fix
    $v\in\ri{D_\mathcal{C}(x^*)}\cap\mathbb{S}_L$.
    Write $x^k:=x^k(\omega)$ and $v_{\mathrm{ex}}^k:=v_{\mathrm{ex}}^k(\omega)$.

    Let $\eta>0$. By density, choose $n,m,i,l$ such that
    \begin{equation*}
        \norm{y_n-x^*}<\eta/4,\qquad
        \delta_m\in(\eta/2,3\eta/4),\qquad
        \norm{q_i-v}<\eta/4,\qquad
        \epsilon_l<\eta/4.
    \end{equation*}
    Since $x^*$ is an accumulation point, $\norm{x^k-x^*}<\eta/4$ for infinitely many $k$.
    For such $k$,
    \begin{equation*}
        \norm{x^k-y_n}\le\norm{x^k-x^*}+\norm{x^*-y_n}<\eta/2<\delta_m.
    \end{equation*}
    Therefore $\norm{x^k-y_n}<\delta_m$ occurs infinitely often.
    Because $\omega\in E_u$ for $u=(n,i,m,l)$, there are infinitely many such $k$ with
    $v_{\mathrm{ex}}^k\in\mathcal{B}[q_i,\epsilon_l]$. For these indices,
    \begin{equation*}
        \norm{x^k-x^*}\le\norm{x^k-y_n}+\norm{y_n-x^*}<\delta_m+\eta/4<\eta,
    \end{equation*}
    and
    \begin{equation*}
        \norm{v_{\mathrm{ex}}^k-v}\le\norm{v_{\mathrm{ex}}^k-q_i}+\norm{q_i-v}<\epsilon_l+\eta/4<\eta/2<\eta.
    \end{equation*}

    Let $\eta_q\downarrow0$. For each $q$, let $K_q$ be the infinite set of indices
    satisfying $\norm{x^k-x^*}<\eta_q$ and $\norm{v_{\mathrm{ex}}^k-v}<\eta_q$.
    Choose inductively $k_q\in K_q$ with $k_q>k_{q-1}$. Then
    $x^{k_q}\to x^*$ and $v_{\mathrm{ex}}^{k_q}\to v$.
\end{proof}
The following lemma records the analytic properties of the
directional derivative that the proof of \Cref{thm:main} relies on.

\begin{lem}[Properties of the directional derivative; cf.\ {\citep[Lemmas~1.1 and~4.4]{greenstein2026achieving}}]\label{lem:h_prime_properties}
    The function
    \[
        d \mapsto h'(x^*; d)
        = \langle \nabla f(x^*), d \rangle - g'(x^*; d)
    \]
    is continuous, positively homogeneous of degree~$1$, and concave on $\real^n$.
\end{lem}

\begin{proof}
    Since $g$ is convex and finite-valued on a neighborhood of $x^*$,
    the directional derivative $g'(x^*; d)$ exists and is finite for every
    $d \in \real^n$ \cite[Theorem~23.1 and Theorem~23.4]{rockafellar1997convex}.

    Since map $d \mapsto g'(x^*; d)$ is convex and positively homogeneous of
    degree~$1$ \cite[Theorem~23.1]{rockafellar1997convex}, and is finite
    everywhere on $\real^n$, it is continuous
    \cite[Theorem~10.1]{rockafellar1997convex}.
    The map $d \mapsto \langle \nabla f(x^*), d \rangle$ is linear.
    Therefore $h'(x^*; \cdot) = \langle \nabla f(x^*), \cdot \rangle - g'(x^*; \cdot)$
    is continuous (linear minus continuous), positively homogeneous (linear minus
    positively homogeneous), and concave (linear minus convex).
\end{proof}
Continuity of $h'(x^*;\, \cdot)$ allows the non-descent inequality to be extended
beyond $\ri{D_\mathcal{C}(x^*)}$ to its closure $\cl{D_\mathcal{C}(x^*)}$, which is the cone
relevant for $d$-stationarity.
\begin{corollary}[Extension from relative interior to tangent cone; cf.\ {\citep[Lemma~4.4]{greenstein2026achieving}}]\label{cor:extension}
    If $h'(x^*;\, v) \geq 0$ for all $v \in \ri{D_\mathcal{C}(x^*)}$,
    then $h'(x^*;\, d) \geq 0$ for all $d \in \cl{D_\mathcal{C}(x^*)}$.
\end{corollary}

\begin{proof}
    Let $d \in \cl{D_\mathcal{C}(x^*)}$.
    Since $\cl{\ri{D_\mathcal{C}(x^*)}} = \cl{D_\mathcal{C}(x^*)}$ (by \cref{lem:foundational properties}),
    there exists a sequence $\{v_j\} \subseteq \ri{D_\mathcal{C}(x^*)}$ with
    $v_j \to d$.
    By continuity of the directional derivative (\cref{lem:h_prime_properties}):
    \[
        h'(x^*;\, d)
        = \lim_{j \to \infty} h'(x^*;\, v_j)
        \geq 0.
    \]
\end{proof}


We can now prove \Cref{thm:main}.
\begin{proof}[Proof of \Cref{thm:main}]
    By \Cref{lem:direction_exists}, there exists an event $\Omega_0$ with
    $P(\Omega_0)=1$ such that, for every $\omega\in\Omega_0$, every accumulation
    point $x^*$ of $\{x_k(\omega)\}$, and every
    $v\in\mathcal{Q}\cap\ri{D_\mathcal{C}(x^*)}$, there exists a subsequence
    $\mathcal{J}\subseteq\N$ with $x_j(\omega)\to x^*$ and $v_{\mathrm{ex}}^j(\omega)\to v$.
    Fix $\omega\in\Omega_0$ and work pathwise.
    Let $x^*$ be an arbitrary accumulation point of $\{x_k(\omega)\}$.
    As in \citep[Theorem~4.1(i)]{greenstein2026achieving}, since the algorithm selects the iterate with the largest decrease,
    comparing against $t = 0$ gives
    \[
        h(x_k) - h(x_{k+1})
        \geq h(x_k) - h(x_k + t_{\mathrm{ex}}^k v_{\mathrm{ex}}^k)
        \geq \tfrac{\gamma}{2} (t_{\mathrm{ex}}^k)^2.
    \]
    Summing over all $k \geq 0$ yields
    \[
        \tfrac{\gamma}{2} \sum_{k=0}^{\infty} (t_{\mathrm{ex}}^k)^2
        \leq \sum_{k=0}^{\infty} \bigl[h(x_k) - h(x_{k+1})\bigr]
        \leq h(x_0) - \inf_{x \in \mathcal{C}} h(x) < \infty,
    \]
    and therefore
    \begin{equation}\label{eq:tex_vanish}
        t_{\mathrm{ex}}^k \to 0 \quad \text{as } k \to \infty.
    \end{equation}
 
    If $\ri{D_\mathcal{C}(x^*)}\cap\mathbb{S}_L=\emptyset$, then
    $\ri{D_\mathcal{C}(x^*)}=\{0\}$ and $h'(x^*;\,0)=0$, so
    $h'(x^*;\,d)\ge0$ for all $d\in\ri{D_\mathcal{C}(x^*)}$; the remainder is immediate.
    Otherwise, let $\{v_i\}_{i=1}^\infty$ be an enumeration of
    $\mathcal{Q}\cap\ri{D_\mathcal{C}(x^*)}$.
    Since $\ri{D_\mathcal{C}(x^*)}\cap\mathbb{S}_L$ is relatively open in
    $\mathbb{S}_L$ and $\mathcal{Q}$ is dense in $\mathbb{S}_L$,
    the set $\{v_i\}_{i=1}^\infty$ is dense in
    $\ri{D_\mathcal{C}(x^*)}\cap\mathbb{S}_L$.
    For each~$v_i$, \Cref{lem:direction_exists} provides a
    subsequence $\mathcal{J}_i \subseteq \mathbb{N}$ with $x_j \to x^*$
    and $v_{\mathrm{ex}}^j \to v_i$.
    Since $v_i \in \ri{D_\mathcal{C}(x^*)}$ and $v_{\mathrm{ex}}^j \in L$
    by construction, \Cref{lem:uniform_step} provides $\bar{\epsilon}_i > 0$
    such that $\tau(x_j;\, v_{\mathrm{ex}}^j) \geq \bar{\epsilon}_i$ for all
    $j \in \mathcal{J}_i$ sufficiently large.
    The step sizes satisfy $t_{\mathrm{ex}}^j \to 0$ by~\eqref{eq:tex_vanish}.
    Applying \Cref{lem:non_descent} yields $h'(x^*;\, v_i) \geq 0$ for all $i\ge1$.
    By continuity of $h'(x^*;\, \cdot)$ (\Cref{lem:h_prime_properties}) and density
    of $\{v_i\}$ in $\ri{D_\mathcal{C}(x^*)} \cap \mathbb{S}_L$, the
    inequality extends to all of
    $\ri{D_\mathcal{C}(x^*)} \cap \mathbb{S}_L$, and hence to all of
    $\ri{D_\mathcal{C}(x^*)}$ by positive homogeneity
    (\Cref{lem:h_prime_properties}).
    By \Cref{cor:extension}, $h'(x^*;\, d) \geq 0$ for all
    $d \in \cl{D_\mathcal{C}(x^*)}$, and $x^*$ is d-stationary.
    Since $x^*$ was an arbitrary accumulation point on this fixed $\omega\in\Omega_0$,
    every accumulation point of $\{x_k(\omega)\}$ is d-stationary.
    Because $P(\Omega_0)=1$, the claim follows.
\end{proof}

\section{Measure-Theoretic Foundations for \cref{lem:always_smooth_almost_surely}}\label{sec:appendix_mt_background}

This appendix makes precise the claim, used in
\cref{sec:noised_alg_properties}, that nonsmooth boundary points are
negligible in the relative boundary. It collects the measure-theoretic
and probabilistic background used there, and states and proves four
auxiliary results together with the full proof of
\cref{lem:always_smooth_almost_surely}:
\begin{enumerate}
    \item smoothness of $\mathcal{C}$ at a boundary point is preserved by
    an isometric flattening onto a full-dimensional convex set
    (\cref{prop:closed_convex_sets_are_almost_everywhere_smooth});
    \item $\mathcal{C}_{ns}$ is a Borel set and is $\mathcal{H}^{\dim\mathcal{C}-1}$-null
    (\cref{prop:nonsmooth_boundary_borel_null});
    \item radial projection from a point of $\mathcal{C}$ maps
    $\mathcal{H}^{\dim\mathcal{C}-1}$-null subsets of $\aff{\mathcal{C}}$ to $\sigma$-null
    subsets of the unit sphere $\mathbb{S}_L$
    (\cref{lem:radial_image_null});
    \item every search direction generated by \cref{alg:construct_directions}
    has a conditional law on $\mathbb{S}_L$ that is absolutely continuous
    with respect to $\sigma$ (\cref{lem:direction_density}).
\end{enumerate}
The full proof of \cref{lem:always_smooth_almost_surely} combines the
last three of these facts and appears at the end of the appendix.
Throughout, ``measure'' means a Borel measure unless stated otherwise, $L := \aff{\mathcal{C} - x}$ for any $x \in \mathcal{C}$ (independent of the choice of $x$ by \cref{lem:foundational properties}), and $\dim\mathcal{C}\ge 1$.

\paragraph{Notation.}
\begin{center}
\renewcommand{\arraystretch}{1.3}
\begin{tabular}{@{}ll@{}}
$\mathbb{S}_L := \{\omega \in L : \norm{\omega} = 1\}$
  & unit sphere of $L$ \\
$\Pi_x(y) := (y - x)/\norm{y - x}$
  & radial projection from $x$, for $y \in (x+L)\setminus\{x\}$ \\
$\mu\!\restriction\! A$
  & restriction: $(\mu\!\restriction\! A)(B) = \mu(A \cap B)$ \\
$\mu \ll \nu$
  & $\mu$ absolutely continuous with respect to $\nu$ \\
$\lambda_M$ (resp.\ $\lambda_L$)
  & Lebesgue measure on an affine (resp.\ linear) subspace \\
$\sigma:=\mathcal{H}^{\dim\mathcal{C}-1}\!\restriction\!\mathbb{S}_L$
  & surface measure on $\mathbb{S}_L$ (\cref{def:surface-measure}) \\
$T_\#\mu$
  & pushforward of $\mu$ under $T$ \\
$\mathcal{L}(X)$, $\mathcal{L}(X \mid \mathcal{G})$
  & law of $X$; conditional law of $X$ given $\mathcal{G}$
\end{tabular}
\end{center}

\paragraph{Lebesgue measure on an affine subspace.}

\begin{defin}\label{def:lebesgue-affine-subspace}
Let $M = p + L_0 \subseteq \real^n$ be an affine subspace of dimension $m$, with direction subspace $L_0$ and base point $p \in M$. For any linear isometry $\iota : \real^m \to L_0$, the \textit{$m$-dimensional Lebesgue measure on $M$}, denoted $\lambda_M$, is the pushforward of $\lambda_m$ on $\real^m$ under $x \mapsto p + \iota(x)$. It is independent of the choice of $p$ and $\iota$, and satisfies $\lambda_M = \mathcal{H}^m\!\restriction\! M$. We write $\lambda_L$ when $M = L$ is linear.

A measure or random vector \textit{has a density on $M$} (equivalently, \textit{is absolutely continuous on $M$}) means it is absolutely continuous with respect to $\lambda_M$.
\end{defin}

\paragraph{Absolute continuity, density, and pushforward.}

\begin{defin}\label{def:absolute-continuity-density}
For measures $\mu, \nu$ on the same measurable space, $\mu \ll \nu$ iff $\nu(B) = 0 \Rightarrow \mu(B) = 0$ for every measurable $B$. When $\mu, \nu$ are $\sigma$-finite, this is equivalent (Radon--Nikodym; \citep{RoydenFitzpatrickRealAnalysis4e}) to the existence of a \textit{density} $p = \tfrac{d\mu}{d\nu} \ge 0$ with $\mu(B) = \int_B p\,d\nu$.
\end{defin}

\begin{defin}\label{def:pushforward-law}
For a measurable map $T$ and a measure $\mu$, the \textit{pushforward} is $T_\#\mu(B) := \mu(T^{-1}(B))$. For a random vector $X$ on $(\Omega, \mathcal{F}, \mathbb{P})$, its \textit{law} is $\mathcal{L}(X) := X_\#\mathbb{P}$; ``$X$ has a density'' means $\mathcal{L}(X) \ll \lambda$.
\end{defin}

\begin{proposition}\label{prop:affine-image-abs-cont}
Let $M$ and $M'$ be affine subspaces of $\real^n$ and let $T : M \to M'$
be a surjective affine map. If $\mu \ll \lambda_M$, then $T_\#\mu \ll \lambda_{M'}$.
In particular, a bijective affine map $T : M \to M$ sends
$\lambda_M$-absolutely continuous measures to $\lambda_M$-absolutely continuous
measures.
\end{proposition}
\begin{proof}
Let $p = \dim M$ and $q = \dim M'$.
Via \cref{def:lebesgue-affine-subspace}, choose linear isometries to identify
$\lambda_M$ with $\lambda_p$ on $\real^p$ and $\lambda_{M'}$ with $\lambda_q$ on
$\real^q$; in these coordinates $T$ becomes an affine map $\hat{T}:\real^p\to\real^q$
with linear part $A:\real^p\to\real^q$ of rank $q$ (surjectivity of $T$ forces this).
Let $E:=(\ker A)^\perp\cong\real^q$, so the restriction $A|_E:E\to\real^q$ is
invertible. For any $B\subseteq\real^q$ with $\lambda_q(B)=0$ and translation
vector $b$ (so $\hat{T}(s)=As+b$),
\begin{equation*}
    \hat{T}^{-1}(B)
    = \bigl\{(u,v)\in\ker A\times E : A|_E(v) \in B-b\bigr\}
    = \ker A\;\times\;(A|_E)^{-1}(B-b).
\end{equation*}
Since $A|_E$ is invertible,
$\lambda_q\bigl((A|_E)^{-1}(B-b)\bigr) = \lvert\det(A|_E)\rvert^{-1}\lambda_q(B-b)=0$.
By Fubini (with $\lambda_p = \lambda_{\ker A}\otimes\lambda_E$),
$\lambda_p\bigl(\hat{T}^{-1}(B)\bigr)=0$, hence $T_\#\mu(B)=\mu(T^{-1}(B))=0$.
\end{proof}

\paragraph{Hausdorff measure and dimension.}

\begin{defin}\label{def:hausdorff-measure}
For $A \subseteq \real^n$, $s \ge 0$, and $\delta > 0$, set
\begin{equation*}
\mathcal{H}^s_\delta(A)
:= \inf\Bigl\{
  \textstyle\sum_{i=1}^\infty \alpha(s)
  \bigl(\tfrac{\diam U_i}{2}\bigr)^s
  : A \subseteq \bigcup_i U_i,\ \diam U_i \le \delta
\Bigr\},
\end{equation*}
and define the \textit{$s$-dimensional Hausdorff measure}
$\mathcal{H}^s(A) := \lim_{\delta \downarrow 0} \mathcal{H}^s_\delta(A)$.
The constant $\alpha(s) = \pi^{s/2}/\Gamma\!\left(\tfrac{s}{2}+1\right)$ normalizes so that $\mathcal{H}^m = \lambda_m$ on $\real^m$ for integer $m$; the normalizing constant is irrelevant for null-set statements. A set $A$ is \textit{$\mathcal{H}^s$-null} if $\mathcal{H}^s(A) = 0$.
\end{defin}

\begin{defin}\label{def:hausdorff-dimension}
The \textit{Hausdorff dimension} of $A$ is
$\dim_{\mathcal{H}}(A) := \inf\{s \ge 0 : \mathcal{H}^s(A) = 0\}$.
\end{defin}

\begin{proposition}[Lipschitz scaling {\citep[Theorem~2.8(i)]{EvansGariepy2015}}]\label{prop:lipschitz-hausdorff}
If $f : A \to \real^{n'}$ is Lipschitz with constant $\operatorname{Lip}(f)$, then for every $s \ge 0$,
\begin{equation*}
\mathcal{H}^s\bigl(f(A)\bigr) \le \operatorname{Lip}(f)^{\,s}\,\mathcal{H}^s(A).
\end{equation*}
Consequently, Lipschitz maps do not increase Hausdorff dimension and map $\mathcal{H}^s$-null sets to $\mathcal{H}^s$-null sets.
\end{proposition}

\paragraph{Surface measure and polar coordinates.}

\begin{defin}\label{def:surface-measure}
For a linear subspace $L$ with $m = \dim L$, the \textit{surface measure} on $\mathbb{S}_L = \{\omega \in L : \norm{\omega} = 1\}$ is $\sigma := \mathcal{H}^{m-1}\!\restriction\!\mathbb{S}_L$. Under any linear isometry $L \cong \real^m$, $\sigma$ is the standard spherical measure on $S^{m-1}$.
\end{defin}

\begin{proposition}[Polar-coordinate factorization; adapted from {\citep[Theorem~3.12]{EvansGariepy2015}}]\label{prop:polar-coordinates}
With $m = \dim L$ and $\sigma$ as above, for every nonnegative Borel function $g$,
\begin{equation*}
\int_L g\,d\lambda_L
= \int_0^\infty\!\!\int_{\mathbb{S}_L} g(r\omega)\,r^{m-1}\,d\sigma(\omega)\,dr.
\end{equation*}
\end{proposition}
Theorem~3.12 in~\citep{EvansGariepy2015} states the formula for $\real^n$ with the inner
integral written as $\int_{\partial B(r)} g\, d\mathcal{H}^{m-1}$; the form above
follows by identifying $L \cong \real^m$ via a linear isometry (which preserves
$\mathcal{H}^s$ by \citep[Theorem~2.2(v)]{EvansGariepy2015}) and by the substitution $x = r\omega$,
which yields $d\mathcal{H}^{m-1}|_{\partial B(r)} = r^{m-1}\, d\sigma(\omega)$
via the scaling property \citep[Theorem~2.2(iv)]{EvansGariepy2015}.

\begin{lem}\label{lem:radial-projection-ac}
Let $R : L \setminus \{0\} \to \mathbb{S}_L$, $R(z) = z/\norm{z}$. If a measure $\nu$ on $L$ satisfies $\nu \ll \lambda_L$, then $R_\#\nu \ll \sigma$.
\end{lem}
\begin{proof}
Let $B \subseteq \mathbb{S}_L$ with $\sigma(B) = 0$. Then $R^{-1}(B) = \{r\omega : r > 0,\,\omega \in B\}$.
Note that $\{0\} \notin R^{-1}(B)$. For $N\in\mathbb{N}$ set $A_N := \{z \in L : 1/N \le \norm{z} \le N\}$; then
$R^{-1}(B) = \bigcup_{N\ge 1} \bigl(R^{-1}(B)\cap A_N\bigr)$. For each $N$,
by \cref{prop:polar-coordinates} applied with $g = \mathbf{1}_{\{r\omega: \omega\in B,\, 1/N\le r\le N\}}$,
\begin{equation*}
\lambda_L\bigl(R^{-1}(B)\cap A_N\bigr)
= \int_{1/N}^{N} r^{m-1}\,dr \cdot \sigma(B) = 0.
\end{equation*}
Hence $\lambda_L(R^{-1}(B)) = 0$ by countable subadditivity, so
$\nu(R^{-1}(B)) = 0$, i.e., $R_\#\nu(B) = 0$.
\end{proof}

\paragraph{Natural filtration of \cref{alg:GRFD}.}

The randomness drawn at iteration $k$ inside \cref{alg:construct_directions} depends
on whether $x_k$ is a spanning point. Let $\bar{\delta}_k\in\real^{2\dim\mathcal{C}}$
denote $\delta_k$ padded with trailing zeros on a branch where $\delta_k$ is drawn,
and set $\bar{\delta}_k:=0$ on the nonsmooth branch where it is not drawn. Define
\begin{equation*}
  \xi_k:=(\bar{\delta}_k,\tilde{v}_{\mathrm{ex}}^k)
  \in\real^{2\dim\mathcal{C}}\times\real^n.
\end{equation*}
The \textit{natural filtration} of \cref{alg:GRFD} is
\begin{equation*}
\mathcal{F}_k := \sigma\bigl(x_0,\,\xi_0, \dots, \xi_{k-1}\bigr),
\qquad
\mathcal{F}_0 = \sigma(x_0).
\end{equation*}
Then $x_k$ is $\mathcal{F}_k$-measurable. Conditionally on $\mathcal{F}_k$:
\begin{itemize}
    \item $\tilde{v}_{\mathrm{ex}}^k \sim \unif{\mathbb{S}_L}$ in both cases, independent of $\mathcal{F}_k$.
    \item On $\{x_k\notin\mathcal{C}_{ns}\}$: $\delta_k \sim \unif{\Delta_{m_k}}$,
    independent of $\tilde{v}_{\mathrm{ex}}^k$, where $m_k := |V^k|$ satisfies
    $m_k = 2\dim\mathcal{C}$ at interior points and $m_k = 2\dim\mathcal{C}-1$ at
    smooth boundary points, and $\bar{\delta}_k$ is this $\delta_k$ padded with trailing zeros.
    \item On $\{x_k\in\mathcal{C}_{ns}\}$: $\bar{\delta}_k=0$.
\end{itemize}
Throughout the proofs below, the deterministic choices in \cref{alg:construct_directions,alg:GRFD}
(the nearest point, inward normal, orthonormal completion, enumeration of $T$ and $U$,
line-search minimizer, and final tie-breaking) are made using fixed Borel selection
and tie-breaking conventions; thus the resulting implementation is a measurable
stochastic process.

Write $N:=4\dim\mathcal{C}-1$ for the maximum number of candidates returned by
\cref{alg:construct_directions}, including exploration. For each $x\in\mathcal{C}$,
enumerate the $T$-updates, then the $U$-updates, then the exploration direction, and
pad any remaining slots as inactive, so that each candidate occupies a fixed index
$i\in[N]$. Let $I_i(x)\in\{0,1\}$ be the Borel activity indicator of slot $i$.
For an active non-exploration slot let $q_i(x,z)$ be the vector inside $\mathrm{normalize}$
in \cref{alg:construct_directions} at state $x$ with noise coordinate $z$; for the
exploration slot put $q_i(x,z)$ equal to the $\mathbb{S}_L$-valued coordinate of $z$.
Give inactive slots any fixed nonzero value in $L$. Only events intersected with
$\{I_i(x_k)=1\}$ are used below, so the dummy values have no probabilistic role.
The maps $I_i$ and $q_i$ are Borel by the selection conventions above.

\paragraph{Conditional law and conditional absolute continuity.}

\begin{defin}\label{def:conditional-law}
Let $\mathcal{X}\in\mathcal{B}(\real^d)$ for some $d\ge 1$, let $X$ be an $\mathcal{X}$-valued random
element, and let $\mathcal{G}\subseteq\mathcal{F}$ be a sub-$\sigma$-algebra. A \textit{regular conditional
distribution} of $X$ given $\mathcal{G}$ is a function $\kappa:\Omega\times\mathcal{B}(\mathcal{X})\to[0,1]$
such that
\begin{enumerate}
  \item[(a)] $\kappa(\omega,\cdot)$ is a probability measure on $\mathcal{B}(\mathcal{X})$ for \textit{every} $\omega\in\Omega$;
  \item[(b)] for every $B\in\mathcal{B}(\mathcal{X})$, $\kappa(\cdot,B)$ is $\mathcal{G}$-measurable and is a version of
  $\mathbb{P}(X\in B\mid\mathcal{G})$, i.e.\ $\kappa(\cdot,B)=\mathbb{P}(X\in B\mid\mathcal{G})$ $\mathbb{P}$-a.s.
\end{enumerate}
Such a $\kappa$ exists for every such $X$
\citep[Chapter~2, \S7, Definition~7, Theorem~5 and Corollary]{shiryaev2016probability1}. We write $\mathcal{L}(X\mid\mathcal{G})=\kappa$.
For a measure $\nu$ on $\mathcal{B}(\mathcal{X})$ we say $\mathcal{L}(X\mid\mathcal{G})\ll\nu$ a.s.\ if some regular
conditional distribution $\kappa$ of $X$ given $\mathcal{G}$ satisfies $\kappa(\omega,\cdot)\ll\nu$ for
$\mathbb{P}$-a.e.\ $\omega$.
\end{defin}

We record the conditional law of the iteration noise in these terms.
For $x\in\mathcal{C}$ let $\varpi(x,\cdot)$ be the law of $(\bar{\delta},v)$, where
$v\sim\unif{\mathbb{S}_L}$ and, independently, $\bar{\delta}$ is the zero-padded
version of $\unif{\Delta_{m(x)}}$, with $m(x)=2\dim\mathcal{C}$ on
$\ri{\mathcal{C}}$ and $m(x)=2\dim\mathcal{C}-1$ on
$\rb{\mathcal{C}}\setminus\mathcal{C}_{ns}$; on $\mathcal{C}_{ns}$ let
$\bar{\delta}=0$. These three sets are Borel
($\mathcal{C}_{ns}$ by \cref{prop:nonsmooth_boundary_borel_null}), so $\varpi$ is a
probability kernel. The conditional sampling rules above state that
$\varpi(x_k(\omega),\cdot)$ is a regular conditional distribution of $\xi_k$ given
$\mathcal{F}_k$ in the sense of \cref{def:conditional-law}, with
$\mathcal{X}=\real^{2\dim\mathcal{C}}\times\real^n$.

\medskip

\paragraph{Isometric flattening and smoothness.}

We first record that the smoothness notion of \cref{def:1} is preserved
by the isometric flattening of $\mathcal{C}$ onto a full-dimensional
convex set in $\real^{\dim\mathcal{C}}$.
\begin{proposition}\label{prop:closed_convex_sets_are_almost_everywhere_smooth}
    Let $\mathcal{C} \subseteq \real^n$ be a nonempty closed convex set with
    $\dim\mathcal{C} \ge 1$. Fix $x_0 \in \mathcal{C}$, let
    $\bar{R} \in \real^{\dim\mathcal{C} \times n}$ be a matrix whose rows form an
    orthonormal basis of the linear subspace $\aff{\mathcal{C}-x_0}$, and
    define
    \begin{equation*}
        G : \aff{\mathcal{C}} \to \real^{\dim\mathcal{C}},
        \qquad
        G(x) := \bar{R}(x - x_0),
        \qquad
        Q := G(\mathcal{C}).
    \end{equation*}
    For $y \in \partial Q$ write
    \begin{equation*}
        \mathcal{N}(Q,y)
        := \left\{
            w \in \real^{\dim\mathcal{C}}
            \;|\;
            \norm{w}=1,\
            \inner{w}{z-y} \le 0
            \ \forall z \in Q
        \right\}
    \end{equation*}
    for the set of unit outer normals to $Q$ at $y$. Then:
    \begin{enumerate}
        \item $G$ is a bijective affine isometry of $\aff{\mathcal{C}}$
        onto $\real^{\dim\mathcal{C}}$, with inverse
        $G^{-1}(y) = x_0 + \bar{R}^T y$. The set $Q$ is closed,
        convex, and of full dimension $\dim\mathcal{C}$ in $\real^{\dim\mathcal{C}}$, and
        $G(\rb{\mathcal{C}}) = \partial Q$.
        \item For every $x \in \rb{\mathcal{C}}$,
        $\mathcal{C}$ is smooth at $x$ in the sense of \cref{def:1} if and
        only if $\left| \mathcal{N}(Q, G(x)) \right| = 1$.
    \end{enumerate}
\end{proposition}
\begin{proof}
    For part~1, since the rows of $\bar{R}$ are orthonormal,
    $\bar{R}\bar{R}^T = I_{\dim\mathcal{C}}$ and
    $\bar{R}^T\bar{R} = P_L$, the orthogonal projection onto
    $L := \aff{\mathcal{C}-x_0}$. Hence for $d \in L$,
    \begin{equation*}
        \norm{\bar{R} d}^2
        = d^T \bar{R}^T\bar{R}\, d
        = d^T P_L d
        = \norm{d}^2,
    \end{equation*}
    so $\bar{R}$ restricted to $L$ is a linear isometry into $\real^{\dim\mathcal{C}}$;
    having rank $\dim\mathcal{C}$ it is onto, and therefore a bijection
    $L \to \real^{\dim\mathcal{C}}$ whose inverse is $\bar{R}^T$. It follows
    that $G$ is a bijective affine isometry of
    $\aff{\mathcal{C}} = x_0 + L$ onto $\real^{\dim\mathcal{C}}$ with the stated inverse.
    Consequently $Q$ is convex of dimension $\dim\mathcal{C}$. Moreover $Q$ is closed,
    since $\mathcal{C}$ is closed in $\aff{\mathcal{C}}$ and $G$ is a
    homeomorphism onto $\real^{\dim\mathcal{C}}$. As $\dim Q = \dim\mathcal{C}$ we have
    $\rb{Q} = \partial Q$. Since $G$ is a homeomorphism of
    $\aff{\mathcal{C}}$ onto $\real^{\dim\mathcal{C}}$, it maps $\ri{\mathcal{C}}$ onto
    $\ri{Q}$ and $\rb{\mathcal{C}}$ onto $\rb{Q} = \partial Q$.

    For part~2, fix $x \in \rb{\mathcal{C}}$ and set $y := G(x) \in \partial Q$.
    By \cref{lem:foundational properties},
    $\aff{\mathcal{C}-x} = \aff{\mathcal{C}-x_0} = L$. We first relate
    the two feasible-direction cones. By
    \cref{lem:feasible_directions_contained_in_affine_hull},
    $D_{\mathcal{C}}(x) \subseteq L$, so for $w \in \real^n$,
    \begin{equation*}
        w \in D_{\mathcal{C}}(x)
        \iff
        w \in L
        \text{ and }
        \exists\,\lambda>0 :\, x + \lambda w \in \mathcal{C}
        \iff
        w \in L
        \text{ and }
        \exists\,\lambda>0 :\, y + \lambda \bar{R} w \in Q,
    \end{equation*}
    where the last equivalence uses
    $G(x+\lambda w) = y + \lambda \bar{R} w$. Hence
    $D_Q(y) = \bar{R}\bigl(D_{\mathcal{C}}(x)\bigr)$, and since
    $\bar{R}$ restricts to a linear homeomorphism $L \to \real^{\dim\mathcal{C}}$, the
    same holds after taking closures:
    \begin{equation}\label{eq:flattening_feasible_cone}
        \cl{D_{\mathcal{C}}(x)}
        = \left\{
            d \in L
            \;|\;
            \bar{R} d \in \cl{D_Q(y)}
        \right\}.
    \end{equation}
    Throughout we use that for a closed convex set, the closure of the
    feasible-directions cone is the tangent cone, which is the polar of
    the normal cone
    (\cite[Theorem 6.9]{RW98};
    \cite[Chapter 3, Corollary 5.2.5]{hiriart1996convex}):
    \begin{equation}\label{eq:tangent_normal_polar}
        \cl{D_Q(y)} = \mathcal{T}_Q(y) = N_Q(y)^\circ,
        \qquad
        N_Q(y) = \mathcal{T}_Q(y)^\circ.
    \end{equation}

    ($\Leftarrow$)
    Suppose $\mathcal{N}(Q,y) = \{w\}$. Since
    $N_Q(y) = \real_+ \cdot \mathcal{N}(Q,y)$, the normal cone is the
    ray $\real_+ w$, and by \eqref{eq:tangent_normal_polar},
    \begin{equation*}
        \cl{D_Q(y)}
        = (\real_+ w)^\circ
        = \left\{ g \in \real^{\dim\mathcal{C}} \;|\; \inner{g}{w} \le 0 \right\}
        = \left\{ g \in \real^{\dim\mathcal{C}} \;|\; \inner{g}{u} \ge 0 \right\},
        \qquad
        u := -w \ne 0,
    \end{equation*}
    where $u$ is the \textit{inward} unit normal. Substituting into
    \eqref{eq:flattening_feasible_cone} and using
    $\inner{\bar{R} d}{u} = \inner{d}{\bar{R}^T u}$,
    \begin{equation*}
        \cl{D_{\mathcal{C}}(x)}
        = \left\{ d \in \real^n \;|\; \inner{d}{v} \ge 0 \right\}
          \cap \aff{\mathcal{C}-x_0},
        \qquad
        v := \bar{R}^T u.
    \end{equation*}
    Here $v$ lies in the column space of $\bar{R}^T$, which equals $L$, and
    $v \ne 0$ because $\bar{R}\bar{R}^T = I_{\dim\mathcal{C}}$ makes
    $\bar{R}^T$ injective. This is precisely the half-space
    condition of \cref{def:1}, so $\mathcal{C}$ is smooth at $x$.

    ($\Rightarrow$)
    Suppose $\mathcal{C}$ is smooth at $x$, so
    $\cl{D_{\mathcal{C}}(x)}
    = \left\{d \;|\; \inner{d}{v} \ge 0\right\} \cap L$
    for some $0 \ne v \in L$. For $d \in L$ we have
    $d = \bar{R}^T\bar{R} d$, hence
    $\inner{d}{v} = \inner{\bar{R} d}{\bar{R} v}$. Setting
    $e := \bar{R} v$, which is nonzero since $\bar{R}$ is injective on
    $L$, \eqref{eq:flattening_feasible_cone} gives
    \begin{equation*}
        \cl{D_Q(y)}
        = \bar{R}\bigl(\cl{D_{\mathcal{C}}(x)}\bigr)
        = \left\{ g \in \real^{\dim\mathcal{C}} \;|\; \inner{g}{e} \ge 0 \right\}.
    \end{equation*}
    Taking polars and using \eqref{eq:tangent_normal_polar},
    $N_Q(y) = {\cl{D_Q(y)}}^\circ = \real_+(-e)$, a ray. Its unit vectors
    form the singleton
    $\mathcal{N}(Q,y) = \bigl\{ -e/\norm{e} \bigr\}$.
\end{proof}

\medskip

\paragraph{Nonsmooth boundary points.}

\begin{proposition}\label{prop:nonsmooth_boundary_borel_null}
    Let $\mathcal{C} \subseteq \real^n$ be a nonempty closed convex set with
    $\dim\mathcal{C}\ge 1$, and let
    $\mathcal{C}_{ns} \subseteq \rb{\mathcal{C}}$ be its set of nonsmooth
    boundary points in the sense of \cref{def:1}. Then $\mathcal{C}_{ns}$
    is a Borel subset of $\real^n$ and
    \begin{equation*}
        \mathcal{H}^{\dim\mathcal{C}-1}(\mathcal{C}_{ns}) = 0.
    \end{equation*}
\end{proposition}
\begin{proof}
    Let $G$, $Q$, and $\mathcal{N}(Q,\cdot)$ be as in
    \cref{prop:closed_convex_sets_are_almost_everywhere_smooth}, and let
    \begin{equation*}
        \Sigma_Q
        := \left\{
            y \in \partial Q
            \;|\;
            \left| \mathcal{N}(Q,y) \right| \ge 2
        \right\}
    \end{equation*}
    denote the set of singular boundary points of $Q$. Every
    $y \in \partial Q$ admits a supporting hyperplane, so
    $\mathcal{N}(Q,y) \ne \emptyset$. Combined with the second part of
    \cref{prop:closed_convex_sets_are_almost_everywhere_smooth}, a
    point $x \in \rb{\mathcal{C}}$ is nonsmooth if and only if
    $\left| \mathcal{N}(Q, G(x)) \right| \ge 2$, and therefore
    $\mathcal{C}_{ns} = G^{-1}(\Sigma_Q)$.

    We first show that $\Sigma_Q$ is Borel. For $j \in \mathbb{N}$ set
    \begin{equation*}
        \Sigma_j
        := \left\{
            y \in \partial Q
            \;|\;
            \exists\, w_1,w_2 \in \mathcal{N}(Q,y)
            \text{ with }
            \norm{w_1-w_2} \ge 1/j
        \right\},
    \end{equation*}
    so that $\Sigma_Q = \bigcup_{j \ge 1} \Sigma_j$. We claim that each
    $\Sigma_j$ is closed. Let $y_i \in \Sigma_j$ with $y_i \to y$; then
    $y \in \partial Q$, since $\partial Q$ is closed. Choose
    $w_1^i, w_2^i \in \mathcal{N}(Q,y_i)$ with
    $\norm{w_1^i - w_2^i} \ge 1/j$. By compactness of the unit sphere,
    pass to a subsequence along which $w_1^i \to w_1$ and
    $w_2^i \to w_2$; both limits are unit vectors and
    $\norm{w_1-w_2} \ge 1/j$. For each fixed $z \in Q$ and each
    $\ell \in \{1,2\}$ we have $\inner{w_\ell^i}{z-y_i} \le 0$, and
    letting $i \to \infty$ gives $\inner{w_\ell}{z-y} \le 0$. Hence
    $w_1,w_2 \in \mathcal{N}(Q,y)$ and $y \in \Sigma_j$. Thus
    $\Sigma_Q$ is $F_\sigma$, and therefore Borel.

    We next show that $\mathcal{H}^{\dim\mathcal{C}-1}(\Sigma_Q)=0$. By the first part
    of \cref{prop:closed_convex_sets_are_almost_everywhere_smooth}, $Q$
    is a closed convex subset of $\real^{\dim\mathcal{C}}$ of dimension $\dim\mathcal{C}$. If $\dim\mathcal{C}=1$,
    then $Q$ is a closed interval, a ray, or a line, and each of its at
    most two boundary points has a unique unit outer normal, so
    $\Sigma_Q = \emptyset$ and the claim is immediate. Assume therefore
    that $\dim\mathcal{C} \ge 2$.

    We identify $\Sigma_Q$ with the first singular stratum of $Q$ in
    \citep[Definition~1.7]{alberti1994structure}. That definition
    applies to any closed convex set $C \subseteq \real^h$ of dimension
    $h$: $\mathcal{N}(C,x)$ is the set of unit outer normals at
    $x \in \partial C$, which agrees with $\mathcal{N}(Q,y)$ in
    \cref{prop:closed_convex_sets_are_almost_everywhere_smooth}, and
    \begin{equation*}
        \Sigma^\ell(C)
        :=
        \bigl\{
            x \in \partial C
            \;|\;
            \dim\mathcal{N}(C,x) \ge \ell
        \bigr\}
    \end{equation*}
    for integers $0 \le \ell < h$. The dimension of
    $\mathcal{N}(C,x)$ is that of a subset of the unit sphere:
    \citep[Definition~1.7]{alberti1994structure} records that
    $\mathcal{N}(C,x)$ is the intersection of the unit sphere with the
    normal cone $N_C(x)$ of $C$ at $x$, and that
    $\dim\mathcal{N}(C,x) = \dim N_C(x) - 1$.

    Thus $\Sigma_Q = \Sigma^1(Q)$. If
    $\dim\mathcal{N}(Q,y) \ge 1$, then
    $\bigl|\mathcal{N}(Q,y)\bigr| \ge 2$, so $y \in \Sigma_Q$.
    Conversely, if $\bigl|\mathcal{N}(Q,y)\bigr| \ge 2$, take distinct
    $w_1,w_2 \in \mathcal{N}(Q,y)$. Convexity of $Q$ makes $N_Q(y)$ a
    convex cone, so $\pos{\{w_1,w_2\}} \subseteq N_Q(y)$. If
    $w_2 \ne -w_1$, then $0 \notin \conv\{w_1,w_2\}$ and this hull is
    one-dimensional; its radial projection onto the unit sphere remains
    one-dimensional, and the cone it generates is two-dimensional.
    The antipodal case would give $\inner{w_1}{z-y} = 0$ for all
    $z \in Q$, placing $Q$ in a hyperplane and contradicting
    $\dim Q = \dim\mathcal{C}$. Hence $\dim N_Q(y) \ge 2$, so
    $\dim\mathcal{N}(Q,y) \ge 1$ and $y \in \Sigma^1(Q)$.

    By \citep[Theorem~3]{alberti1994structure}, which is stated for
    closed convex sets of full dimension, applied to
    $Q \subseteq \real^{\dim\mathcal{C}}$ with Alberti's ambient dimension $h$ equal to
    the present $\dim\mathcal{C}$ and singularity index $1$ (admissible since
    $0 \le 1 < \dim\mathcal{C}$), the set $\Sigma^1(Q)$ is an
    $(\mathcal{H}^{\dim\mathcal{C}-2}, \dim\mathcal{C}-2)$-rectifiable subset of $\real^{\dim\mathcal{C}}$ of class
    $C^2$. By \citep[Definition~1.1]{alberti1994structure}, such a set
    is a countable union of sets of finite $\mathcal{H}^{\dim\mathcal{C}-2}$ measure
    and therefore has Hausdorff dimension at most $\dim\mathcal{C}-2$. Since
    $\dim\mathcal{C}-1 > \dim\mathcal{C}-2$, \cref{def:hausdorff-dimension} gives
    $\mathcal{H}^{\dim\mathcal{C}-1}(\Sigma_Q) = 0$.

    It remains to transfer both conclusions to $\mathcal{C}$. By the
    first part of
    \cref{prop:closed_convex_sets_are_almost_everywhere_smooth},
    $G^{-1}$ is an isometry, and in particular Lipschitz with constant
    $1$. By \cref{prop:lipschitz-hausdorff},
    \begin{equation*}
        \mathcal{H}^{\dim\mathcal{C}-1}(\mathcal{C}_{ns})
        = \mathcal{H}^{\dim\mathcal{C}-1}\bigl(G^{-1}(\Sigma_Q)\bigr)
        \le 1^{\dim\mathcal{C}-1} \cdot \mathcal{H}^{\dim\mathcal{C}-1}(\Sigma_Q)
        = 0.
    \end{equation*}
    Finally, $\mathcal{C}_{ns}=G^{-1}(\Sigma_Q)$ is the preimage of the $F_\sigma$ set $\Sigma_Q$ under the
    continuous map $G:\aff{\mathcal{C}}\to\real^{\dim\mathcal{C}}$, hence $F_\sigma$ in
    $\aff{\mathcal{C}}$. Since $\aff{\mathcal{C}}$ is closed in $\real^n$, $\mathcal{C}_{ns}$ is an
    $F_\sigma$ (in particular Borel) subset of $\real^n$, and therefore a countable union of compact
    sets, since every closed subset of $\real^n$ is one.
\end{proof}

\begin{remark}\label{rem:nonsmooth_dim_one}
    When $\dim\mathcal{C}=1$ the conclusion reads $\mathcal{H}^0(\mathcal{C}_{ns})=0$,
    i.e.\ $\mathcal{C}_{ns}=\emptyset$; this case is treated separately in the proofs of
    \cref{prop:nonsmooth_boundary_borel_null} and \cref{lem:direction_density}.
\end{remark}

\begin{lem}[Radial image of a null set]\label{lem:radial_image_null}
    Let $L\subseteq\real^n$ be a linear subspace, let
    $x\in\real^n$, and write $\mathbb{S}_L:=\{\omega\in L:\norm{\omega}=1\}$ and
    $\Pi_x(y):=(y-x)/\norm{y-x}$ for $y\in (x+L)\setminus\{x\}$. If
    $A\subseteq (x+L)\setminus\{x\}$ satisfies $\mathcal{H}^{\dim L-1}(A)=0$, then
    $\mathcal{H}^{\dim L-1}\!\left(\Pi_x(A)\right)=0$ in $\mathbb{S}_L$.
\end{lem}
\begin{proof}
    Write $r(z):=z/\norm{z}$ for $z\in L\setminus\{0\}$, so that
    $\Pi_x=r\circ\tau$ with $\tau(y):=y-x$ an isometry. For any
    $z_1,z_2\in L\setminus\{0\}$,
    \begin{equation*}
        z_1\norm{z_2}-z_2\norm{z_1}
        = \norm{z_1}(z_1-z_2)+z_1\bigl(\norm{z_2}-\norm{z_1}\bigr),
    \end{equation*}
    so by the triangle and reverse triangle inequalities
    $\norm{z_1\norm{z_2}-z_2\norm{z_1}}\le 2\norm{z_1}\,\norm{z_1-z_2}$.
    Dividing by $\norm{z_1}\norm{z_2}$ gives $\norm{r(z_1)-r(z_2)}\le
    2\norm{z_1-z_2}/\norm{z_2}$. By symmetry (swapping $z_1$ and $z_2$)
    we also have $\norm{r(z_1)-r(z_2)}\le
    2\norm{z_1-z_2}/\norm{z_1}$. Taking the minimum of the two bounds,
    \begin{equation*}
        \norm{r(z_1)-r(z_2)}
        \;\le\; \frac{2\,\norm{z_1-z_2}}{\max(\norm{z_1},\norm{z_2})}.
    \end{equation*}
    In particular $r$, and hence $\Pi_x$, is Lipschitz with constant $2m$ on
    $\{\,\norm{z}\ge 1/m\,\}$. For $m\in\mathbb{N}$ set
    $A_m:=A\cap\{\,y:\norm{y-x}\ge 1/m\,\}$. A Lipschitz map does not increase
    Hausdorff dimension and scales $\mathcal{H}^{\dim L-1}$ by at most the
    $(\dim L-1)$-th power of its constant, so
    $\mathcal{H}^{\dim L-1}(\Pi_x(A_m))\le (2m)^{\,\dim L-1}\,\mathcal{H}^{\dim L-1}(A_m)=0$.
    Since $\Pi_x(A)=\bigcup_{m\ge1}\Pi_x(A_m)$ is a countable union of null
    sets, it is null.
\end{proof}

\begin{lem}[Conditional density of the step direction]\label{lem:direction_density}
    Fix an iteration number $k$, suppose $\alpha\in(0,1)$ and $\eta\in(0,1]$, and condition on
    $\mathcal{F}_k$ on the event $x_k\notin\mathcal{C}_{ns}$. Set
    $L:=\aff{\mathcal{C}-x_k}$. For every $i\in[s_k]$, let $q_i$ be the vector inside
    $\mathrm{normalize}$ that produces $\tilde{v}_i^k$, with
    $q_i:=\tilde{v}_{\mathrm{ex}}^k$ for the exploration slot. Then $q_i\neq0$ for every
    realization of the iteration noise, $\tilde{v}_i^k=q_i/\norm{q_i}$ lies in
    $\mathbb{S}_L$, and the conditional law of
    $\omega_i:=\tilde{v}_i^k$ given $\mathcal{F}_k$ is
    absolutely continuous with respect to the surface measure
    $\sigma=\mathcal{H}^{\dim\mathcal{C}-1}\!\restriction\mathbb{S}_L$.
\end{lem}
\begin{proof}
    The index $i$ ranges over the active slots of
    the enumeration fixed above, so $[s_k]=\{i\in[N] : I_i(x_k)=1\}$.

    We first verify that every active generator is nonzero pointwise. Write
    $d:=\dim\mathcal{C}$ and recall that
    $\bar{\eta}=\eta/(2\sqrt{d})\le 1/2$. Since
    $\zeta_k=\sum_{j=1}^{m_k}\delta_{k,j}v_{k,j}$, where
    $\delta_k\in\Delta_{m_k}$ and every $v_{k,j}\in U\cup T$ has unit norm,
    \begin{equation*}
      \norm{\zeta_k}
      \le \sum_{j=1}^{m_k}\delta_{k,j}\norm{v_{k,j}}
      =1.
    \end{equation*}
    For an active $U$-update, the reverse triangle inequality gives
    \begin{equation*}
      \norm{q_i}
      =\norm{\alpha\bar{\eta}\,\zeta_k+(1-\alpha\bar{\eta})v}
      \ge 1-2\alpha\bar{\eta}
      =1-\frac{\alpha\eta}{\sqrt{d}}
      \ge 1-\alpha>0.
    \end{equation*}
    For an active $T$-update, $u$ and $v$ are orthogonal unit vectors, and hence
    \begin{align*}
      \norm{q_i}
      &\ge
      \norm{(1-\alpha)\bar{\eta}\,u+(1-\bar{\eta})v}
      -\alpha\bar{\eta}\norm{\zeta_k} \\
      &\ge 1-(1+\alpha)\bar{\eta}
      \ge \frac{1-\alpha}{2}>0.
    \end{align*}
    For the exploration slot, $\norm{q_i}=\norm{\tilde{v}_{\mathrm{ex}}^k}=1$.
    Thus $q_i\neq0$ for every admissible realization of the iteration noise, and
    $\tilde{v}_i^k=q_i/\norm{q_i}\in\mathbb{S}_L$.

    For the kernel construction, let $R(q):=q/\norm{q}$ for $q\neq 0$ and assign
    an arbitrary fixed value in $\mathbb{S}_L$ at $q=0$. This extension only
    specifies $R(q_i(x,z))$ at off-support noise values and does not affect an
    active algorithmic direction. For each slot $i\in[N]$ define
    \begin{equation*}
      \kappa_i(\omega,B)
      :=\varpi\bigl(x_k(\omega),
          \{z:R(q_i(x_k(\omega),z))\in B\}\bigr).
    \end{equation*}
    The map $z\mapsto R(q_i(x_k,z))$ is Borel, so $\kappa_i$ is a probability kernel as the
    image of $\varpi(x_k,\cdot)$ under that map. It is a regular conditional distribution of
    $R(q_i(x_k,\xi_k))$ given $\mathcal{F}_k$: for bounded Borel $f$,
    \begin{equation*}
      \mathbb{E}\!\left[f\!\left(R(q_i(x_k,\xi_k))\right)\mid\mathcal{F}_k\right]
      =\int f\!\left(R(q_i(x_k,z))\right)\varpi(x_k,dz)
      \qquad\text{a.s.}
    \end{equation*}
    The identity holds first for indicator functions of rectangles in the joint
    variables $(z,x_k)$ by the definition of the conditional kernel and
    \citep[Chapter~2, \S7, property~K$^{*}$]{shiryaev2016probability1}, then for bounded
    Borel functions by the monotone-class theorem.

    If $i$ is an active exploration slot, then by
    \cref{alg:construct_directions} $\omega_i=\tilde{v}_{\mathrm{ex}}^k$ is uniform on
    $\mathbb{S}_L$, whose law is $\sigma/\sigma(\mathbb{S}_L)\ll\sigma$, and therefore
    $\kappa_i(\omega,\cdot)\ll\sigma$ on this event.

    If $i$ is an active non-exploration slot, freeze $x_k=x$ and write the pre-normalized
    generator as
    \begin{equation*}
      q_i=b_i+\alpha\bar{\eta}\,\zeta_k,
    \end{equation*}
    where $b_i$ is determined by $x_k$ via the selection conventions above, hence
    $\mathcal{F}_k$-measurable, and $\alpha\bar{\eta}>0$. We first show
    $\zeta_k$ has a conditional density on $L$, i.e.\ $\zeta_k\ll\lambda_L$.
    Write $m_k:=\lvert V^k\rvert$ for the number of base spanning directions. Recall
    $\zeta_k=\sum_{j=1}^{m_k}\delta_{k,j}\,v_{k,j}$ with
    $V^k=\{v_{k,1},\dots,v_{k,m_k}\}\subseteq L$.

    \noindent Since $x_k\notin\mathcal{C}_{ns}$, the
    point $x_k$ is interior or a smooth boundary point. If $x_k$ is a smooth
    boundary point with $\dim\mathcal{C}=1$, then $V^k=\{u\}$, so
    $q_i=u$ is deterministic and nonzero, $\mathbb{S}_L=\{u,-u\}$,
    and $\sigma=\mathcal{H}^{0}\!\restriction\mathbb{S}_L$ is counting measure;
    the law of $\omega_i=u$ is then trivially absolutely continuous with
    respect to $\sigma$, hence $\kappa_i(\omega,\cdot)\ll\sigma$ on this event, and the claim holds. In every other case the
    difference vectors $\{v_{k,j}-v_{k,1}\}_{j}$ span $L$.
    Since $x_k\notin\mathcal{C}_{ns}$, \cref{alg:construct_directions} takes the first branch ($\rho(x_k)>\rho_0$), or the third branch at an interior point ($0<\rho(x_k)\le\rho_0$), or the third branch at a smooth boundary point ($\rho(x_k)=0$). In the first case, $V^k=\{\pm w_1,\dots,\pm w_{\dim\mathcal{C}}\}$.
    In the second, $V^k=U\cup T$ as a set equals
    $U=\{\pm u,\pm w_1,\dots,\pm w_{\dim\mathcal{C}-1}\}$.
    In the third, with $\dim\mathcal{C}\ge 2$,
    $V^k=\{u,\pm w_1,\dots,\pm w_{\dim\mathcal{C}-1}\}$.
    The first two contain an antipodal pair in every basis
    direction of $L$, and the third contains the antipodal pair $\pm w_1$
    together with $u$.

    \noindent Hence the linear map $T:\delta\mapsto\sum_j\delta_j v_{k,j}$, restricted to the affine hull of
    $\Delta_{m_k}$, has differential of rank $\dim\mathcal{C}$ onto $L$. As
    $\delta_k$ has a density on $\Delta_{m_k}$ (\cref{alg:construct_directions}), its
    pushforward $\zeta_k=T(\delta_k)$ is absolutely continuous with respect
    to $\lambda_L$ by \cref{prop:affine-image-abs-cont}.

    Consequently $q_i=b_i+\alpha\bar{\eta}\zeta_k$, with $b_i\in L$,
    is a bijective affine image of $\zeta_k$ within $L$, hence
    $q_i\ll\lambda_L$ by the bijective case of \cref{prop:affine-image-abs-cont}.
    By the pointwise nonvanishing established above,
    $\omega_i=\tilde{v}_i^k=R(q_i)$ surely.
    Since $q_i\ll\lambda_L$, \cref{lem:radial-projection-ac} gives
    $R(q_i)\ll\sigma$. As the frozen state
    $x\in\mathcal{C}\setminus\mathcal{C}_{ns}$ was arbitrary,
    $\kappa_i(\omega,\cdot)\ll\sigma$ whenever
    $I_i(x_k(\omega))=1$ and $x_k(\omega)\notin\mathcal{C}_{ns}$.
\end{proof}

\begin{proof}[Proof of \cref{lem:always_smooth_almost_surely}]
    Let $\{\mathcal{F}_k\}_{k=0}^{K}$ be the natural filtration of the
    algorithm, generated by $x_0$ together with the randomness of iterations
    $0,\dots,k-1$; thus $x_k$ is $\mathcal{F}_k$-measurable. Write
    $L:=\aff{\mathcal{C}-x_k}$ (independent of $x_k$). We argue by induction on $k$; the base case
    holds by assumption.

    Assume $\mathbb{P}(x_j\in\mathcal{C}_{ns})=0$ for all $j\le k$. We work
    conditionally on $\mathcal{F}_k$, on the event $\{x_k\notin\mathcal{C}_{ns}\}$,
    which has full probability by the inductive assumption. On this event
    $x_k$ is interior or a smooth boundary point, so every direction in
    $\tilde{V}^k$ is either produced by the $T$- and $U$-updates of
    \cref{alg:construct_directions} or is the exploration
    direction, and \cref{lem:direction_density} applies to each.

    Fix $i\in[N]$ and let $A_i:=\{I_i(x_k)=1\}$, which is $\mathcal{F}_k$-measurable.
    On $A_i$ set $Y_i:=x_k+t_i^k\tilde{v}_i^k$ and
    $\omega_i:=\tilde{v}_i^k$; define $Y_i$ arbitrarily off $A_i$. We claim
    \begin{equation}\label{eq:endpoint_implies_radial}
        A_i\cap\{Y_i\in\mathcal{C}_{ns}\}
        \;\subseteq\;
        A_i\cap\left\{\omega_i\in
          \Pi_{x_k}(\mathcal{C}_{ns}\setminus\{x_k\})\right\}.
    \end{equation}
    Indeed, on $A_i\cap\{Y_i\in\mathcal{C}_{ns}\}$ we have
    $Y_i\in\mathcal{C}_{ns}\subseteq\rb{\mathcal{C}}$; since
    $x_k\notin\mathcal{C}_{ns}$, necessarily $Y_i\neq x_k$, whence
    $t_i^k\norm{\tilde{v}_i^k}=\norm{Y_i-x_k}>0$ and
    $\Pi_{x_k}(Y_i)=(Y_i-x_k)/\norm{Y_i-x_k}=\omega_i$ with
    $Y_i\in\mathcal{C}_{ns}\setminus\{x_k\}$. This proves
    \eqref{eq:endpoint_implies_radial}. The inclusion \eqref{eq:endpoint_implies_radial} holds for every value of the step $t_i^k$.

    For $y\in\mathcal{C}$ let $D_y:=\Pi_y(\mathcal{C}_{ns}\setminus\{y\})\subseteq\mathbb{S}_L$ and
    \begin{equation*}
      D:=\{(\theta,y)\in\mathbb{S}_L\times\mathcal{C}:\ \theta\in D_y\}
       =\{(\theta,y)\in\mathbb{S}_L\times\mathcal{C}:\ y+r\theta\in\mathcal{C}_{ns}\text{ for some }r>0\}.
    \end{equation*}
    By construction $D_{x_k}=\Pi_{x_k}(\mathcal{C}_{ns}\setminus\{x_k\})$, so the right-hand
    event of \eqref{eq:endpoint_implies_radial} is $A_i\cap\{(\omega_i,x_k)\in D\}$.
    \textit{(i) Every section is $\sigma$-null.} By \cref{prop:nonsmooth_boundary_borel_null}, $\mathcal{H}^{\dim\mathcal{C}-1}(\mathcal{C}_{ns})=0$,
    so \cref{lem:radial_image_null} (applied to the subspace $L$ of dimension $\dim\mathcal{C}$, the point $y$, and
    $A=\mathcal{C}_{ns}\setminus\{y\}\subseteq(y+L)\setminus\{y\}$) gives $\sigma(D_y)=0$ for every $y\in\mathcal{C}$.

    \textit{(ii) $D$ is Borel.} By the last paragraph of the proof of \cref{prop:nonsmooth_boundary_borel_null},
    $\mathcal{C}_{ns}=\bigcup_{m\ge1}K_m$ with $K_m$ compact. For $m,j\in\mathbb{N}$ the set
    $\{(\theta,y,r)\in\mathbb{S}_L\times(\mathcal{C}\cap\mathcal{B}[0,j])\times[1/j,j]:\ y+r\theta\in K_m\}$ is compact, hence so is its
    projection $D_{m,j}$ onto the $(\theta,y)$-coordinates, and $D=\bigcup_{m,j}D_{m,j}$ (given $(\theta,y)\in D$ with
    $y+r\theta\in K_m$, take $j\ge\max\{\norm{y},r,1/r\}$). Thus $D$ is a countable union of compact sets, hence Borel;
    in particular every $D_y$ is Borel.

    \textit{(iii) Conclusion.} Let $\kappa_i$ be the regular conditional distribution constructed in the proof of
    \cref{lem:direction_density}. For every Borel $E\subseteq\mathbb{S}_L\times\mathcal{C}$,
    \begin{equation*}
      \mathbb{P}\big((\omega_i,x_k)\in E\mid\mathcal{F}_k\big)(\omega)
      =\kappa_i\big(\omega,E_{x_k(\omega)}\big)
      \qquad\mathbb{P}\text{-a.s.}
    \end{equation*}
    Indeed, the identity holds for rectangles by property K$^{*}$ and extends to Borel $E$ by the
    monotone-class theorem. On $A_i\cap\{x_k\notin\mathcal{C}_{ns}\}$,
    \cref{lem:direction_density} gives $\kappa_i(\omega,\cdot)\ll\sigma$, while (i) gives
    $\sigma(D_{x_k(\omega)})=0$. Since $A_i\in\mathcal{F}_k$, property K$^{*}$ and the inductive
    assumption therefore give
    \begin{equation*}
        \mathbb{P}\!\left(A_i\cap\{Y_i\in\mathcal{C}_{ns}\}\mid\mathcal{F}_k\right)
        \;\le\;
        \mathbf{1}_{A_i}\,
        \kappa_i\big(\cdot,D_{x_k}\big)
        \;=\;0
        \qquad\text{a.s.}
    \end{equation*}

    Since $x_{k+1}$ is one of the active candidates,
    $\{x_{k+1}\in\mathcal{C}_{ns}\}
    \subseteq\bigcup_{i=1}^{N}\left(A_i\cap\{Y_i\in\mathcal{C}_{ns}\}\right)$,
    and a union bound yields
    \begin{equation*}
        \mathbb{P}\!\left(x_{k+1}\in\mathcal{C}_{ns}\mid\mathcal{F}_k\right)
        \;\le\;
        \sum_{i=1}^{N}\mathbb{P}\!\left(A_i\cap\{Y_i\in\mathcal{C}_{ns}\}\mid\mathcal{F}_k\right)
        \;=\;0
        \qquad\text{a.s.}
    \end{equation*}
    Taking expectations, and using that the conditioning event
    $\{x_k\notin\mathcal{C}_{ns}\}$ has full probability by the inductive
    assumption, gives
    $\mathbb{P}(x_{k+1}\in\mathcal{C}_{ns})=0$, which completes the induction.
\end{proof}

\end{document}